\documentclass[12pt]{amsart}
\usepackage[T1]{fontenc}

\usepackage[square,numbers]{natbib}
\usepackage{amssymb}
\usepackage{hyperref}
\usepackage{graphicx}
\usepackage{tikz-cd}
\usepackage{amsmath}
\usepackage{xfrac}

\usepackage{caption}
\usepackage{subcaption}

\newtheorem{theorem}{Theorem}[section]
\newtheorem{lemma}[theorem]{Lemma}
\newtheorem{corollary}[theorem]{Corollary}

\theoremstyle{definition}

\theoremstyle{remark}
\newtheorem{remark}{Remark}

\newcommand{\R}{\mathbb{R}}
\newcommand{\C}{\mathbb{C}}
\newcommand{\Z}{\mathbb{Z}}
\newcommand{\Q}{\mathbb{Q}}

\newcommand{\Parn}{{\rm Par}_n}
\newcommand{\Parj}{{\rm Par}_J}

\newcommand{\hsp}{\hspace{0.2em}}
\newcommand{\Comp}{{\rm Comp}_n}

\title{Spherical Schubert varieties in $GL_n / P$}
\author{Jaewon Min}
\address{Dept.~of Mathematics, University of Illinois at Urbana-Champaign, Urbana, IL 61801}
\email{jaewonm2@illinois.edu}
\date{October 30, 2026}

\begin{document}

\maketitle

\vspace{-1ex}

\begin{abstract}
Let a general linear group act on a partial flag variety transitively, and consider the induced action of a Levi subgroup $L$. A Schubert variety is \emph{$L$-spherical} if it is stable under the $L$-action and contains a dense open orbit of a Borel subgroup of $L$. It is \emph{toric} if, in addition, $L$ is a torus. In this paper, we give a combinatorial classification of $L$-spherical Schubert varieties in partial flag varieties. This is achieved by applying Lakshmibai-Seshadri paths forming crystal graphs of Demazure modules. Previously, Hodges and Lakshmibai classified $L$-spherical Schubert varieties in Grassmannians. $L$-spherical Schubert varieties in complete flag varieties were treated independently by Gao--Hodges--Yong and Can--Saha. Our classification generalizes and unifies these results.
\end{abstract}

\setcounter{tocdepth}{2}


\section{Introduction}\label{sub1.1}

\subsection{History and motivation}

Our main goal is to classify Levi-spherical Schubert varieties in partial flag varieties $GL_n (\C) /P$, generalizing and unifying results of \cite{MR3044412}, \cite{MR4517644}, \cite{MR4690502} and \cite{MR4643501}. This topic sits in the broader context of sphericality in Lie theory. Notably, the space of harmonic polynomials on a sphere is $SO_3 (\R)$-multiplicity-free \cite[Corollary 17.17]{MR3112817}. E.~B.~Vinberg pointed out this is the origin of \textit{spherical homogeneous spaces} \cite[Ch.~I, \S 3]{MR1845642}. Another example is from enumerative geometry, \textit{Steiner's problem}. This concerns the number of smooth plane conics tangent to five given conics, computed by M. Chasles \cite[Ch.~0]{MR3617981}. Its parameter space is a $PGL_3(\C)$-spherical homogeneous space. In the modern solution, the key is to take a compactification, which is a \textit{spherical variety} with $PGL_3(\C)$-action; see M.~Brion \cite[\S 2.4, \S3.6]{MR1040855}.

Our proof is representation-theoretic and is inspired by the Borel-Weil theorem. From line bundles over a Schubert variety, construct \textit{Demazure modules} \cite[\S 5.4]{MR0354697}. A tensor product of a one-dimensional Demazure module with an arbitrary one has an \textit{excellent filtration}, named by P.~Polo \cite{MR1021515}. The existence of such filtrations was suggested by A.~Joseph and then studied by O.~Mathieu and many others over decades; for the history, see \cite{MR2214240}, edited by J.~Bernstein \textit{et al}. A Schubert variety is spherical precisely when all the relevant excellent filtrations are multiplicity-free. The combinatorial description of the filtration was given by V.~Lakshmibai, P.~Littelmann and P.~Magyar, applying standard monomial theory \cite{MR1887117}. A basis of a Demazure module is parametrized by \textit{Lakshmibai-Seshadri paths}, which in turn form a \textit{crystal graph}; see, \textit{e.g.}, the work of M.~Kashiwara \cite{MR1240605} and G.~Lusztig \cite{MR1227098}.

\subsection{Background}

We refer to A.~Borel's textbook \cite[\S 6]{MR1102012} for the framework. Let $G$ be a connected, complex, reductive algebraic group. For each closed subgroup $H$, $G/H$ is called a \textbf{homogeneous space}, which is a smooth quasi-projective variety. $H$ is a \textbf{parabolic subgroup} if, in addition, $G/H$ is projective \cite[\S 11]{MR1102012}. \textbf{Borel subgroups} are minimal parabolic subgroups, which are conjugate to each other in $G$. Denote a parabolic subgroup and a Borel subgroup by $P$ and $B$, respectively. $G/P$ is called a \textbf{partial flag variety}.

A normal $G$-variety $X$ is \textbf{$G$-spherical} if it has an open, and therefore dense, orbit of a Borel subgroup $B$ of $G$. In N.~Perrin's survey, $X$ is $G$-spherical if and only if it has finitely many $B$-orbits \cite[Theorem 2.1.2]{MR3177371}. Assuming that $X$ is quasi-projective, this is also equivalent to the condition that for any $G$-linearized line bundle $L$ over $X$, the space of global sections $H^0 (X , L)$ is $G$-multiplicity-free.

Automatically, a $G$-spherical variety has finitely many $G$-orbits. Here, each $G$-orbit is $G$-isomorphic to a $G$-spherical homogeneous space. For instance, $G/P$ with the left $G$-action is a $G$-spherical homogeneous space. Indeed, consider a $G$-linearized line bundle $L$ over $G/P$. From the Borel-Weil theorem, $H^0(G/P , L)$ is either irreducible or zero, and hence $G$-multiplicity-free.

Since $G/P$ is $G$-spherical, it has finitely many $B$-orbits, called \textbf{Bruhat cells} \cite[\S 14]{MR1102012}. To describe the orbits, fix a Borel subgroup and a maximal torus $T \subseteq B$. Let $(W,S)$ be the Coxeter system composed of the Weyl group $W= N_G(T)/T$ and the minimal set of generators $S = \{s_1 , \cdots , s_n\}$. Parabolic subgroups $P_J$ containing $B$ are indexed by $J \subseteq [n]$. Here, $N_{P_J}(T)/T$ is isomorphic to $W_J = \langle s_j \mid j \in J \rangle$. From the \textbf{Bruhat decomposition},
\begin{equation}
G/P_J = \bigsqcup_{wW_J \in W/W_J} BwP_J / P_J , \quad X(wW_J):= \overline{BwP_J / P_J }.
\end{equation}
Subvarieties $X(wW_J)$ of $G/P_J$ are called \textbf{Schubert varieties}.

Recall that $G/P_J$ has the left $G$-action. However, the subvariety $X(wW_J)$ is not stable under the $G$-action in general. In fact, there exists $I \subseteq [n]$ such that ${\rm stab}_G(X(wW_J)) = P_I$; see, \textit{e.g.}, \cite[\S 1.2]{MR4394321}. Equivalently, $X(wW_J)$ is stable under $P_I$ if and only if for some $v \in wW_J$, $v$ is the maximal representative of $W_I v W_J$ with respect to the Coxeter length $\ell$.

Since $P_I$ is not reductive for $I \neq [n]$, we choose a maximal reductive subgroup, called a \textbf{Levi subgroup} \cite[\S 11]{MR1102012}. Levi subgroups are conjugate to each other. Choose one containing $T$ and denote it by $L_I$. $P_I$ has a \textbf{Levi decomposition}
\begin{equation}
P_I = L_I \ltimes R_u(P_I),
\end{equation}
where $R_u(P_I)$ is the unipotent radical. Since $R_u(P_I) \subseteq B$, $X(wW_J)$ is stable under $L_I$ if and only if it is stable under $P_I$.

A Schubert variety is a normal variety \cite[Theorem 4.7]{MR0630769}. For any Schubert variety with a Levi-action, we want to determine Levi-sphericality. More specifically, let $I , J \subseteq [n]$ and $w \in W$ be the maximal representative of $W_I w W_J$. There are $| W_I \backslash W / W_J|$ distinct classification problems, generalizing the question posed in \cite[Definition 1.8]{MR4394321}. 
\begin{equation}\label{eqn1.9}
\textit{Which Schubert varieties $X(wW_J)$ are $L_I$-spherical?}
\end{equation}
We propose a framework for addressing this problem and then give an explicit combinatorial classification when $G$ is of Lie type A.

\subsection{Main results} 

Let $G = GL_n$, a complex general linear group. The Weyl group $W$ is the symmetric group $S_n$. For any $I , J \subseteq [n-1]$, the double coset decomposition of $W$ is
\begin{equation}
W = \bigsqcup_{u \in {}^I W^J} W_I u W_J,
\end{equation}
where ${}^I W^J$ is the set of minimal representatives of $W_I \backslash W / W_J$ with respect to $\ell$; see, \textit{e.g.}, \cite[\S 1]{MR2120109}.

For each $J \subseteq [n-1]$, define a sequence $b_1 , b_2 ,  \cdots , b_n$ inductively as follows. Start with $b_1 = 1$. $b_{j+1} = b_j$ if $j \in J$ and $b_{j+1} = 1+b_j$ if $j \notin J$. Let $w \in W$. Write $x_{j} := b_{w^{-1}(j)}$ and denote the sequence of $x_j$ by $[w,J] := x_1 x_2 \cdots x_n$. In addition, let $I \subseteq [n-1]$. Insert $i$ between $x_j$ and $x_{j+1}$ whenever $j \in [n-1]-I$. $[w,J]$ with $i$ insertions, which is a sequence of positive integers and symbols $i$, is denoted by $[w,I,J]$.

Let $\mathbf{a}$ and $\mathbf{b}$ be finite sequences of positive integers and symbols $i$. For instance, they can be taken as $[w,I,J]$. We say $\mathbf{a}$ has $\mathbf{b}$ as a \textbf{pattern} if there is a subsequence $\mathbf{c}$ of $\mathbf{a}$ with two conditions satisfied. First of all, placements of symbols $i$ are exactly the same for $\mathbf{b}$ and $\mathbf{c}$. Secondly, the remaining numerical entries of $\mathbf{b}$ and $\mathbf{c}$ have the same relative order, including equalities. In other words, if the remaining numerical entries are $b_1, \cdots , b_p$ and $c_1 , \cdots , c_q$, then $p= q$ and for any $1 \leq r < s \leq p$, only three cases are possible: (1) $b_r - b_s , c_r - c_s >0$, (2) $b_r - b_s , c_r - c_s =0$ and (3) $b_r - b_s , c_r - c_s <0$. In particular, if $I = [n-1]$ and $J = \emptyset$, then $[w^{-1},I,J]$ is the one-line notation of a permutation $w$. Therefore, $w$ has $v$ as the usual \textbf{permutation pattern} if $[w^{-1} , I , J]$ has $[v^{-1} , I , J]$ as a pattern.

For instance, if  $n=7$ and $J=\{ 2,4,6\}$, then $b_1b_2 b_3 b_4 b_5 b_6 b_7 = 1223344$. If in one-line notation $w = 1452673$, \textit{i.e.}, $w^{-1} = 1472356$, then $[w,J] = 1342234$. If $I= \{ 1,4,5\}$, then $[w,I,J] = 13i4i223i4$. $\underline{1}\underline{3}\underline{i}\underline{4}i22\underline{3}\underline{i}\underline{4}$ has $12i32i3$ as a pattern due to the underlined subsequence $13i43i4$. It avoids a pattern $231i1$ since there is no subsequence $bcaia$ where the numerical entries satisfy $a < b< c$. Our first theorem, which classifies spherical Schubert varieties in partial flag varieties $GL_n/P_J$, is stated below.

\begin{theorem}\label{thm1.2}
Let $G= GL_n$, $W= S_n$ and $I, J \subseteq [n-1]$. Let $u \in {}^I W^J$ and $w$ be the maximal representative of $W_I u W_J$. Then $X(wW_J)$ is $L_I$-spherical if and only if $[u,I,J]$ avoids the following patterns.
\begin{equation}\label{eqn1.7}
\begin{matrix}
321, & 3412, & 3i312, & 231i1, & 2i21i1, & 244113,\\
233112, & 3331122, & 2233111, & 22211i11, & 22i22111, & 22ii1i1, \\
33ii12, & 44123, & 3312i2, & 2i2ii11, & 23ii11, & 23411, \\
2i2311, & 221211, & 331322, & 221311, & 331422.
\end{matrix}
\end{equation}
\end{theorem}

Previously, the case $I= \emptyset$ was classified by P.~Karuppuchamy in \cite{MR3044412}. The case $|J| = n-2$ was treated by R.~Hodges and V.~Lakshmibai in \cite{MR4517644}. The case $J=\emptyset$ was completed by Y.~Gao, R.~Hodges and A.~Yong in \cite{MR4615612} and by M.~B.~Can and P.~Saha in \cite{MR4643501}, separately; see \S \ref{sub1.2} for more details. As a corollary, we conclude sphericality is symmetric between index sets $I$ and $J$.

\begin{corollary}\label{coro1.3}
Let $G= GL_n$, $W= S_n$ and $I,J \subseteq [n-1]$. Assume that $w \in W$ is the maximal representative of $W_I w W_J$. Then
\begin{equation}
\textit{$X(wW_J)$ is $L_I$-spherical} \quad \Leftrightarrow \quad \textit{$X(w^{-1}W_I)$ is $L_J$-spherical}.
\end{equation}
\end{corollary}

We prove the theorem and the corollary in \S \ref{sub6.2}. In particular, if $J = \emptyset$, then $P_J = B$ and $L_J  = T$. $G/B$ is called a \textbf{complete flag variety}. A normal variety with a $T$-action is \textbf{toric} if there exists a dense open orbit of $T$. Therefore, when $w$ is the maximal representative of $W_Iw$,
\begin{equation}
\textit{$X(w)$ is $L_I$-spherical} \quad \Leftrightarrow \quad \textit{$X(w^{-1}W_I)$ is toric}.
\end{equation}
In other words, classifications of toric Schubert varieties in partial flag varieties and spherical Schubert varieties in complete flag varieties imply each other.

The classification of toric $X(wW_J) \subseteq G/P_J$ follows from the special case $X(w) \subseteq G/B$ due to $T$-equivariant birational equivalence. As \cite[Section 4.3]{gao2024orbitstructurescomplexityschubert} pointed out, however, the same argument cannot be applied to determine sphericality. For instance, consider $G = GL_4$, $W = S_4$, $I = \{1,3\}$, $J = \{2\}$ and $w = 4213$. Then the minimal representative of $wW_J$ is $v = 4123$ and $X(v)$ is birationally equivalent to $X(wW_J)$; see, \textit{e.g.}, \cite[Section 1.2]{MR2143072}. Unlike $X(wW_J)$, $X(v)$ is not stable under the $L_I$-action since $v$ is not the maximal representative of $W_I v$. Therefore, we cannot deduce sphericality of Schubert varieties in $G/P_J$ from that in $G/B$.

Fix a permutation $u \in W$. The \textbf{left ascent set} (resp. \textbf{right ascent set}) of $u$ is the set of indices $j \in [n-1]$ such that $\ell (s_j u) > \ell (u)$ (resp. $\ell (u s_j ) > \ell (u)$). We have another corollary, whose proof appears in \S \ref{sub6.2}.

\begin{corollary}\label{coro1.4}
Let $G = GL_n$, $W = S_n$ and $u \in W$. If $u$ has one of the permutation patterns
\begin{equation}\label{eqn1.12}
321, \quad 246135, \quad 351624, \quad 356124, \quad 415263, \quad 451623, \quad 2561374, \quad 4157236,
\end{equation}
then for any $I, J \subseteq [n-1]$ satisfying $u \in {}^I W^J$, $X(wW_J)$ is not $L_I$-spherical. Conversely, if $u$ avoids all the patterns, then by taking $I$ and $J$ to be the left and the right ascent set of $u$, $X(wW_J)$ is $L_I$-spherical. Here, $w$ is always chosen as the maximal representative of $W_I u W_J$.
\end{corollary}

A similar result can be found in \cite[Theorem 1.4]{MR4376215}. In that paper, C.~Gaetz proved that for each $w \in S_n$ with the \textit{left descent set} $I$, $X(w)$ is $L_I$-spherical if and only if $w$ avoids twenty-one permutation patterns. In this case, we call $w$ \textbf{maximally spherical}; for more details, see \cite[\S 3]{MR4394321}.

On the way to the proof of Theorem \ref{thm1.2}, we give another classification when $X(w)$ is in $GL_n/B$. Write
\begin{equation}\label{eqn1.3}
\Parn := \{ \lambda = (\lambda_1 , \cdots , \lambda_n) \mid \lambda_1 \geq \cdots \lambda_n \geq 0, \hspace{0.5em} \forall \lambda_i \in \Z\}.
\end{equation}
For each $\lambda \in \Parn$, it can be regarded as a character $\lambda : B \rightarrow \C^*$. Define an action of $B$ on $G \times \C_\lambda$ by $b(g,t):= (gb^{-1}, \lambda(b)t)$. Construct a $G$-linearized line bundle
\begin{equation}\label{eqn1.4}
L_\lambda := (G\times \C_{-\lambda}) / B
\end{equation}
over $G/B$. If $X(w)$ has an $L_I$-action, so does the set of global sections $H^0(X(w),L_\lambda)$.

\begin{theorem}[{\cite[Conjecture 4.19]{MR4394321}}]\label{thm1.1}
Let $G=GL_n$ and $W = S_n$. Let $I \subseteq [n-1]$ and $\lambda \in \Parn$ be strictly decreasing. Assume $w \in W$ is the maximal representative of $W_Iw$. Then $X(w)$ is $L_I$-spherical if and only if $H^0(X(w) , L_\lambda)$ is $L_I$-multiplicity-free.
\end{theorem}

We postpone the proof until \S \ref{sub4.2}. Previously, it was known that $X(w)$ is $L_I$-spherical if and only if $H^0(X(w) , L_\lambda)$ is $L_I$-multiplicity-free for \textit{all} $\lambda \in \Parn$; see Theorem \ref{lem2.4}. Our theorem asserts that it is sufficient to check the multiplicity for a single $\lambda$ which is strictly decreasing.

The $T$-character of $H^0(X(w), L_\lambda)^*$ is given by the \textbf{key polynomial} $\kappa_{w\lambda}$; see, \textit{e.g.}, \cite[\S 4.4]{MR4394321}. To demonstrate, let $n=4$, $w = 2413$ and $I = \{1,3\}$. Choose a strictly decreasing $\lambda = (3,2,1,0)$ and expand $\kappa_{w\lambda}$ into a sum of products of Schur polynomials
\begin{equation}\label{eqn1.5}
\kappa_{w\lambda}(x_1 , x_2 , x_3 , x_4) = s_{3,1}(x_1 , x_2) s_{2,0}(x_3 , x_4) + s_{3,2}(x_1 , x_2) s_{1,0}(x_3 , x_4).
\end{equation}
Since each product of Schur polynomials appears at most once, $H^0 (X(w), L_\lambda)$ is $L_I$-multiplicity-free. From Theorem \ref{thm1.1}, $X(w)$ is $L_I$-spherical.

\subsection{Comparisons with the literature}\label{sub1.2}

We first recall some known results about spherical Schubert varieties overlapping Theorem \ref{thm1.2}. For $G= GL_n$, there are classifications for toric $X(wW_J)$ in partial flag varieties, $L_I$-spherical $X(wW_J)$ in Grassmannians and $L_I$-spherical $X(w)$ in complete flag varieties.

\begin{theorem}[{\cite{MR3044412}, \cite{MR4599567}, \cite{gao2024orbitstructurescomplexityschubert}}]\label{thm2.1}
Let $G = GL_n$, $W=S_n$, $I=\emptyset$, $J \subseteq [n-1]$ and $w \in W$. Then $X(wW_J)$ is toric if and only if $[w,J]$ avoids $321$, $3412$, $3312$, $2311$ and $2211$.
\end{theorem}
\begin{proof}
In \cite[Theorems 2 and 4]{MR3044412}, a classification of toric $X(w) \subseteq G/B$ was given. Using birational equivalence, this extends to $X(wW_J) \subseteq G/P_J$. Another proof can be given from a classification of multiplicity-free key polynomials in \cite[Theorem 1.1]{MR4599567}. For a generalization to other Lie types, see \cite[Corollary 4.12]{gao2024orbitstructurescomplexityschubert}.
\end{proof}

Recently, a formula for computing the minimal codimension of $T$-orbits of $X(w)$ in complete flag varieties of general Lie types was shown in \cite[Corollary 1.3]{gao2024orbitstructurescomplexityschubert}.

\begin{theorem}[{\cite[Theorem 1.1.1]{MR4517644}}]\label{thm2.2}
Let $G=GL_n$, $W=S_n$ and $I ,J \subseteq [n-1]$. Assume that $|J| = n-2$ so that $G/P_J$ is a Grassmannian. Let $u \in {}^I W^J$ and $w$ be the maximal representative of $W_I u W_J$. Then $X(wW_J)$ is $L_I$-spherical if and only if $[u,I,J]$ avoids the following patterns.
\begin{equation}\label{eqn2.1}
222i11i11, \quad 22i22i111, \quad 22i12i11, \quad 2i2i1i1, \quad 22ii1i1, \quad 2i2ii11.
\end{equation}
\end{theorem}
\begin{proof}
The first three patterns are due to the case $b_L = 3$ in \cite[Theorem 1.1.1]{MR4517644}, while the remaining patterns are due to the case $b_L=4$.
\end{proof}

\begin{theorem}[{\cite[Theorem 1.3]{MR4615612}}]\label{thm2.3}
Let $G = GL_n$, $W=S_n$, $I \subseteq [n-1]$ and $u,w \in W$. Assume that $u$ and $w$ are the minimal and maximal representatives of $W_I u$, respectively. Then $X(w)$ is $L_I$-spherical if and only if $u$ avoids $321$ and $3412$ permutation patterns.
\end{theorem}
\begin{proof}
In \cite[\S 1.1]{MR4615612}, $X(w)$ is $L_I$-spherical if and only if $u$ is a \textit{standard Coxeter element}. This is equivalent to avoiding $321$ and $3412$ permutation patterns; see \cite[Theorem 4.3]{MR2333139}. For generalizations to other Lie types, see \cite[Theorem 1.3]{MR4690502} or \cite[Theorem 1.1]{MR4643501}.
\end{proof}

Before ending this subsection, we discuss some topics related to the results and the proof techniques in this paper. The representation theoretic approach to sphericality is Theorem \ref{lem2.4}, which can be applied to other varieties relevant to Schubert calculus. For instance, we point out the classification of sphericality for the product of flag varieties with diagonal action. This was studied by P.~Littelmann \cite{MR1276821}, P.~Magyar, J.~Weyman, A.~Zelevinsky \cite{MR1667147} and J.~Stembridge \cite{MR2017064}. Here, the classification is related to multiplicity-free tensor products of finite-dimensional irreducible modules.

Fix a $G$-spherical homogeneous space $G/H$. Luna--Vust theory gives us a bijection between embeddings of $G/H$ into $G$-spherical varieties and combinatorial objects called \textit{colored fans}; see, \textit{e.g.}, \cite[Theorem 3.1.7]{MR3177371}. This relation generalizes fans used to classify toric varieties. To classify $G/H$, D.~Luna initiated in \cite{MR1896179} a classification program for \textit{$G$-wonderful varieties}, recently completed in \cite[Theorem 1.2.3]{MR3473657}. These varieties are $G$-spherical, smooth and projective; for the formal definition, see, \textit{e.g.}, \cite[\S 3.4.4]{MR3177371}. For instance, partial flag varieties $G/P$ are $G$-wonderful. For a connection to Schubert varieties, see, \textit{e.g.}, \cite[Theorem 1.5]{MR4643501}.

In Lemma \ref{lem3.4}, we multiply a key polynomial $\kappa_\alpha$ with a monomial $\kappa_\beta$ and expand $\kappa_\alpha \cdot \kappa_\beta = \sum_\gamma \kappa_\gamma$. A similar computation arises in other topics. P.~Magyar combined it with \textit{Demazure operators} to compute Schubert polynomials and called this process \textit{orthodontia} in \cite{MR1617793}. For the recent development in this direction, see, \textit{e.g.}, \cite{MR4369248}. Analogous to the relation between Littlewood-Richardson coefficients and Kostka numbers, we can also use it to expand $\kappa_\alpha$ into monomials. This leads to \textit{semi-skyline augmented fillings} which are related to \textit{nonsymmetric Macdonald polynomials}; see \cite{MR2496309}.

For Schubert varieties, other properties are characterized using permutation patterns. For instance, in \cite{MR1051089}, $X(w)$ of type A is smooth if and only if $w$ avoids $3412$ and $4231$ permutation patterns. In fact, $X(w)$ is smooth if and only if $X(w^{-1})$ is smooth for general Lie types, which motivates Corollary \ref{coro1.3}; see \cite{MR2827039}. Another classification of Schubert varieties when they are \textit{Gorenstein} is given by patterns in \cite{MR2264071}. For other examples related to Schubert varieties and permutation patterns, see, \textit{e.g.}, \cite{MR3012380}.

\subsection{Organization}

In \S \ref{sec3}, we discuss Demazure modules over finite-dimensional semisimple Lie algebras. In particular, we define Lakshmibai-Seshadri paths and their crystal graphs. From \S \ref{sec4}, we restrict to the case $G =GL_n$. By definition, $[u,I,J]$ is \textit{spherical} if and only if $X(wW_J)$ is $L_I$-spherical; see the introduction to \S \ref{sec4}. In \S \ref{sec5}, we prove that the patterns suggested in Theorem \ref{thm2.1} and Theorem \ref{thm2.3} imply non-sphericality. In \S \ref{sec6}, we reduce the remaining cases to the \textit{bigrassmannian} case, which is possible due to the technical Lemma \ref{lem6.2}.

The proof of Theorem \ref{thm1.1} is in \S \ref{sub4.2}; the reader only interested in Theorem \ref{thm1.2} may skip this subsection. At the end of the paper, we prove Theorem \ref{thm1.2} and its corollaries. We use Theorem \ref{thm2.2} in the proof of Theorem \ref{thm6.4} to simplify the argument. The proof of Theorem \ref{thm1.2} is independent of Theorem \ref{thm2.1} and Theorem \ref{thm2.3}. In fact, we can reprove Theorem \ref{thm2.3} using Lemma \ref{lem6.2}; see Remark \ref{rmk6.1}.

\section{Demazure modules}\label{sec3}

On a polynomial ring $\C [x_1 , \cdots , x_n]$, the \textbf{Demazure operator} is defined by
\begin{equation}\label{eqn3.31}
\Delta_j : \C [x_1 , \cdots , x_n] \rightarrow \C [x_1 , \cdots , x_n], \quad f \mapsto \frac{x_j f - x_{j+1} s_j f}{x_j - x_{j+1}},
\end{equation}
where $(s_j f)(x_1 , \cdots , x_n) := f(x_1, \cdots , x_{j+1} , x_j , \cdots , x_n)$. For each $w \in S_n$, choose a reduced word expression $w = s_{i_1} \cdots s_{i_l}$ and $\lambda = (\lambda_1 , \cdots , \lambda_n ) \in \Parn$. From \cite[Theorem 5]{MR1887117},
\begin{equation}
\kappa_{w\lambda} = \Delta_{i_1} \cdots \Delta_{i_l} (x_1^{\lambda_1} \cdots x_n^{\lambda_n}).
\end{equation}

Define $\Comp: = \{ \alpha = (\alpha_1 , \cdots , \alpha_n ) \mid \forall \alpha_i \in \Z_{\geq 0}\}$. Since $\kappa_\beta (x_1 , \cdots , x_m) \kappa_{\gamma}(y_1 , \cdots , y_n)$ forms a linear basis of $\C [x_1 , \cdots , x_m , y_1 , \cdots , y_n]$,
\begin{equation}
\kappa_\alpha (x_1 , \cdots , x_m , y_1 , \cdots , y_n) = \sum_{\beta \in {\rm Comp}_m, \gamma \in {\rm Comp}_n} c_{\beta , \gamma}^\alpha \kappa_\beta (x_1 , \cdots , x_m) \kappa_\gamma (y_1 , \cdots , y_n),
\end{equation}
for any $\alpha \in {\rm Comp}_{m+n}$ and for some $c_{\beta,\gamma}^\alpha \in \R$. In fact, the coefficient $c_{\beta , \gamma}^\alpha$ is a non-negative integer according to Lemma \ref{lem3.2} or Lemma \ref{lem3.4}. For a large enough positive integer $N$, $c_{\beta, \gamma}^\alpha$ also appears in
\begin{subequations}
\begin{equation}
\kappa_{\alpha} \cdot \kappa_{N \xi} = \sum_{\beta \in {\rm Comp}_m, \gamma \in {\rm Comp}_n} c_{\beta , \gamma}^\alpha \cdot  \kappa_{\tau_1 (\beta ) + \tau_2 (\gamma) + N \xi }, \quad \xi:=( \overbrace{1,1, \cdots ,1}^{\# m},\overbrace{0,0, \cdots ,0}^{\# n}),
\end{equation}
\begin{equation}
\tau_1 : {\rm Comp}_{m} \rightarrow {\rm Comp}_{m+n}, \quad (\beta_1 , \cdots , \beta_m) \mapsto (\beta_1 , \cdots , \beta_m , 0 , \cdots , 0),
\end{equation}
\begin{equation}
\tau_2 : {\rm Comp}_{n} \rightarrow {\rm Comp}_{m+n}, \quad (\gamma_1 , \cdots , \gamma_n) \mapsto (0 , \cdots , 0, \gamma_1 , \cdots , \gamma_n).
\end{equation}
\end{subequations}
In particular, if $\alpha , \beta , \gamma$ are weakly increasing \textit{i.e.} reversed partitions $\lambda , \mu , \nu$, then $c_{\beta , \gamma}^\alpha$ is the \textbf{Littlewood-Richardson coefficient} of $\lambda , \mu , \nu$. 

A combinatorial rule of $c_{\beta , \gamma}^\alpha$ is given by using the crystal structure of Demazure modules as in Theorem \ref{thm3.3}. We use this to determine the case when $c_{\beta , \gamma}^\alpha$ is either $0$ or $1$ in order to prove sphericality; see \cite[Section 4.4]{MR4615612}.

\subsection{Annihilators of weight vectors}

Let $\mathfrak{g}$ be a semisimple Lie algebra over $\C$ of rank $n$ and $\mathfrak{h}$ be a Cartan subalgebra of $\mathfrak{g}$. Let $\Phi \subseteq \mathfrak{h}^*$ be the set of roots and choose subsets of simple roots and positive roots, denoted by $\Delta$ and $\Phi^+$, respectively. Let $\Lambda \subseteq \mathfrak{h}^*$ be the lattice of integral weights and $\Lambda^+$ be the subset of dominant integral weights. In particular, write the simple root as $\alpha_i$ and the fundamental weight as $\varpi_i$ satisfying $( \varpi_i , \alpha_j^\vee ) = \delta_{i,j}$, where $\alpha^\vee$ denotes the coroot of $\alpha$. Write a Chevalley basis of $\mathfrak{g}$ as $\{H_i\}_{1 \leq i \leq n} \cup \{X_\alpha\}_{\alpha \in \Phi}$.

Let $W$ be the Weyl group associated to $\mathfrak{g}$. Let $\mathfrak{b}$ be a Borel subalgebra of $\mathfrak{g}$ generated by $\{H_i \}_{1 \leq i \leq n} \cup \{X_\alpha \}_{\alpha \in \Phi^+}$. Finite-dimensional irreducible $\mathfrak{g}$-modules are characterized by highest weights $\lambda \in \Lambda^+$. We write an irreducible module as $V_\lambda$ of highest weight $\lambda$ with highest weight vector $v_\lambda^+$.

For each $w \in W$, there exists the unique weight vector of weight $w\lambda$ in $V_\lambda$ up to scalar multiplication. We call it an \textbf{extremal weight vector} and denote it by $v^w_\lambda$. Define a $\mathfrak{b}$-module
\begin{equation}\label{eqn3.1}
V_\lambda (w):= \mathcal{U} (\mathfrak{b} )  v_\lambda^w.
\end{equation}
$V_\lambda (w)$ is called a \textbf{Demazure module}.

\begin{theorem}[{Theorem 3.4 \cite{MR0826100}}]\label{thm3.1}
Let $\lambda \in \Lambda^+$ and $w \in W$. For each $\alpha \in \Phi^+$, define
\begin{equation}\label{eqn3.2}
k_\alpha :=
\begin{cases}
- ( w \lambda , \alpha^\vee ) , & ( w \lambda , \alpha^\vee )  \leq 0 \\
0, & ( w \lambda , \alpha^\vee )  > 0
\end{cases}
.
\end{equation}
Then
\begin{equation}\label{eqn3.3}
{\rm Ann \hspace{0.2em}}_{\mathcal{U}(\mathfrak{b})} v^w_\lambda = \sum_{\alpha \in \Phi^+} \mathcal{U}(\mathfrak{b}) X_\alpha^{k_\alpha +1} + \sum_{1 \leq i \leq n} \mathcal{U}(\mathfrak{b}) \left( H_i - ( w \lambda , \alpha_i^\vee ) \cdot 1  \right).
\end{equation}
\end{theorem}

Let $I \subseteq [n]$. Write
\begin{equation}\label{eqn3.4}
\Delta_I := \{\alpha_i \mid i \in I\}, \quad \Phi_I:= \{ \alpha \in \Phi \mid \alpha = \sum_{i \in I} c_i \alpha_i\}, \quad \Phi_I^+ := \Phi_I \cap \Phi^+.
\end{equation}
Define Lie subalgebras of $\mathfrak{g}$
\begin{equation}\label{eqn3.5}
\mathfrak{g}_I := \langle H_i \mid i \in I \rangle \oplus \langle X_\alpha \mid \alpha \in \Phi_I \rangle , \quad \mathfrak{b}_I := \langle H_i \mid i \in I \rangle \oplus \langle X_\alpha \mid \alpha \in \Phi_I^+ \rangle.
\end{equation}
Write $[n]-I = \{y_1 < y_2 < \cdots < y_{g-1}\}$, $y_0 :=0$ and $y_g:=n+1$. Rearrange simple roots so that for each $1 \leq p < q \leq n$ such that $p, q \in I$, if there exists $y \in [n]-I$ such that $p < y < q $, then $(\alpha_p , \alpha_q) = 0$. Indeed, recall that in the \textit{Dynkin diagram} of $\mathfrak{g}$, each $\alpha_i$ corresponds to a vertex $i$. We rearrange simple roots so that after removing vertices $y \in [n]-I$, indices of vertices in each connected component are between some $y_i$ and $y_{i+1}$.

For each $1 \leq i \leq g$, let $\Phi_i := \Phi \cap \langle \alpha_j \mid y_{i-1} < j < y_i \rangle$ and $\Phi_i^+ := \Phi_i \cap \Phi^+$ and
\begin{equation}
\mathfrak{g}_i := \langle H_j \mid y_{i-1} < j < y_{i} \rangle \oplus \langle X_\alpha \mid \alpha \in \Phi_i \rangle, \hsp \mathfrak{b}_i := \langle H_j \mid y_{i-1} < j < y_{i} \rangle \oplus \langle X_\alpha \mid \alpha \in \Phi_i^+ \rangle.
\end{equation}
In particular, each $\mathfrak{g}_i$ is a semisimple Lie algebra. Also, $\mathfrak{g}_I = \mathfrak{g}_1 \oplus \cdots \oplus \mathfrak{g}_g$ and $\mathfrak{b}_I = \mathfrak{b}_1 \oplus \cdots \oplus \mathfrak{b}_g$. Define
\begin{equation}\label{eqn3.6}
\Lambda_i := \bigoplus_{j=y_{i-1}+1}^{y_i-1} \Z \varpi_{j}, \quad \Lambda_i^+ := \Lambda_i \cap \Lambda^+,
\end{equation}
and the projection $p_i : \Lambda \rightarrow \Lambda_i$. Also, $W_i$ is a subgroup of $W$ generated as $\langle s_{j} \mid y_{i-1} < j < y_i \rangle$.

Let $\lambda \in \Lambda^+$ and $w \in W$. Let $N$ be a large enough positive integer. Define
\begin{equation}\label{eqn3.7}
\xi:= \sum_{i \notin I} \varpi_i.
\end{equation}
Write the identity of $W$ as $e$. Due to \cite[Theorem~1]{MR1016897}, $V_\lambda (w) \otimes V_{N \xi} (e)$ has a filtration whose subquotients are Demazure modules. In other words,
\begin{equation}\label{eqn3.8}
V_\lambda (w) \otimes V_{N \xi} (e) = \mathcal{F}_m \supseteq \mathcal{F}_{m-1} \supseteq \cdots \supseteq \mathcal{F}_0 = 0, \quad \mathcal{F}_j / \mathcal{F}_{j-1} \cong V_{\lambda_j}(w_j) \quad 1 \leq j \leq m,
\end{equation}
for some $w_j \in W$ and $\lambda_j \in \Lambda^+$. On the other hand,
\begin{equation}\label{eqn3.9}
{\rm res}_{\mathfrak{b}_I}^\mathfrak{b} V_\lambda (w) \cong {\rm res}_{\mathfrak{b}_I}^\mathfrak{b} \left( V_\lambda (w) \otimes V_{N \xi}(e) \right),
\end{equation}
since $V_{N \xi}(e)$ is a trivial module after restriction.

Let $M$ be a finite-dimensional $\mathfrak{h}$-module admitting a weight-space decomposition with weights in $\Lambda$. Then define the formal character ${\rm ch} (M)$ in the group ring $\Z [\Lambda]$ with basis $e^\lambda$ where $\lambda \in \Lambda$. Define a ring homomorphism
\begin{equation}\label{eqn3.10}
q_I: \Z [\Lambda] \rightarrow \Z [\Lambda ] , \quad e^{\varpi_j} \mapsto e^{\varpi_j} \quad (j \in I), \quad  e^{\varpi_j} \mapsto 1 \quad (j \notin I).
\end{equation}

\begin{lemma}\label{lem3.2}
Let $\mathfrak{g}$ be a semisimple Lie algebra of rank $n$ with Weyl group $W$. Let $\lambda \in \Lambda^+$, $w \in W$ and $I \subseteq [n]$. Then there exists a sequence of $\lambda_j \in \Lambda^+$ and $w_j \in W$ such that
\begin{equation}
q_I \left( {\rm ch} \left( V_\lambda (w) \right) \right) = \sum_j \prod_{1 \leq i \leq g} {\rm ch} \left( V_{\lambda_{i,j}}(w_{i,j} ) \right),
\end{equation}
where
\begin{equation}
{\rm ch}\left( V_\lambda (w) \right) \cdot {\rm ch}\left( V_{N\xi} (e) \right) = \sum_{j} {\rm ch}\left( V_{\lambda_j} (w_j) \right),
\end{equation}
for a large enough integer $N$ and $\xi := \sum_{i \notin I} \varpi_i$. Here, $\lambda_{i,j} \in \Lambda_i^+$ and $w_{i,j} \in W_i$ are determined by $p_i (w_{i,j} \lambda_{i,j} )= p_i (w_j \lambda_j)$.
\end{lemma}

If $y_{i-1}+1 > y_i - 1$, ${\rm ch} \left( V_{\lambda_{i,j}}(w_{i,j} )  \right) = 1$ as a convention.

\begin{proof}
Since $p_i (w_j \lambda_j) \in \Lambda_i$, there exists $w_{i,j} \in W_i$ such that $p_i(w_{i,j}^{-1} p_i (w_j \lambda_j)) \in \Lambda_i^+$. Let $\lambda_{i,j} = p_i(w_{i,j}^{-1} p_i (w_j \lambda_j))$.

We need to show that if $N$ is a large enough integer, then as $\mathfrak{b}_I$-modules,
\begin{equation}\label{eqn3.11}
{\rm res}_{\mathfrak{b}_I}^\mathfrak{b} V_{\lambda_j}(w_j) \cong V_{\lambda_{1,j}}(w_{1,j}) \boxtimes \cdots \boxtimes V_{\lambda_{g,j}}(w_{g,j}).
\end{equation}
Here, the right hand side is an external tensor product of representations of $\mathfrak{b}_I = \mathfrak{b}_1 \oplus \cdots \oplus \mathfrak{b}_g$.
It is enough to prove that
\begin{equation}\label{eqn3.12}
\mathcal{U} (\mathfrak{b}) v^{w_j}_{\lambda_j} = \mathcal{U} (\mathfrak{b}_I) v^{w_j}_{\lambda_j}.
\end{equation}
Let $\alpha \in \Phi^+$ and write $\alpha = \sum_{i=1}^{n} c_i \alpha_i$. Then
\begin{equation}\label{eqn3.13}
( w_j \lambda_j , \alpha^\vee ) = \frac{2}{( \alpha , \alpha )}\sum_{i=1}^{n} c_i ( w_j \lambda_j , \alpha_i ).
\end{equation}
Assume that $\alpha \notin \Phi_I^+$. Then there exists a $y \notin I$ such that $c_y >0$. Since $w_j \lambda_j$ appears as a weight in $V_\lambda (w) \otimes V_{N\xi} (e)$, $( w_j \lambda_j , \alpha_y )$ is sufficiently large to make $( w_j \lambda_j , \alpha^\vee )$ positive. Due to Theorem \ref{thm3.1},
\begin{equation}\label{eqn3.14}
X_\alpha v^{w_j}_{\lambda_j} = 0.
\end{equation}
Together with the PBW Theorem, this leads to \eqref{eqn3.12}. From Theorem \ref{thm3.1}, we have \eqref{eqn3.11}.
\end{proof}

\subsection{Lakshmibai-Seshadri paths}

In this subsection, we introduce Lakshmibai-Seshadri(LS) paths in \cite{MR1887117} for a semisimple Lie algebra $\mathfrak{g}$. Let $\Lambda_\R := \Lambda \otimes_\Z \R$. A \textbf{path} is a piecewise-linear map $\pi : [0,1]\rightarrow \Lambda_\R$ up to reparametrization with $\pi (0) = 0$. Write $\pi_1 * \pi_2$ as a concatenation of two paths. Formally,
\begin{equation}\label{eqn3.15}
(\pi_1 * \pi_2) (t) =
\begin{cases}
\pi_1 (2t) & 0\leq t \leq 1/2 \\
\pi_2 (2t-1) + \pi_1 (1) & 1/2 \leq t \leq 1
\end{cases}
.
\end{equation}
A path $\pi$ is \textbf{dominant} if $( \pi (t) , \alpha_i^\vee ) \geq 0$ for all $t \in [0,1]$ and $1 \leq i \leq n$.

The \textbf{weight} of a path $\pi$ is $\pi(1)$ and denoted by ${\sf wt} (\pi)$. For each $\lambda \in \Lambda$, write the path $\pi^\lambda : [0,1] \rightarrow \Lambda_\R$ as a straight line \textit{i.e.} $\pi^\lambda (t) := t \lambda$. Let $\lambda \in \Lambda^+$. Write the stabilizer subgroup $W_\lambda$ of $W$ for $\lambda$. An \textbf{LS-chain} of $\lambda$ is a pair of lists
\begin{equation}\label{eqn3.16}
(\tau_1 > \tau_2 > \cdots > \tau_q ; \hspace{0.5em}  0 = a_0 < a_1 < \cdots < a_q = 1),
\end{equation}
where $\tau_j \in W/W_\lambda$ and $a_j \in \Q$. Here, we have an additional condition that for each $j$, there exists a chain of Bruhat order $\tau_j = \sigma_0 > \sigma_1 > \cdots > \sigma_p = \tau_{j+1}$ with $\ell (\sigma_{k} ) = \ell (\sigma_{k+1})+1$ and $a_j (\sigma_{k+1} \lambda - \sigma_k \lambda)$ contained in the root lattice $\bigoplus_{i=1}^{n} \Z \alpha_i$.

Construct a path $\pi : [0,1] \rightarrow \Lambda_\R$ from the LS-chain
\begin{equation}\label{eqn3.17}
\pi (t) := \sum_{j=1}^{k-1} (a_j - a_{j-1})\tau_j \lambda  + (t - a_{k-1})\tau_k \lambda, \quad a_{k-1} \leq t \leq a_k.
\end{equation}
We call it an \textbf{LS-path} of $\lambda$. Due to the construction, an LS-path is defined by a unique LS-chain. We write an LS-path $\pi = (\tau_1 > \cdots > \tau_q ; a_0 < \cdots < a_q)$. For instance, $\pi^\lambda$ is also an LS-path of $\lambda$ \textit{i.e.} $\pi^\lambda = (e;0<1)$.

From now on, we only consider a path which is a concatenation of some LS-paths. For each $1 \leq i \leq n$, we want to define a \textbf{lowering root operator} $f_i$ acting on a path $\pi$. Let $Q$ be the minimum value of $( \pi (t ) , \alpha_i^{\vee} )$ for $t \in [0,1]$. Let $t_1$ be the largest $t$ satisfying $Q = ( \pi (t ) , \alpha_i^{\vee} )$. Let $t_2$ be the smallest $t > t_1$ satisfying $Q+1 =  ( \pi (t ) , \alpha_i^{\vee} )$, if it does exist. Split the path $\pi$ into $\pi_1 * \pi_2 * \pi_3$ in accordance with the three intervals $[0,t_1]$, $[t_1 , t_2]$ and $[t_2 , 1]$. Then $f_i (\pi) := \pi_1 * s_i \pi_2 * \pi_3$ where $s_i \pi_2 (t):= s_i (\pi_2(t))$. If such $t_2$ does not exist, then $f_i (\pi)$ is undefined. If $f_i (\pi) = \pi'$, then ${\sf wt}(\pi) -\alpha_i= {\sf wt}(\pi')$.

Similarly, define a \textbf{raising root operator} $e_i$ as follows. Let $Q$ be defined the same as before. Let $t_2'$ be the smallest $t$ satisfying $Q = ( \pi (t ) , \alpha_i^{\vee} )$. Choose $t_1'$ to be the largest $t < t_2'$ with $Q+1 =  ( \pi (t ) , \alpha_i^{\vee} )$, if it does exist. Split the path $\pi$ into three pieces using $[0,t_1']$, $[t_1' , t_2']$ and $[t_2',1]$. Same as before, define $e_i (\pi) := \pi_1 * s_i \pi_2 * \pi_3$. If such $t_1'$ does not exist, then $e_i (\pi)$ is undefined. If $e_i (\pi) = \pi'$, then ${\sf wt}(\pi) +\alpha_i= {\sf wt}(\pi')$.

Let $w \in W$ and $\lambda \in \Lambda^+$. Fix a reduced word expression of $w = s_{i_1} s_{i_2} \cdots s_{i_l}$. Record this into $\mathbf{i} := (i_1 , i_2 , \cdots , i_l)$. For each subset $J \subseteq [l]$, write $J = \{ j_1  , \cdots , j_m \}$ where $j_1 < \cdots < j_m$. Write $\mathbf{i}(J) := (i_{j_1} , i_{j_2}, \cdots , i_{j_m})$. Also, denote $w(\mathbf{i}(J)):= s_{i_{j_1}} s_{i_{j_2}} \cdots s_{i_{j_m}}$.

Write $\lambda = \sum_{i} l_i \varpi_i$ and for each $1 \leq j \leq l$, define
\begin{equation}\label{eqn3.18}
m_j = 
\begin{cases}
l_k & \text{if }i_j = k \text{ and }i_{j'} \neq k \text{ for all }j' > j \\
0 & \text{otherwise}
\end{cases}
.
\end{equation}

Write $\mathbf{m} = (m_1 , m_2 , \cdots , m_l)$ and $\lambda_j:= m_j \varpi_{i_j}$. Let $\pi_j$ be an LS-path of shape $\lambda_j$ for each $1 \leq j \leq l$. From the corresponding LS-chain, write $\pi_j$ as $(\tau_{j,1} > \tau_{j,2} > \cdots > \tau_{j,p_j} ; a_{j,0} < a_{j,1} < \cdots < a_{j,p_j})$ where each $\tau_{j,p} \in W/W_{m_j \varpi_{i_j}}$. A concatenation $\Pi = \pi_1 * \cdots * \pi_l$ is called an \textbf{LS-monomial} of $(\mathbf{i} , \mathbf{m})$ if there exists a chain of index sets
\begin{equation}\label{eqn3.19}
[l] \supseteq J_{1,1} \supseteq \cdots \supseteq J_{1, p_1} \supseteq \cdots \supseteq J_{l,1} \supseteq \cdots \supseteq J_{l,p_l}
\end{equation}
such that $\mathbf{i}(J_{j,p} \cap [j])$ is a reduced expression and
\begin{equation}\label{eqn3.20}
w \left(  \mathbf{i}(J_{j,p} \cap [j]) \right) \equiv \tau_{j,p} \quad {\rm modulo} \hspace{0.2em} W_{m_j \varpi_{i_j}}.
\end{equation}

Suppose there are finitely many paths, each of them corresponding to a vertex. If a pair of paths $\Pi$ and $\Pi'$ satisfies $f_i (\Pi) = \Pi'$, we connect them with a directed edge labelled by $i$. The graph composed of vertices with weights and directed edges labelled by $i$ is called a \textbf{crystal graph} of paths.

Let $\lambda \in \Lambda^+$, $w \in W$ and choose $(\mathbf{i} , \mathbf{m})$ as before. Define $K := \{ k \in [n] \mid \exists 1 \leq j \leq l \hsp \text{such that} \hsp k = i_j \} $. Then $\sum_{j=1}^l \lambda_j= \sum_{i \in K} l_i \varpi_i$. Let $\nu := \sum_{i \notin K} l_i \varpi_i$. For a moment, we call $\pi^\nu * \Pi$ an LS-monomial of $(\mathbf{i} , \mathbf{m} ; \lambda)$. From LS-monomials of $(\mathbf{i} , \mathbf{m} ; \lambda)$, we get the crystal graph of the Demazure module $V_\lambda (w)$, independent of the choices of $(\mathbf{i} , \mathbf{m})$. This graph is connected. There is a weight-preserving bijection between a linear basis of $V_\lambda (w)$ and the set of LS-monomials of $(\mathbf{i} , \mathbf{m} ;\lambda)$.

From a crystal graph of LS-monomials $\Pi$ of $(\mathbf{i} , \mathbf{m})$, add the weight $\nu$ to each vertex. The claim is that it is exactly the same as the crystal graph of LS-monomials $\pi^\nu * \Pi$ of $(\mathbf{i} , \mathbf{m}; \lambda)$ by mapping $\Pi$ to $\pi^\nu * \Pi$. To show this, note that if $i \notin K$, then $s_i$ does not appear in a reduced expression of $w$. Therefore, ${\sf wt}(\Pi) \in \{ \xi \in \Lambda \mid \lambda - \nu - \xi \in \oplus_{i \in K} \Z_{\geq 0} \alpha_i \}$. In other words, for any pair of LS-monomials $\Pi$ and $\Pi'$ of $(\mathbf{i} , \mathbf{m})$ and $i \notin K$, ${\sf wt}(\Pi) - {\sf wt}(\Pi') \neq \alpha_i$. Therefore, there is no edge labelled by $i$ between $\Pi$ and $\Pi'$ or between $\pi^\nu * \Pi$ and $\pi^\nu *\Pi'$. On the other hand, suppose $i \in K$. Then $s_i \nu = \nu$ since $(\nu , \alpha_i^\vee ) = 0$. Therefore, $\pi^\nu$ does not change the computation of $f_i$. In other words, if $i \in K$, then $f_i (\Pi) = \Pi'$ if and only if $f_i (\pi^\nu * \Pi) = \pi^\nu * \Pi'$. This shows that the claim is true.

In fact, $V_{\lambda}(w) \cong V_{\lambda - \nu}(w) \otimes V_\nu (e)$ as $\mathfrak{b}$-modules; see \cite[\S 1.1]{MR1887117}. Similarly, for any $\mu \in \Lambda^+$, $1 \leq i \leq n$ and a pair of LS-monomials $\Pi$ and $\Pi'$ of $(\mathbf{i} , \mathbf{m})$, $f_i (\pi^\mu * \Pi) = \pi^\mu * \Pi'$ if and only if $f_i(\pi^{\mu  + \nu} * \Pi) = \pi^{\mu + \nu} * \Pi'$.

\begin{theorem}[{\cite[Theorem 2]{MR1887117}}]\label{thm3.3}
Choose $(\mathbf{i} , \mathbf{m})$ from $\lambda \in \Lambda^+$ and $w \in W$ as before. Let $\mu \in \Lambda^+$. Then the crystal graph on the set of $\{\pi^\mu * \Pi \}$ is a disjoint union of crystal graphs of Demazure modules where $\Pi$ runs over LS-monomials of $(\mathbf{i} , \mathbf{m})$. Also, there is a bijection between each copy of a Demazure module and a dominant path $\pi^\mu * \Pi $. In particular, a Demazure module $V_{\lambda'} (w')$ is mapped to a dominant path with weight $\lambda'$.
\end{theorem}
\begin{proof}
See \cite[Proposition 12]{MR1887117}. Each copy of a Demazure module can be bijectively mapped to a path $\pi^\mu * \Pi $ where $e_i $ is undefined for all $i$. This is equivalent to $( (\pi^\mu * \Pi)(t) , \alpha_i^\vee ) \geq 0$ so that $Q = 0$ and $t_2 ' =0$.
\end{proof}

As before, let $\nu = \sum_{i \notin K} l_i \varpi_i$. $\{ \pi^\mu * \Pi \}$ is a crystal graph of $V_{\lambda - \nu}(w) \otimes V_\mu (e)$. On the other hand, $\{ \pi^{\mu+ \nu} * \Pi \}$ is a crystal graph of $V_{\lambda - \nu}(w) \otimes V_{\mu + \nu} (e) \cong V_{\lambda- \nu} (w) \otimes V_\mu (e) \otimes V_\nu (e) \cong V_\lambda(w) \otimes V_\mu (e)$. For each $V_\lambda' (w')$ appearing in $\{ \pi^\mu * \Pi \}$, $w'$ may be chosen in $W_K$; see the construction of $w'$ in \cite[Proposition 12]{MR1887117}. Under the bijection between $\{ \pi^\mu * \Pi \}$ and $\{ \pi^{\mu + \nu} * \Pi \}$, a crystal graph of $V_{\lambda'}(w')$ is mapped to a crystal graph of $V_{\lambda' + \nu} (w')$.

\section{Levi subgroup multiplicity}\label{sec4} 

From now on, let $G = GL_n$, a complex general linear group. Let $W= S_n$. In this paper, a \textbf{minimal sequence} of $n$ is a finite sequence with numerical entries $1,2, \cdots , n$ and entries equal to the symbol $i$ with conditions as follows.
\begin{itemize}
\item There are exactly $n$ numerical entries.
\item Each entry with $i$ has a successor and a predecessor, both of which are numbers.
\item For consecutive number entries, the successor is more than or equal to the predecessor.
\item If the number $k$ appears, either $k-1$ also appears or $k=1$.
\end{itemize}
For instance, there are thirty-three minimal sequences of $3$
\begin{equation}\label{eqn1.10}
\begin{matrix}
123, & 112, & 122, & 111, & 1i23, & 2i13, & 3i12, & 1i12, & 2i11, & 1i22, \\
2i12, & 1i11, & 12i3, & 13i2, & 23i1, & 11i2, & 12i1, & 12i2, & 22i1, & 11i1, \\
1i2i3, & 1i3i2, & 2i1i3, & 2i3i1, & 3i1i2, & 3i2i1, & 1i1i2, & 1i2i1, & 2i1i1, & 1i2i2, \\
2i1i2, & 2i2i1, & 1i1i1.
\end{matrix}
\end{equation}

Let $I , J \subseteq [n-1]$ and $u \in {}^I W^J$. $(u,I,J) \mapsto [u,I,J]$ bijectively maps a triple $(u,I,J)$ to a minimal sequence of $n$. Let $w$ be the maximal representative of $W_I u W_J$. If $X(wW_J)$ is $L_I$-spherical, then we say $[u,I,J]$ is \textbf{spherical}. Otherwise, $[u,I,J]$ is \textbf{non-spherical}. For instance, all minimal sequences in \eqref{eqn1.10} are spherical except for $3i2i1$. The question \eqref{eqn1.9} can now be restated.
\begin{equation}\label{eqn1.11}
\textit{Which minimal sequences are spherical?}
\end{equation}

\subsection{Sortings} 

Let $G = GL_n$ and $W=S_n$. Let $S= \{ s_1 , \cdots , s_{n-1}\}$ be a set of minimal generators of $W$ with relations $s_i s_{i+1} s_i = s_{i+1} s_i s_{i+1}$ and $s_i s_j = s_j s_i$ if $|i-j| >1$.

Let $I,J \subseteq [n-1]$. Write
\begin{subequations}\label{eqn1.8}
\begin{equation}
[n-1]-I = \{ y_1 , \cdots , y_{g-1}\}, \quad y_1 < \cdots < y_{g-1}, \quad y_0:=0, \quad y_g:=n.
\end{equation}
\begin{equation}
[n-1]-J = \{ z_1 , \cdots , z_{h-1}\}, \quad z_1 < \cdots < z_{h-1}, \quad z_0:=0, \quad z_h:=n.
\end{equation}
\end{subequations}
If we choose a Borel subgroup $B$ as the set of upper triangular matrices and a maximal torus $T$ as the set of diagonal matrices, then a Levi subgroup is chosen by
\begin{equation}\label{eqn1.2}
L_I:= GL_{y_1 - y_0} \times \cdots \times GL_{y_{g} - y_{g-1}}.
\end{equation}
Also, the parabolic subgroup $P_J$ is given by the set of block upper triangular matrices of size $z_1 - z_0 , \cdots z_{h} - z_{h-1}$. 

Write the Grassmannian ${\rm Gr}(d,n)$ as the set of $d$-dimensional subspaces of $\C^n$. Using the Pl{\" u}cker embedding, it is a projective variety. In general, let $(d_1 , \cdots , d_m)$ be a sequence of positive integers such that $n = \sum_j d_j$. A \textbf{partial flag} of type $(d_1 , \cdots , d_m)$ in $\C^n$ is a nested sequence of vector spaces
\begin{equation}
0 = V_0 \subseteq V_1 \subseteq \cdots \subseteq V_m = \C^n, \quad {\rm dim}_{\C} (V_j / V_{j-1}) = d_j \hsp \hsp (1 \leq j \leq m).
\end{equation}
The set of flags of type $(d_1 , \cdots , d_m)$ can be identified with $G/ P_J$ by choosing $J \subseteq [n-1]$
\begin{equation}
J := [n-1] - \{ d_1 , d_1 + d_2 , d_1+d_2+d_3 , \cdots \}.
\end{equation}
$G/P_J$ is a subvariety of a product of Grassmannians.

Write
\begin{equation}\label{eqn2.2}
\Parj := \{ \lambda \in \Parn \mid \forall j \in J, \hspace{0.5em} \lambda_j = \lambda_{j+1}\}.
\end{equation}
The \textbf{character} of $P_J$ is a group homomorphism $P_J \rightarrow \C^*$ which is an algebraic map. Consider the Levi decomposition $P_J = L_J \ltimes R_u(P_J)$. The character of $L_J$ can be constructed as a product of powers of determinants using $\lambda \in \Parj$. This can be extended as a character of $P_J$ by using $P_J \rightarrow L_J$. Write the character as $\lambda : P_J \rightarrow \C^*$. Define an action of $P_J$ on $G \times \C_\lambda$ by $p(g,t) := (g p^{-1} , \lambda(p)t)$ for each $p \in P_J$. Then we have a $G$-linearized line bundle $(G \times \C_{-\lambda})/P_J$ over $G/P_J$, also denoted by $L_\lambda$; see, \textit{e.g.}, \cite{MR1782635} or \cite{MR2015057}. $H^0(X(wW_J),L_\lambda)^*$ gives another construction of Demazure modules as $B$-modules. 

\begin{theorem}[{\cite[Theorem 4.13]{MR4394321}}]\label{lem2.4}
Let $G= GL_n$ and $W= S_n$. Let $I, J \subseteq [n-1]$. Let $w \in W$ be the maximal representative of $W_I w W_J$. Then $X(wW_J)$ is $L_I$-spherical if and only if for all $\lambda \in \Parj$, $H^0(X(wW_J),L_\lambda)$ is $L_I$-multiplicity-free.
\end{theorem}
\begin{proof}
In \cite[Theorem~4.13]{MR4394321}, the statement is true when $J=\emptyset$. This generalizes to an arbitrary subset $J \subseteq [n-1]$. In \cite{MR0630769}, $X(wW_J)$ is a normal variety. Therefore, $X(wW_J)$ is $L_I$-spherical if and only if for any $L_I$-linearized line bundle $L$ over $X(wW_J)$, $H^0(X(wW_J),L)$ is $L_I$-multiplicity-free; see \cite[Theorem~2.1.2]{MR3177371}. This proves ``$\Rightarrow$'' direction. 

For ``$\Leftarrow$'' direction, choose $\lambda \in \Parj$ such that $\lambda_j > \lambda_{j+1}$ if $j \notin J$. Then $L_\lambda$ is a very ample line bundle over $X(wW_J)$; see, \textit{e.g.}, \cite[II \S 8]{MR2015057}. In the proof of \cite[Theorem~2.1.2]{MR3177371}, if $H^0(X(wW_J) , L_{n \lambda} )$ is $L_I$-multiplicity-free for all $n \in \Z_{>0}$, then $X(wW_J)$ is $L_I$-spherical. Since $n \lambda \in \Parj$, we are done.
\end{proof}

Recall that ${}^I W^J$ is the set of $w \in W$ such that $w$ is the minimal representative of $W_I w W_J$. Similarly, define ${}^I W$ and $W^J$ by replacing $W_I w W_J$ with $W_I w$ and $wW_J$, respectively. Define \textbf{Bruhat order} on $W/W_I$
\begin{equation}\label{eqn2.3}
uW_I \leq vW_I \quad \Leftrightarrow \quad X(uW_I) \subseteq X(vW_I).
\end{equation}

\begin{lemma}\label{lem3.0}
Let $W=S_n$ and $I,J \subseteq [n-1]$.
\begin{itemize}
\item[(a)] Let $u, v \in {}^I W^J$. Then $u \leq v$ if and only if there exist $u' \in W_I u W_J$ and $v' \in W_I v W_J$ such that $u' \leq v'$.
\item[(b)] Let $u \in {}^I W^J$. Then for each $w \in W_I u W_J$, there exist $w_I \in W_I$ and $w_J \in W_J$ such that $w = w_I u w_J$ and $\ell (w) = \ell (w_I) + \ell (u) + \ell (w_J)$.
\item[(c)] For each ${}^I w \in {}^I W$ and $w^J \in W^J$, the minimal representatives of ${}^I w W_J$ and $W_I w^J$ are in ${}^I W^J$.
\end{itemize}
For (a) and (c), the analogous statements are true for maximal representatives.
\end{lemma}
\begin{proof}
(a) See \cite[Lemma 3.5]{MR0435249}.

(b), (c) See \cite[Proposition 1.8]{MR2120109}.
\end{proof}

\begin{lemma}\label{lem2.5}
Let $G = GL_n$ and $W = S_n$. Let $I, J \subseteq [n-1]$, $\lambda \in \Parj$ and assume that $u , v \in W$ are maximal representatives of $W_I u$ and $W_I v$, respectively. Then if $uW_J \leq vW_J$, there is a natural $L_I$-module monomorphism
\begin{equation}\label{eqn2.4}
H^0 (X(uW_J),L_\lambda)^* \hookrightarrow H^0(X(vW_J),L_\lambda)^*.
\end{equation}
\end{lemma}
\begin{proof}
The natural map is defined from the inclusion $X(uW_J) \subseteq X(vW_J)$. Injectivity follows from the fact that $H^0(X(uW_J),L_\lambda )^* \rightarrow H^0(G/P_J , L_\lambda)^*$ is injective; see, \textit{e.g.}, \cite[II Proposition~14.15]{MR2015057}.
\end{proof}

\begin{lemma}\label{lem2.6}
Let $G=GL_n$, $W=S_n$, $I,J \subseteq [n-1]$ and $u \in {}^I W$. Denote $[u,J] = x_1 \cdots x_n$. Whenever we have a pair of $1\leq p<q \leq n$ such that $x_p > x_q$, then there is an $i$ inserted between $x_p$ and $x_q$ in $[u,I,J]$ \textit{i.e.}
\begin{equation}\label{eqn2.5}
[u,I,J] = x_1 \cdots x_p \cdots i \cdots x_q \cdots x_n.
\end{equation}
\end{lemma}
\begin{proof}
Suppose there is no $i$ inserted between $x_p$ and $x_q$. This means that $s_{p} , s_{p+1} , \cdots , s_{q-1} \in W_I$. Therefore, $u$ is not a minimal representative of $W_I u$, contradicting the assumption.
\end{proof}

For each $a , b \in \Z$ with $a \leq b$, write $[a,b]:=\{ a , a+1 , \cdots , b\}$. As a convention, $[a,b] := \emptyset$ if $a>b$.

\begin{lemma}[$I$-sorting]\label{lem2.7}
Let $G = GL_n$, $W= S_n$, $I,J \subseteq [n-1]$, $\lambda \in \Parj$ and $u \in {}^I W^J$. As in \eqref{eqn1.8}, choose $1 \leq p<q \leq g$ and
\begin{equation}\label{eqn2.6}
K:= [y_{p-1}+1 , y_{q}-1]\subseteq [n-1].
\end{equation}
Let $u'$ be the minimal representative of $W_K u$. Let $w$ and $w'$ be the maximal representatives of $W_I u$ and $W_I u'$, respectively. Then if $H^0 (X(w'W_J),L_\lambda)^*$ is not $L_I$-multiplicity-free, neither is $H^0 (X(wW_J),L_\lambda)^*$. In particular, if $[u,I,J]$ is spherical, so is $[u',I,J]$.
\end{lemma}
\begin{proof}
$u' \leq u$ leads to $w' \leq w$. Due to Lemma \ref{lem2.5}, we have the first statement. From Theorem \ref{lem2.4}, we have the second statement.
\end{proof}

Intuitively, by sorting the numbers between $(p-1)$th $i$ and the $q$th $i$ of $[u,I,J]$ into weakly increasing order, we obtain $[u',I,J]$. For instance, consider $13i4i223i4$ which is a minimal sequence of $7$. Choose $p=2$ and $q=3$. Since it has $I = \{ 1,4,5\}$, $y_1 = 2$ and $y_3 = 6$, we have $K = [3,5]$. Then the underlined numbers between the first $i$ and the third $i$ in $13\underline{i4i223i}4$ are sorted in weakly increasing order and converted into $13i2i234i4$. If $13i4i223i4$ is spherical, so is $13i2i234i4$ due to $I$-sorting.

\begin{lemma}[$J$-sorting]\label{lem2.8}
Let $G = GL_n$, $W= S_n$, $I,J \subseteq [n-1]$, $\lambda \in \Parj$ and $u \in {}^I W^J$. As in \eqref{eqn1.8}, choose $1 \leq p<q \leq h$ and
\begin{equation}\label{eqn2.7}
K:= [ z_{p-1} + 1 ,  z_q - 1 ] \subseteq [n-1].
\end{equation}
Let $u'$ be the minimal representative of $uW_K$. Let $w$ and $w'$ be the maximal representatives of $W_I u$ and $W_I u'$, respectively. Then if $H^0 (X(w'W_J),L_\lambda)^*$ is not $L_I$-multiplicity-free, neither is $H^0 (X(wW_J),L_\lambda)^*$. In particular, if $[u,I,J]$ is spherical, so is $[u',I,J]$.
\end{lemma}
\begin{proof}
$u' \leq u$ leads to $w' \leq w$. Due to Lemma \ref{lem2.5}, we have the first statement. From Theorem \ref{lem2.4}, we have the second statement.
\end{proof}

Intuitively, by sorting the numbers $p , p+1 , \cdots , q$ of $[u,I,J]$ into weakly increasing order, we obtain $[u' ,I,J]$. For instance, consider $13i4i223i4$ which is a minimal sequence of $7$. Choose $p=2$ and $q = 3$. Since it has $J = \{ 2,4,6 \}$, $z_1 = 1$ and $z_3 = 5$, we have $K = [2,4]$. Then the underlined numbers $2$ and $3$ in $1 \underline{3} i 4 i \underline{223} i 4$ are sorted in weakly increasing order and converted into $12i4i233i4$
. If $13i4i223i4$ is spherical, so is $12i4i233i4$ due to $J$-sorting.

\subsection{Expansions}

Let $G = GL_n$ and $W=S_n$. Then $W$ gives a left action on $\R^n$ in a way that for each $w \in W$ and $\alpha = (\alpha_1 , \cdots , \alpha_n) \in \R^n$,
\begin{equation}
w \cdot \alpha := (\alpha_{w^{-1}(1)}, \alpha_{w^{-1}(2)} , \cdots , \alpha_{w^{-1}(n)}).
\end{equation}

Let $\epsilon_1 := (1 , 0 , \cdots , 0), \epsilon_2 := (0 , 1 , \cdots , 0) , \cdots , \epsilon_n := (0 , 0 , \cdots , 1)$ be the standard basis of $\R^n$. Write the fundamental weights $\varpi_1 , \cdots , \varpi_{n-1} \in \Parn$
\begin{equation}
\varpi_{j} := \sum_{i=1}^j \epsilon_i, \quad 1 \leq j \leq n-1.
\end{equation}

Let $I \subseteq [n-1]$, $w \in W$ and $\lambda \in \Parn$. Let $\alpha = w\lambda$. Consider $\kappa_{\alpha} (x_1 , \cdots , x_n)$, the $T$-character of $H^0 (X(w), L_\lambda)^*$. Write $x_1 , \cdots , x_n$ as variables of polynomials. From \eqref{eqn1.8}, write $X_i = \{ x_j \mid y_{i-1}+1 \leq j \leq y_i \}$ for $1 \leq i \leq g$. Key polynomials form a linear basis of polynomial rings; see, 
\textit{e.g.}, \cite[\S 1]{MR4599567}. Therefore, there uniquely exists $c_{\beta^{(1)} , \cdots , \beta^{(g)}}^\alpha \in \Q$ such that
\begin{equation}
\kappa_\alpha (x_1 , \cdots , x_n) = \sum_{\beta^{(i)} \in {\rm Comp}_{y_i - y_{i-1}}} c_{\beta^{(1)} , \cdots , \beta^{(g)}}^\alpha \cdot \kappa_{\beta^{(1)}}(X_1) \cdots \kappa_{\beta^{(g)}}(X_g),
\end{equation}
by separating variables $\{ x_1 , \cdots , x_n \}$ into $X_1 \cup \cdots \cup X_g$.

On the other hand, let $N \in \Z_{\geq 0}$. Let $\xi = \sum_{i \in [n-1]-I} \varpi_i$. From \cite[Theorem~1]{MR1016897}, there exists $d_{\gamma , N}^{\alpha} \in \Z_{\geq 0}$ such that
\begin{equation}\label{eqn3.30}
\kappa_{\alpha}(x_1 , \cdots , x_n) \cdot \kappa_{N \xi} (x_1 , \cdots , x_n) = \sum_{\gamma \in {\rm Comp}_n}  d_{\gamma , N}^{\alpha} \kappa_{\gamma} (x_1 , \cdots , x_n).
\end{equation}
In particular, there is no cancellation between monomials in \eqref{eqn3.30}. Since $\kappa_{N\xi}(x_1 , \cdots , x_n)$ is a monomial, $d_{\gamma , N}^{\alpha} \neq 0$ implies $\kappa_{\gamma}$ is divisible by $\kappa_{N\xi}$, leading to $\gamma - N \xi \in {\rm Comp}_n$. Write
\begin{equation}
{\sf proj}_i : {\rm Comp}_n \rightarrow {\rm Comp}_{y_{i} - y_{i-1}}, \quad \alpha = (\alpha_1, \cdots , \alpha_n) \mapsto (\alpha_{y_{i-1}+1} , \alpha_{y_{i-1}+2} ,\cdots , \alpha_{y_i}).
\end{equation}

\begin{lemma}\label{lem3.4}
Let $I \subseteq [n-1]$. Let $\alpha \in {\rm Comp}_n$, $\xi = \sum_{i \in [n-1]-I} \varpi_i$ and $N \in \Z_{\geq 0}$. If $N$ is large enough, then for any $\gamma \in {\rm Comp}_n$ such that $\gamma - N \xi \in {\rm Comp}_n$,
\begin{equation}
c_{\beta^{(1)} , \cdots , \beta^{(g)}}^\alpha = d_{\gamma , N}^\alpha, \quad {\sf proj}_i (\gamma - N \xi) = \beta^{(i)} , \quad 1 \leq i \leq g.
\end{equation}
\end{lemma}

\begin{proof}
The proof is analogous to the proof of Lemma \ref{lem3.2}. We only need to show that
\begin{equation}\label{eqn3.32}
\kappa_{\beta^{(1)}} (X_1) \cdots \kappa_{\beta^{(g)}}(X_g) \cdot \kappa_{N\xi} (x_1 , \cdots , x_n) = \kappa_{\gamma} (x_1 , \cdots , x_n).
\end{equation}
Let $\gamma = (\gamma_1 , \cdots , \gamma_n)$ and $d_{\gamma , N}^\alpha \neq 0$. If $N$ is large enough, then for any monomial $x_1^{k_1} \cdots x_n^{k_n}$ that appears in \eqref{eqn3.30}, $k_p > k_q$ whenever there exists $1 \leq i \leq g$ such that $p \leq y_i < q$. This leads to $\gamma_p > \gamma_q$. In other words, there exist $\lambda \in \Parn$ and $w \in W_I$ such that $\gamma = w\lambda$.

From a reduced word $(i_1 , \cdots ,i_l)$ of $w$, let $\Delta_w = \Delta_{i_1} \circ \cdots \circ \Delta_{i_l}$ using \eqref{eqn3.31}. From $\Delta_w (x_1^{\lambda_1} \cdots x_n^{\lambda_n}) = \kappa_\gamma (x_1 , \cdots , x_n)$ and $w \in W_I$, we have
\begin{equation}
\kappa_\gamma (x_1 , \cdots , x_n) = \kappa_{{\sf proj}_1 (\gamma) }(X_1) \cdots \kappa_{{\sf proj}_g (\gamma) }(X_g).
\end{equation}
Since $\kappa_\gamma (x_1 , \cdots , x_n)$ is divisible by a monomial $\kappa_{N \xi} (x_1 , \cdots , x_n)$, each $\kappa_{{\sf proj}_i (\gamma) }(X_i)$ is divisible by $\kappa_{{\sf proj}_i (N\xi) }(X_i)$, which is $(x_{y_{i-1}+1} \cdots x_{y_i})^{N(g-i)}$. Therefore, we have \eqref{eqn3.32}. Replacing the right-hand side of \eqref{eqn3.30} with \eqref{eqn3.32}, we have the theorem due to the uniqueness of the expansion.
\end{proof}

Let $w \in W$ and $w = s_{i_1} \cdots s_{i_l}$ be a reduced word expression of $w$. Record it into $\mathbf{i} = (i_1 , \cdots, i_l)$. Let $\lambda =(\lambda_1 , \cdots , \lambda_n) \in \Parn$ and for each $1 \leq j \leq l$, define
\begin{equation}\label{eqn4.3}
m_j =
\begin{cases}
\lambda_k - \lambda_{k+1} & \text{if }i_j = k \text{ and }i_{j'} \neq k \text{ for all }j' > j \\
0 & \text{otherwise}
\end{cases}
.
\end{equation}

In \citep{MR1644050}, LS-monomials of $GL_n$ were defined as $\{ \alpha_{j,k}\}_{1 \leq j \leq l, 1 \leq k \leq m_j}$ associated to $(\mathbf{i}, \mathbf{m})$ where $\alpha_{j,k} \in \Z^n$. The condition imposed on $\alpha_{j,k}$ is that there exists a nested sequence of index sets
\begin{equation}\label{eqn4.4}
[l] \supseteq J_{1,1} \supseteq \cdots \supseteq J_{1, m_1} \supseteq \cdots \supseteq J_{l,1} \supseteq \cdots \supseteq J_{l,m_l},
\end{equation}
where $\mathbf{i}(J_{j,k}\cap[j])$ is a reduced word expression and
\begin{equation}\label{eqn4.5}
w(\mathbf{i}(J_{j,k}\cap [j])) \varpi_{i_j}= \alpha_{j,k}.
\end{equation}
As before, the weight of an LS-monomial $\Pi = \{ \alpha_{j,k}\}_{1 \leq j \leq l, 1 \leq k \leq m_j}$ is defined as
\begin{equation}\label{eqn4.7}
{\sf wt}\left( \Pi \right) = \sum_{1 \leq j \leq l} \left( \sum_{1 \leq k \leq m_j} \alpha_{j,k} \right).
\end{equation}

We say that an LS-monomial $\Pi = \{ \alpha_{j,k} \}_{1 \leq j \leq l, 1 \leq k \leq m_j} $ is \textbf{$I$-dominant} if for any $i \in I$, $1 \leq j \leq l$ and $1 \leq k \leq m_j$, $\beta_i^{(j,k)} \geq \beta_{i+1}^{(j,k)}$ where $\beta^{(j,k)} \in \Z^n$ is defined as
\begin{equation}\label{eqn4.6}
\beta^{(j,k)}:= \sum_{1 \leq j' \leq j-1} \left( \sum_{1 \leq k' \leq m_{j'}} \alpha_{j',k'} \right)+ \sum_{1 \leq k' \leq k} \alpha_{j,k'}.
\end{equation}

\begin{theorem}\label{thm4.1}
Let $G=GL_n$, $W=S_n$, $I , J \subseteq [n-1]$ and $v \in W$. Let $w$ be the maximal representative of $W_I v W_J$. Choose $( \mathbf{i}, \mathbf{m})$ from a reduced word expression of $v$ and $\lambda \in \Parj$. Then $H^0(X(wW_J),L_\lambda)^*$ is $L_I$-multiplicity-free if and only if there are no two distinct $I$-dominant LS-monomials of $(\mathbf{i}, \mathbf{m})$ which have the same weight.
\end{theorem}
\begin{proof}
Choose $(\mathbf{i'}, \mathbf{m'})$ from a reduced word expression of $w$ and $\lambda$. Let $\xi := \sum_{i \in [n-1]- I} \varpi_i$ and $N \in \Z_{\geq 0}$ be large enough. $H^0(X(wW_J) , L_\lambda)^*$ is $L_I$-multiplicity-free if and only if $c_{\beta^{(1)} , \cdots , \beta^{(g)}}^{w \lambda}$ is either $1$ or $0$. From Lemma \ref{lem3.4}, this is equivalent to the condition that $d_{\gamma , N}^{w\lambda}$ is either $1$ or $0$ in the expansion of $\kappa_{w\lambda} \cdot \kappa_{N \xi}$. Applying Theorem \ref{thm3.3} to LS-monomials of $GL_n$ in \citep{MR1644050}, this is also equivalent to the following statement: there are no two distinct $I$-dominant LS-monomials of $(\mathbf{i'}, \mathbf{m'})$ of $w$ and $\lambda$ which have the same weight. Indeed, if $N$ is large enough, then a LS-monomial $\Pi$ is $I$-dominant if and only if $\pi^{N\xi} * \Pi$ is dominant.

Write $w_0$ as the longest element in $W_I$. Choose a reduced word of $w_0$ and construct a Demazure operator $\Delta_{w_0}$. Applying $\Delta_{w_0}$ to $\kappa_{v \lambda} \cdot \kappa_{N \xi}$, we get $\kappa_{w \lambda} \cdot \kappa_{N \xi}$ due to $\Delta_{w_0} (\kappa_{N\xi}) = \kappa_{N \xi}$. If $\kappa_{w' \lambda'}$ appears in $\kappa_{v \lambda} \cdot \kappa_{N\xi}$ with $\lambda' \in \Parn$, then $w'$ may be chosen in $W_I$ assuming $N$ is large enough. Then $\kappa_{w_0 \lambda'}$ appears in $\kappa_{w \lambda} \cdot \kappa_{N\xi}$. There exist two distinct $I$-dominant LS-monomials of $(\mathbf{i} , \mathbf{m})$ with the same weight if and only if there exist $w_{1}' , w_2' \in W_I$ and $\lambda'$ such that $\kappa_{w_1' \lambda'} + \kappa_{w_2' \lambda'}$ appears in the expansion of $\kappa_{v \lambda} \cdot \kappa_{N\xi}$. Equivalently, $\kappa_{w_0 \lambda'} + \kappa_{w_0 \lambda'}$ appears in the expansion of $\kappa_{w \lambda} \cdot \kappa_{N \xi}$. Therefore, there exist two distinct $I$-dominant LS-monomials of $(\mathbf{i}, \mathbf{m})$ with the same weight if and only if the same condition holds for $(\mathbf{i'}, \mathbf{m'})$.
\end{proof}

\begin{lemma}[$I$-expansion]\label{lem4.4}
Let $G = GL_n$, $W=S_n$, $I \subseteq I' \subseteq [n-1]$, $J \subseteq [n-1]$ and $u \in {}^I W^J$. Let $u'$ be the minimal representative of $W_{I'} u W_J$. Let $w$ and $w'$ be the maximal representatives of $W_I u W_J$ and $W_{I'} u' W_J$, respectively. Then if $X(wW_J)$ is $L_I$-spherical, $X(w'W_J)$ is $L_{I'}$-spherical. In other words, if $[u,I,J]$ is spherical, so is $[u',I',J]$.
\end{lemma}

\begin{proof}
Suppose $X(w'W_J)$ is not $L_{I'}$-spherical. Then due to Theorem \ref{lem2.4}, there exists a $\lambda \in \Parj$ such that $H^0(X(w'W_J),L_\lambda)^*$ is not $L_{I'}$-multiplicity-free. Choose a reduced expression of $u'$ 
\begin{equation}\label{eqn4.34}
u' = s_{i_1'} s_{i_2'} \cdots s_{i_{l'}'},
\end{equation}
and record it into $(\mathbf{i'}, \mathbf{m'})$ from $\lambda$. Due to Theorem \ref{thm4.1}, there exist $I'$-dominant LS-monomials $\{\alpha_{j,k}'\}$ and $\{\beta_{j,k}'\}$ with the same weight. Write the corresponding chains of indices $\{J_{j,k}'\}$ and $\{K_{j,k}'\}$.

Since $u \in W_{I'}u'$, we may choose a reduced expression of $u$
\begin{equation}\label{eqn4.35}
u = s_{i_1} s_{i_2} \cdots s_{i_l} s_{i_1'} s_{i_2 ' } \cdots s_{i_{l'}'}. 
\end{equation}
By computing $\mathbf{m}$ from $\mathbf{i}$ and $\lambda$,
\begin{equation}\label{eqn4.36}
m_{j+l} = m_j' \quad j \in [l'].
\end{equation}
For each index set $J = \{ j_1 , \cdots , j_m\}$ and $a \in \Z$, denote $a + J := \{ a+ j_1 , \cdots , a+ j_m\}$. Define
\begin{equation}\label{eqb4.37}
J_{j,k} = 
\begin{cases}
[l+1 , l+l'] & j \in [l] \\
l+J_{j-l,k}' & j \in [l+1 , l+l']
\end{cases}
, \quad K_{j,k} = 
\begin{cases}
[l+1 , l+l'] & j \in [l] \\
l+K_{j-l,k}' & j \in [l+1 , l+l']
\end{cases}
.
\end{equation}
From $\{J_{j,k}\}$ and $\{K_{j,k}\}$, we have LS-monomials $\{\alpha_{j,k}\}$ and $\{\beta_{j,k}\}$, respectively. Due to $I \subseteq I'$, they are $I$-dominant. Therefore, $H^0 (X(wW_J),L_\lambda)^*$ is not $L_I$-multiplicity-free, proving the given statement.
\end{proof}

Intuitively, suppose we have $1 \leq p \leq g-1$ such that $y_p \in  (I'-I)$ using notation in \eqref{eqn1.8}. Then $I$-expansion eliminates the $p$th $i$ in $[u,I,J]$. Then some numbers are reordered so that numbers between successive occurrences of $i$ are still weakly increasing. For instance, consider $13i4 i223i4$ which is a minimal sequence of $7$. From $I = \{1,4,5 \}$, choose $I' := \{ 1, 3, 4, 5 \}$. Then the second $i$ is removed; we have $13i\underline{4223}i4$. Reorder the underlined numbers into weakly increasing order and we have $13i2234i4$. If $13i4i223i4$ is spherical, so is $13i2234i4$, due to $I$-expansion.

\begin{lemma}[$J$-expansion]\label{lem4.5}
Let $G = GL_n$, $W=S_n$, $I \subseteq [n-1]$, $J \subseteq J' \subseteq [n-1]$ and $u \in {}^I W^J$. Let $u'$ be the minimal representative of $W_I u W_{J'}$. Let $w$ and $w'$ be the maximal representatives of $W_I u W_J$ and $W_I u' W_{J'}$, respectively. Then if $X(wW_{J})$ is $L_I$-spherical, $X(wW_{J'})$ is $L_I$-spherical. In other words, if $[u,I,J]$ is spherical, so is $[u',I , J']$.
\end{lemma}
\begin{proof}
Suppose $X(wW_{J'})$ is not $L_I$-spherical. Then due to Theorem \ref{lem2.4}, there exists a $\lambda \in {\rm Par}_{J'}$ such that $H^0 (X(wW_{J'}), L_\lambda)^*$ is not $L_I$-multiplicity-free.

Let $u'$ be the minimal representative of $W_I w W_{J'}$. Choose a reduced expression of $u'$
\begin{equation}
u' = s_{i_1'} s_{i_2'} \cdots s_{i_{l'}'},
\end{equation}
and record it into $(\mathbf{i'}, \mathbf{m}')$ from $\lambda$. Due to Theorem \ref{thm4.1}, there exist $I$-dominant LS-monomials $\{\alpha_{j,k}'\}$ and $\{\beta_{j,k}'\}$ with the same weight. Write the corresponding chains of indices $\{J_{j,k}'\}$ and $\{K_{j,k}'\}$.

Let $u$ be the minimal representative of $W_I w W_J$. Since $J \subseteq J'$, we may choose a reduced expression of $u$
\begin{equation}
s_{i_1 ' }s_{i_2 ' } \cdots s_{i_{l'}'} s_{i_1} s_{i_2} \cdots s_{i_l}.
\end{equation}
Due to $J \subseteq J'$, $i_k \in J'$ for $1 \leq k \leq l$. Record it into $(\mathbf{i}, \mathbf{m})$ from $\lambda$ which is also in $\Parj$. Then
\begin{equation}
m_j = 
\begin{cases}
m_j ' & j \in [l'] \\
0 & j \in [l'+1 , l' + l]
\end{cases}
.
\end{equation}
Define
\begin{equation}
J_{j,k}:= J_{j,k}' \quad K_{j,k}:= K_{j,k}'.
\end{equation}
From $\{J_{j,k}\}$ and $\{K_{j,k}\}$, we have $LS$-monomials $\{\alpha_{j,k}\}$ and $\{\beta_{j,k}\}$, respectively. Also, they are $I$-dominant. Therefore, $H^0(X(wW_J), L_\lambda)^*$ is not $L_I$-multiplicity-free, proving the given statement.
\end{proof}

Intuitively, suppose we have $1 \leq p \leq h-1$ such that $z_p \in (J' - J)$ using notation in \eqref{eqn1.8}. Then $J$-expansion equalizes the numbers $p$ and $p+1$ into $p$. The remaining numbers $p+2, p+3 , \cdots , h$ are shifted to $p+1 , p+2 , \cdots , h-1$. For instance, consider $13i4i223i4$ which is a minimal sequence of $7$. From $J = \{2,4,6\}$, choose $J' := \{2,3,4,6 \}$. Then the underlined numbers in $1\underline{3}i4i\underline{223}i4$ are equalized into the same number and we have $12i3i222i3$. If $13i4i223i4$ is spherical, so is $12i3i222i3$ due to $J$-expansion.

As a convention, we assume that an empty product of simple reflections is the identity.

\begin{lemma}\label{lem3.3}
Let $G= GL_n$, $W= S_n$, $I , J \subseteq [n-1]$, $u \in {}^I W^J$ with a reduced word expression $u = s_{i_1} \cdots s_{i_l}$. Assume there exists a $p \in [n-1]$ such that for any $j \in [l]$, $i_j \neq  p$. Let
\begin{equation}
u' := \prod_{i_j \leq p-1}s_{i_j} ,  \quad u'' := \prod_{i_j \geq p+1} s_{i_j}.
\end{equation}
\begin{itemize}
\item[(a)] $[u,I,J]$ avoids patterns \eqref{eqn1.7} if and only if $[u',I,J]$ and $[u'',I,J]$ avoid patterns \eqref{eqn1.7}.
\item[(b)] $[u,I,J]$ is spherical if and only if $[u' , I, J]$ and $[u'',I,J]$ are spherical.
\end{itemize}
\end{lemma}
\begin{proof}
``$\Rightarrow$ of (a)'': Suppose a pattern appears in $[u',I,J]$ or $[u'',I,J]$. Let $[u',J] = x_1' \cdots x_n'$ and $[u'',J] = x_1'' \cdots x_n''$. The patterns in \eqref{eqn1.7} have the following common property. The leftmost number entry is not the minimal number and the rightmost number entry is not the maximal number. Therefore, whenever a pattern appears, it should be contained in $x_1' \cdots x_p'$ or $x_{p+1}'' \cdots x_n''$ with $i$ insertions. Then it appears at the same position of $[u,I,J]$.

``$\Leftarrow$ of (a)'': Let $[u,J] = x_1 \cdots x_n$. Then for all $i \leq p < p+1 \leq j$, $x_i \leq x_j$. Suppose $[u,I,J]$ has one of the patterns \eqref{eqn1.7}, of which numerical entries are $x_{j_1} x_{j_2} \cdots x_{j_k}$. Our claim is that $j_k \leq p$ or $j_1 \geq p+1$. Otherwise, $j_1 \leq p < p+1 \leq  j_k$. Since $x_{j_1} \leq x_{j_k}$, the pattern should be either $233112$ or $244113$. However, $x_{j_2} > x_{j_k}$ and $x_{j_5} < x_{j_1}$ imply $j_2 \geq p+1$ and $j_5 \leq p$, leading to a contradiction. From our claim, $[u',I,J]$ has the same pattern if $j_k \leq p$. If $j_1\geq p+1$, $[u'',I,J]$ has the same pattern.

``$\Rightarrow$ of (b)'': Without loss of generality, assume that $[u',I,J]$ is non-spherical. Let $u' = s_{i_1'} \cdots s_{i_{l'}'}$ be a reduced word expression. Record it into $\mathbf{i'} = (i_1' , \cdots , i_{l'}')$. Due to Theorem \ref{thm4.1}, there exist a $\lambda \in \Parj$, a corresponding $\mathbf{m'}$ and two distinct $I$-dominant LS-monomials $\{\alpha_{j,k}' \}$ and $\{\beta_{j,k}' \}$ with the same weight. Let $\{ J_{j,k}' \}$ and $\{K_{j,k}'\}$ be the corresponding LS-chains.

On the other hand, choose a reduced word expression of $u = s_{i_1 '} \cdots s_{i_{l'}'} s_{i_1 ''} \cdots s_{i_{l''}''}$. Record it into $\mathbf{i} = (i_1' , \cdots , i_{l'}', i_{1}'', \cdots , i_{l''}'')$. Let $\mathbf{m}$ be defined from $\mathbf{i}$ and $\lambda$. For $1 \leq j \leq l'+l''$ and $1 \leq k \leq m_j$, define
\begin{equation}
J_{j,k} :=
\begin{cases}
J_{j,k}' & j \in [l'] \\
\emptyset & j \in [l'+1 , l'+l'']
\end{cases}
, \quad K_{j,k}:=
\begin{cases}
K_{j,k}' & j \in [l'] \\
\emptyset & j \in [l'+1 , l'+l'']
\end{cases}
.
\end{equation}
From LS-chains $\{J_{j,k} \}$ and $\{K_{j,k}\}$, we have LS-monomials $\{ \alpha_{j,k} \}$ and $\{\beta_{j,k}\}$, which are the same as $\{\alpha_{j,k}' \}$ and $\{\beta_{j,k}' \}$ with some $\varpi_{i_j}$ added. This proves that $[u,I,J]$ is non-spherical.

``$\Leftarrow$ of (b)'': Choose a reduced word expression of $u  = s_{i_1 '} \cdots s_{i_{l'}'} s_{i_1 ''} \cdots s_{i_{l''}''}$ so that $i_j ' \leq p-1$ and $i_j'' \geq p+1$. Record this into $\mathbf{i}$. Suppose $[u,I,J]$ is non-spherical. Then there exist a $\lambda \in \Parj$, $\mathbf{m}$ and $I$-dominant LS-monomials $\{\alpha_{j,k}\}$ and $\{\beta_{j,k}\}$ with the same weight due to Theorem \ref{thm4.1}. Since they are distinct, there exists a $q \in [n]$ such that the $q$th components of some $\alpha_{j,k}$ and $\beta_{j,k}$ are different. Then $q \leq p$ or $q \geq p+1$.

If $q \leq p$, then choose a reduced word expression of $u' = s_{i_1'} \cdots s_{i_{l'}'}$ with $\mathbf{m'}$ from $\lambda$. For each $1 \leq j \leq l'$ and $1 \leq k \leq m_j'$, define
\begin{equation}
J_{j,k}' := J_{j,k} \cap [l'], \quad K_{j,k}':= K_{j,k}\cap [l'].
\end{equation}
Then the corresponding LS-monomials $\{ \alpha_{j,k}' \}$ and $\{ \beta_{j,k}'\}$ are distinct, $I$-dominant with the same weight.

If $q \geq p+1$, then choose a reduced word expression of $u'' = s_{i_1''} \cdots s_{i_{l''}''}$ with $\mathbf{m''}$ from $\lambda$. For each $1 \leq j \leq l''$ and $1 \leq k \leq m_j''$, define
\begin{equation}
J_{j,k}'' := (-l')+(J_{j+l',k} \cap [l'+1 , l'+l'']) , \quad K_{j,k}'' := (-l')+(K_{j+l',k} \cap [l'+1 , l'+l'']) 
\end{equation}
Then the corresponding LS-monomials $\{ \alpha_{j,k}'' \}$ and $\{ \beta_{j,k}''\}$ are distinct, $I$-dominant with the same weight.
\end{proof}

\subsection{Line bundles criterion}\label{sub4.2}

We prove \citep[Conjecture 4.19]{MR4394321} and give another criterion for $L_I$-spherical $X(w)$ in $GL_n / B$.

\begin{lemma}\label{lem4.7}
Let $G= GL_n$ and $W= S_n$. Let $v \in W$, $I \subseteq [n-1]$ and $J=\emptyset$. Write $u$ and $w$ as the minimal and maximal representatives of $W_I v$, respectively. Then if $[u,I,J]$ has a $3412$ pattern, $H^0(X(w),L_\lambda)^*$ is not $L_I$-multiplicity-free for any strictly decreasing partition $\lambda$.
\end{lemma}
\begin{proof}
Write $[u,J] = x_1 \cdots x_n$. Let $a < b < c < d$ such that $x_c < x_d < x_a < x_b$. Using Lemma \ref{lem2.6}, there exists a $y \in [n-1]-I$ such that $b \leq y < c$. Use $I$-sorting Lemma \ref{lem2.7} twice on $K = [1,y-1]$ and $K = [y+1,n-1]$. Then we have $u'$ and $w'$ which are the minimal and maximal representatives of $W_Iu' = W_I w'$ with the properties described in Lemma \ref{lem2.7}. Denote $[u',J]= x_1' \cdots x_n '$. We still have $1 \leq a' < b' \leq y < c' < d' \leq n$ such that $x_{a'}' = x_a$, $x_{b'}' = x_b$, $x_{c'}' = x_c$ and $x_{d'}' = x_d$. Therefore, $[u',J]$ also has a $3412$ pattern. We may assume
\begin{equation}\label{eqn4.10}
x_1 < x_2 < \cdots < x_{y}, \quad x_{y+1} < x_{y+2} < \cdots < x_n, \quad y=a+1 = b = c-1 = d-2
\end{equation}
without loss of generality. Our claim is that it is sufficient to prove the case when $x_b - x_c = 3$.

(Case 1) $x_d - x_c >1$ : Write $x_c = p$ and $x_d = q$. Choose
\begin{equation}\label{eqn4.9}
K:= [ p ,  q-2 ] \subseteq [n-1].
\end{equation}
Apply $J$-sorting Lemma \ref{lem2.8} from $K$ to obtain $u'$ and $w'$. Then $p , p+1 , \cdots , q-1$ parts of $[u,J]$ are rearranged into increasing order which is identical to $[u',J]$. Replace $u$ and $w$ with $u'$ and $w'$ and we get $x_d - x_c $ reduced by increasing $x_c$.

(Case 2) $x_b - x_a>1$ : Analogous to (Case 1).

(Case 3) $x_a - x_d >1$ : Write $x_a = q$ and $x_d = p$. Use $J$-sorting Lemma \ref{lem2.8} twice by setting $K = [p, q-2]$ and $K = [p+1 , q-1]$. In other words, apply $J$-sorting on $p, \cdots , q-1$ and then $p+1, \cdots , q$ to obtain $u'$ and $w'$. Denote $[u',J] = x_1' \cdots x_n'$. Due to the previously applied $J$-sortings, $x_a' \leq x_a$ and $x_d' \geq x_d$. Our claim is that either $x_a ' < x_a$ or $x_d' > x_d$. If neither of them happens, then for any $j \in [n]$ with $p < x_j <q$, $j >d$ and $j <a$, leading to a contradiction. Therefore, $x_a - x_d$ can be reduced.

Our conclusion is that we may assume
\begin{equation}\label{eqn4.12}
x_b - x_a = x_a - x_d = x_d - x_c = 1
\end{equation}
without loss of generality. Combined with \eqref{eqn4.10}, $u = s_y s_{y-1}s_{y+1} s_y$. Let $\mathbf{i} = (y , y-1 , y+1 , y)$. Then $\mathbf{m} = (m_1 , m_2 , m_3 , m_4)$ is computed as
\begin{equation}\label{eqn4.37}
m_1 = 0, \quad m_2 = \lambda_{y-1}-\lambda_y , \quad m_3 = \lambda_{y+1} - \lambda_{y+2} , \quad m_4 = \lambda_{y} - \lambda_{y+1}.
\end{equation}
Since $\lambda$ is strictly decreasing, $m_2 , m_3 , m_4 \geq 1$. Using this, define $J_{j,k}$ and $K_{j,k}$
\begin{equation}\label{eqn4.38}
J_{j,k}=
\begin{cases}
[4] & 1 \leq j \leq 2 \\
\{1,3,4\} & j=3 , \hsp k=1 \\
\emptyset & j=3 , \hsp k>1 \\
\emptyset & j=4 \hsp k=1 \\
\emptyset & j=4 , \hsp k>1
\end{cases}
, \quad K_{j,k} = 
\begin{cases}
[4] & 1 \leq j \leq 2 \\
[2,4] & j=3 , \hsp k=1 \\
\{2,4\} & j=3 \hsp k>1 \\
\{4\} & j=4 , \hsp k=1 \\
\emptyset & j=4 , \hsp k>1
\end{cases}
.
\end{equation}
Then the corresponding $\alpha_{j,k}$ and $\beta_{j,k}$ are
\begin{subequations}\label{eqn4.41}
\begin{equation}\label{eqn4.39}
\alpha_{j,k}=
\begin{cases}
\varpi_{y-1}-\epsilon_{y-1}+\epsilon_{y+1} & j=2 \\
\varpi_{y+1}-\epsilon_y+\epsilon_{y+2} & j=3 , \hsp k=1 \\
\varpi_{y+1} & j=3 , \hsp k>1 \\
\varpi_{y} & j=4 \hsp k=1 \\
\varpi_y & j=4 , \hsp k>1
\end{cases}
\end{equation}
\begin{equation}\label{eqn4.40}
\beta_{j,k} = 
\begin{cases}
\varpi_{y-1}-\epsilon_{y-1}+\epsilon_{y+1} & j=2 \\
\varpi_{y+1}-\epsilon_{y+1}+\epsilon_{y+2} & j=3 , \hsp k=1 \\
\varpi_{y+1} & j=3 \hsp k>1 \\
\varpi_{y}-\epsilon_{y}+\epsilon_{y+1} & j=4 , \hsp k=1 \\
\varpi_y & j=4 , \hsp k>1
\end{cases}
.
\end{equation}
\end{subequations}
$\{\alpha_{j,k}\}$ and $\{\beta_{j,k}\}$ are $I$-dominant LS-monomials with the same weight. Due to Theorem \ref{thm4.1}, this proves $H^0(X(w) , L_\lambda)^*$ is not $L_I$-multiplicity-free.
\end{proof}

\begin{lemma}\label{lem4.8}
Let $G= GL_n$ and $W= S_n$. Let $v \in W$, $I \subseteq [n-1]$ and $J=\emptyset$. Write $u$ and $w$ as the minimal and maximal representatives of $W_I v$, respectively. Then if $[u,I,J]$ has a $321$ pattern, $H^0(X(w),L_\lambda)^*$ is not $L_I$-multiplicity-free for any strictly decreasing partition $\lambda$.
\end{lemma}
\begin{proof}
Write $[u,J] = x_1 \cdots x_n$. Since it has a $321$ pattern, there exist $1 \leq a < b < c \leq n$ such that $x_c < x_b < x_a$. Due to Lemma \ref{lem2.6}, there exist $y , y' \in [n-1]-I$ such that $a \leq y < b$ and $b \leq y' < c$. Use $I$-sorting Lemma \ref{lem2.7} three times on $K = [1,y-1]$, $K=[y+1 , y'-1] $ and $K=[y'+1 , n-1]$. This process rearranges $x_1 , \cdots , x_{y}$ and $x_{y+1} , \cdots , x_{y'}$ and  $x_{y'+1} , \cdots , x_n$ so that they are in increasing order while still maintaining a $321$ pattern. Therefore, we may assume that
\begin{equation}\label{eqn4.27}
x_1 < x_2 < \cdots < x_y , \quad x_{y+1} < x_{y+2} < \cdots < x_{y'} , \quad x_{y'+1} < x_{y'+2} < \cdots < x_n,
\end{equation}
\begin{equation}\label{eqn4.28}
a = y , \quad y+1 \leq b \leq y' ,\quad c=y'+1.
\end{equation}

Our claim is that we may assume that $x_a - x_c = 2$. Suppose $x_a - x_c >2$. Without loss of generality, assume that $x_a - x_b > 1$. Write $p:=x_b$ and $q:=x_a$. Choose $K := [p +1, q -1]$. Apply $J$-sorting Lemma \ref{lem2.8} and rearrange $p+1 , p+2 , \cdots , q$ which appear in $[u,J]=x_1 \cdots x_n$ into increasing order $x_1' \cdots x_n'$. If $x_a' = q$, then there exists $f<a$ such that $x_f' = p+1$ due to $J$-sorting. Then $f<a<b<c$ and $x_c' < x_b' < x_f' < x_a'$, meaning that there is a $3421$ pattern. Using $I$-sorting to $K = [y+1 , n-1]$, we get a $3412$ pattern. Using Lemma \ref{lem4.7}, this shows us that $H^0(X(w) , L_\lambda)$ is not $L_I$-multiplicity-free. On the other hand, if $x_a' < q$, then we reduce $x_a - x_c$. Therefore, we may assume that
\begin{equation}\label{eqn4.29}
x_a = p+1 , \quad x_b = p , \quad x_c = p-1,
\end{equation}
without loss of generality. Lastly, apply $J$-sorting Lemma \ref{lem2.8} on $K = [1, p-3]$ and $K=[p+2 , n-1]$ and rearrange $1,2, \cdots , p-2$ and $p+2 , p+3 , \cdots , n$ in $[u,I,J]$ in increasing order. In particular, $p+2 , p+3, \cdots $ should appear after $p$ whereas $\cdots , p-3 , p-2$ appear before $p$ due to \eqref{eqn4.27}. Therefore, $p = b$ and
\begin{equation}\label{eqn4.42}
u = s_a s_{a+1} \cdots s_{b-2}s_{b-1} s_{c-1} s_{c-2} \cdots s_{b}s_{b-1}.
\end{equation}
Let $l:= \ell (u)$ and define $\mathbf{i}$ from \eqref{eqn4.42}. Then from $\lambda$, compute $\mathbf{m}$
\begin{equation}\label{eqn4.43}
m_j = 
\begin{cases}
\lambda_{j+a-1}- \lambda_{j+a} & j \in [b-a-1] \\
0 & j =b-a \\
\lambda_{c+b-a-j} - \lambda_{c+b-a-j+1} & j \in [b-a+1 , l]
\end{cases}
.
\end{equation}
Since $\lambda$ is strictly decreasing, $m_j \geq 1$ except for $j = b-a$. Define
\begin{equation}\label{eqn4.44}
J_{j,k} = 
\begin{cases}
[l] & j \in [l-2] \\
[l] & j =l-1, \hsp k=1 \\
\emptyset & j = l-1 , \hsp k>1 \\
\emptyset & j = l, \hsp k=1 \\
\emptyset & j=l, \hsp k>1
\end{cases}
, \quad K_{j,k} = 
\begin{cases}
[l] & j \in [l-2] \\
[l]-\{ b-a \} & j =l-1, \hsp k=1 \\
[b-a-1] \cup \{l\} & j = l-1 , \hsp k>1 \\
[b-a-1] \cup \{l\} & j = l, \hsp k=1 \\
\emptyset & j=l, \hsp k>1
\end{cases}
.
\end{equation}
Then the corresponding $\alpha_{j,k}$ and $\beta_{j,k}$ are
\begin{subequations}\label{eqn4.45}
\begin{equation}\label{eqn4.46}
\alpha_{j,k} = 
\begin{cases}
\varpi_{j+a-1} +\epsilon_{j+a}-\epsilon_a & j \in [b-a-1] \\
\varpi_{c+b-a-j}-\epsilon_{c+b-a-j}+\epsilon_c & j \in [b-a+1 , l-2] \\
\varpi_b -\epsilon_a +\epsilon_c & j =l-1, \hsp k=1 \\
\varpi_b & j = l-1 , \hsp k>1 \\
\varpi_{b-1} & j = l, \hsp k=1 \\
\varpi_{b-1} & j=l, \hsp k>1
\end{cases}
,
\end{equation}
\begin{equation}\label{eqn4.47}
\beta_{j,k}=
\begin{cases}
\varpi_{j+a-1} +\epsilon_{j+a}-\epsilon_a & j \in [b-a-1] \\
\varpi_{c+b-a-j}-\epsilon_{c+b-a-j}+\epsilon_c & j \in [b-a+1 , l-2] \\
\varpi_b - \epsilon_b + \epsilon_c & j = l-1 , \hsp k=1 \\
\varpi_b & j = l-1 , \hsp k>1 \\
\varpi_{b-1}- \epsilon_a + \epsilon_{b} & j = l , \hsp k=1 \\
\varpi_{b-1} & j = l , \hsp k>1
\end{cases}
.
\end{equation}
\end{subequations}
$\{\alpha_{j,k}\}$ and $\{\beta_{j,k}\}$ are $I$-dominant LS-monomials with the same weight. Due to Theorem \ref{thm4.1}, this proves $H^0(X(w) , L_\lambda)^*$ is not $L_I$-multiplicity-free.
\end{proof}

\begin{proof}[Proof of Theorem \ref{thm1.1}]
Suppose $X(w)$ is $L_I$-spherical. Then from Theorem \ref{lem2.4}, $H^0(X(w),L_\lambda)$ is $L_I$-multiplicity-free.

Conversely, suppose $H^0(X(w),L_\lambda)$ is $L_I$-multiplicity-free. Let $u$ be the minimal representative of $W_I w$ and $J =\emptyset$. Then due to Lemma \ref{lem4.7} and \ref{lem4.8}, $[u,I,J]$ does not have $3412$ or $321$ patterns. Due to Theorem \ref{thm2.3}, $X(w)$ is $L_I$-spherical.
\end{proof}

\begin{remark}\label{rmk3.1}
According to \cite[Theorem~2.1.2]{MR3177371}, a normal quasi-projective $G$-variety $X$ is $G$-spherical if and only if for any $G$-linearized line bundle $L$, $H^0 (X , L)$ is $G$-multiplicity-free. We improved this criterion when $X = X(w)$ and $G = L_I$; it is sufficient to check multiplicity for a single line bundle $L = L_\lambda$ where $\lambda$ is strictly decreasing. Geometrically, $L_\lambda$ over $G/B$ is very ample if and only if $\lambda$ is strictly decreasing. However, an analogous statement is not true for other Lie types or in partial flag varieties.

For instance, let $n=4$, $I = \emptyset$, $J = \{1,3\}$ and $u = 3412$. Then $w = 4321$ is the maximal representative of $W_I u W_J$. Let $\lambda = (1,1,0,0) \in \Parj$. Even though $H^0 (X(wW_J),L_\lambda)$ is $L_I$-multiplicity-free, $H^0 (X(wW_J),L_{2\lambda})$ is not $L_I$-multiplicity-free. In other words, $\kappa_{0,0,1,1}$ is monomial multiplicity-free whereas $\kappa_{0,0,2,2}$ is not; see \cite[Theorem 1.1]{MR4599567}. Therefore, an analog of Theorem \ref{thm1.1} does not hold in partial flag varieties.
\end{remark}

\section{Sphericality of Coxeter elements}\label{sec5}

\begin{figure}
\centering
\begin{subfigure}[t]{0.45\textwidth}
\centering
	\includegraphics[width=\textwidth]{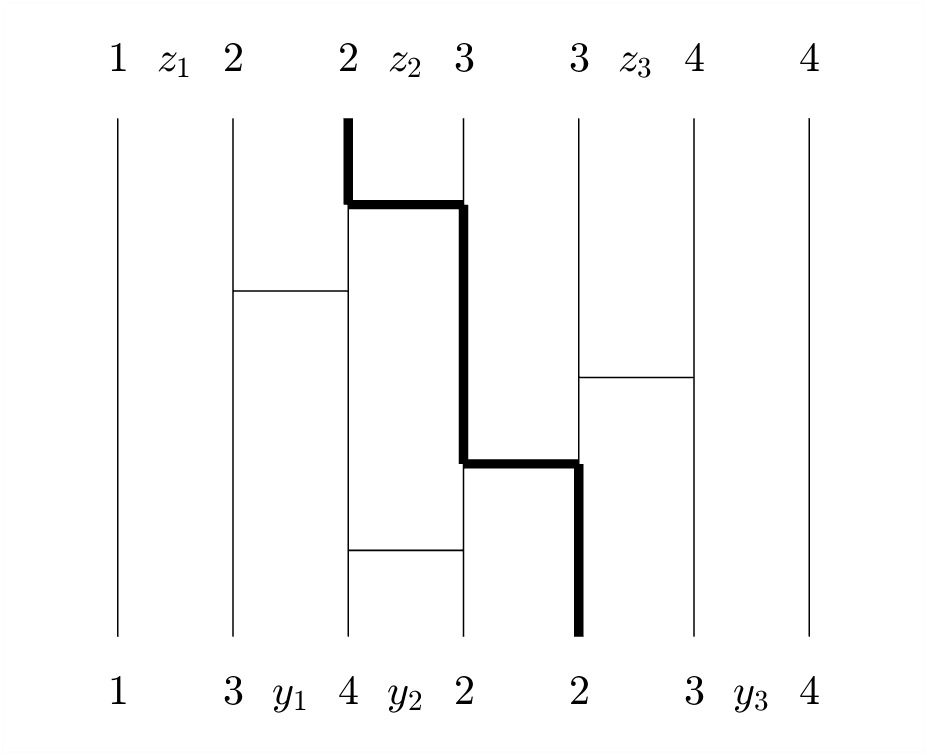}
	\caption{$\{ y_1 , y_2 , y_3 \}$ and $\{z_1 , z_2 , z_3 \}$}
	\label{fig0}
\end{subfigure}
\hfill
\begin{subfigure}[t]{0.45\textwidth}
\centering
\includegraphics[width=\textwidth]{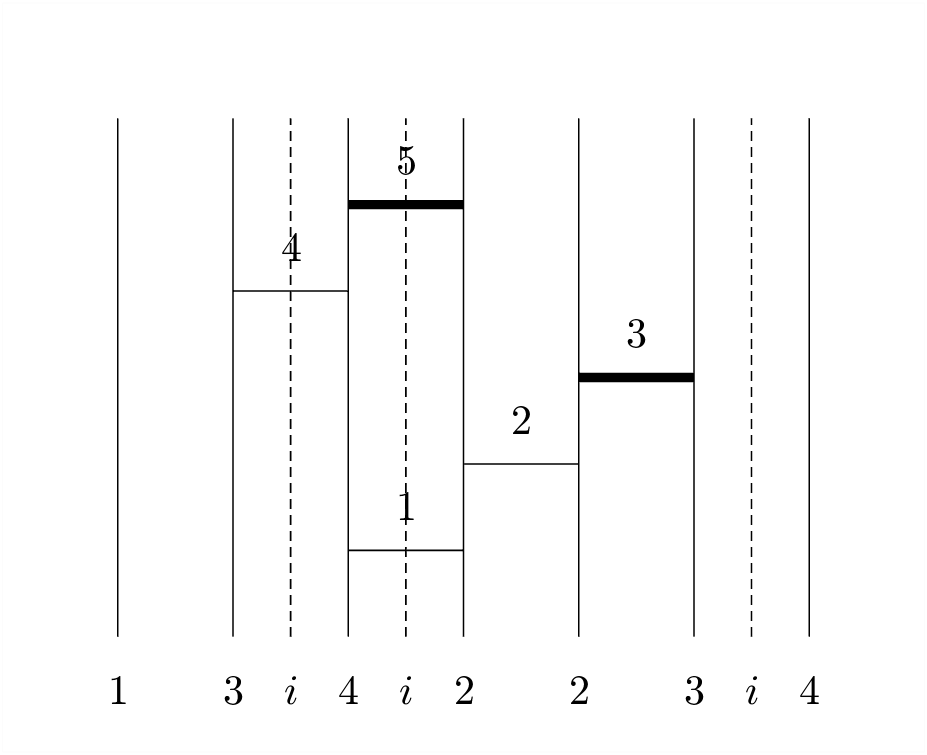}
	\caption{$(i_1, i_2 , i_3 , i_4 , i_5) = (3,4,5,2,3)$}
	\label{fig00}
\end{subfigure}
\caption{Wiring diagram of $[u,I,J] = 13i4i223i4$}
\label{fig000}
\end{figure}

To compute LS-monomials, we use \textit{wiring diagrams} of permutations, also known as \textit{ghost legs} or \textit{ladder games}. For instance, consider $13i4i223i4$ which is a minimal sequence of $7$. Since there are $i$ insertions at $2,3,6$, we have $[6] - I = \{2,3,6\}$, leading to $y_1 = 2$, $y_2 = 3$ and $y_3 = 6$. On the other hand, by reordering number entries into weakly increasing order, we have $1223344$. Since the numbers are strictly increasing at $1,3,5$, we have $[6]-J = \{1,3,5\}$, leading to $z_1 = 1$, $z_2 = 3$ and $z_3 = 5$. See the numbers at the top and bottom of Figure \ref{fig0}.

By choosing $u = s_3 s_4 s_5 s_2 s_3$, $[u,I,J] = 13i4i223i4$. We can verify this by using the wiring diagram in Figure \ref{fig0}; each horizontal line in the $j$th column represents a simple transposition $s_j$. While going downward along a vertical line, cross to another vertical line whenever we meet a horizontal line. For instance, the bold lines in Figure \ref{fig0} depict the trajectory of the second $2$ on the top.

Choose $\lambda = (4,3,3,2,2,0,0) \in \Parj$. From $\mathbf{i} = (3,4,5,2,3)$, we have $\mathbf{m} = (0,0,2,0,1)$. We denote this by using bold lines for the third and the fifth horizontal lines in Figure \ref{fig00}. By setting $J_{3,1} = \{1,2,3,4,5\}$, $J_{3,2} = \{1,2,4,5\}$ and $J_{5,1} = \{ 4 ,5 \}$, we have $\alpha_{3,1} = (1,1,0,1,1,1,0)$, $\alpha_{3,2} = (1,1,1,1,1,0,0)$ and $\alpha_{5,1} = (1,0,1,1,0,0,0)$. We draw dotted lines to denote the $i$ insertions as in Figure \ref{fig00}. From there, $\{\alpha_{3,1} , \alpha_{3,2}, \alpha_{5,1}\}$ is $I$-dominant with weight $(3,2,2,3,2,1,0)$.

\subsection{Sphericality from complete flags}

In this subsection, we discuss $321$ and $3412$ patterns, observed in Theorem \ref{thm2.3}.

\begin{lemma}\label{lem5.5}
Let $W = S_n$ and $u \in W$ with a reduced word expression $u = s_{i_1} \cdots s_{i_l}$. Choose a subset $K \subseteq [n]$. Define a subsequence $i_{j_1} , i_{j_2} , \cdots , i_{j_m}$ of $i_1 , \cdots , i_l$ by selecting $i_j$ for each $1 \leq j \leq l$ satisfying
\begin{equation}
s_{i_l} s_{i_{l-1}} \cdots s_{i_{j+1}} (i_j) \in K \quad \text{or} \quad s_{i_l} s_{i_{l-1}} \cdots s_{i_{j+1}} (1+ i_j ) \in K.
\end{equation}
Write $u' := s_{i_{j_1}} \cdots s_{i_{j_m}}$ and $u^{-1} = x_1 \cdots x_n$, $u'^{-1} = x_1' \cdots x_n'$ using one-line notation. Then
\begin{itemize}
\item $u' \leq u$.
\item For each $1 \leq t \leq n$, if $x_t' \in K$, then $x_t = x_t'$.
\item For each $1 \leq s < t \leq n$, if $x_s' , x_t' \notin K$, then $x_s' < x_t'$.
\end{itemize}
\end{lemma}
\begin{proof}
$u' \leq u$ is due to the construction of $u'$. Suppose $x_t' \in K$. Then whenever
\begin{equation}
s_{i_j} s_{i_{j+1}} \cdots s_{i_l}(x_t') \neq  s_{i_{j+1}} \cdots s_{i_l}(x_t'),
\end{equation}
$s_{i_j}$ is included in the construction of $u'$.

On the other hand, suppose $x_s ' , x_t' \notin K$ and $s<t$ but $x_s' > x_t'$. Since $s < t$, this means that there exists a $1 \leq k \leq m$ such that
\begin{equation}
s_{i_{j_k}} s_{i_{j_{k+1}}} \cdots s_{i_{j_m}}(x_s') > s_{i_{j_k}} s_{i_{j_{k+1}}} \cdots s_{i_{j_m}}(x_t'), \quad s_{i_{j_{k+1}}} \cdots s_{i_{j_m}}(x_s') < s_{i_{j_{k+1}}} \cdots s_{i_{j_m}}(x_t').
\end{equation}
This implies that $s_{i_{j_k}}$ is switching $x_s'$ and $x_t'$ which are not in $K$, leading to a contradiction.
\end{proof}

Intuitively, select $i_j$ whenever the $j$th horizontal line of the wiring diagram is included in the trajectories of some elements in $K$. For instance, if we choose $K = \{3\}$ in Figure \ref{fig0}, then the second and the fifth horizontal lines are included in the trajectory. Therefore, the corresponding subsequence is $i_2 , i_5$ and $u' = s_{i_2} s_{i_5} = s_4 s_3$.

\begin{figure}
\centering
\begin{subfigure}[t]{0.4\textwidth}
\centering
\includegraphics[width=\textwidth]{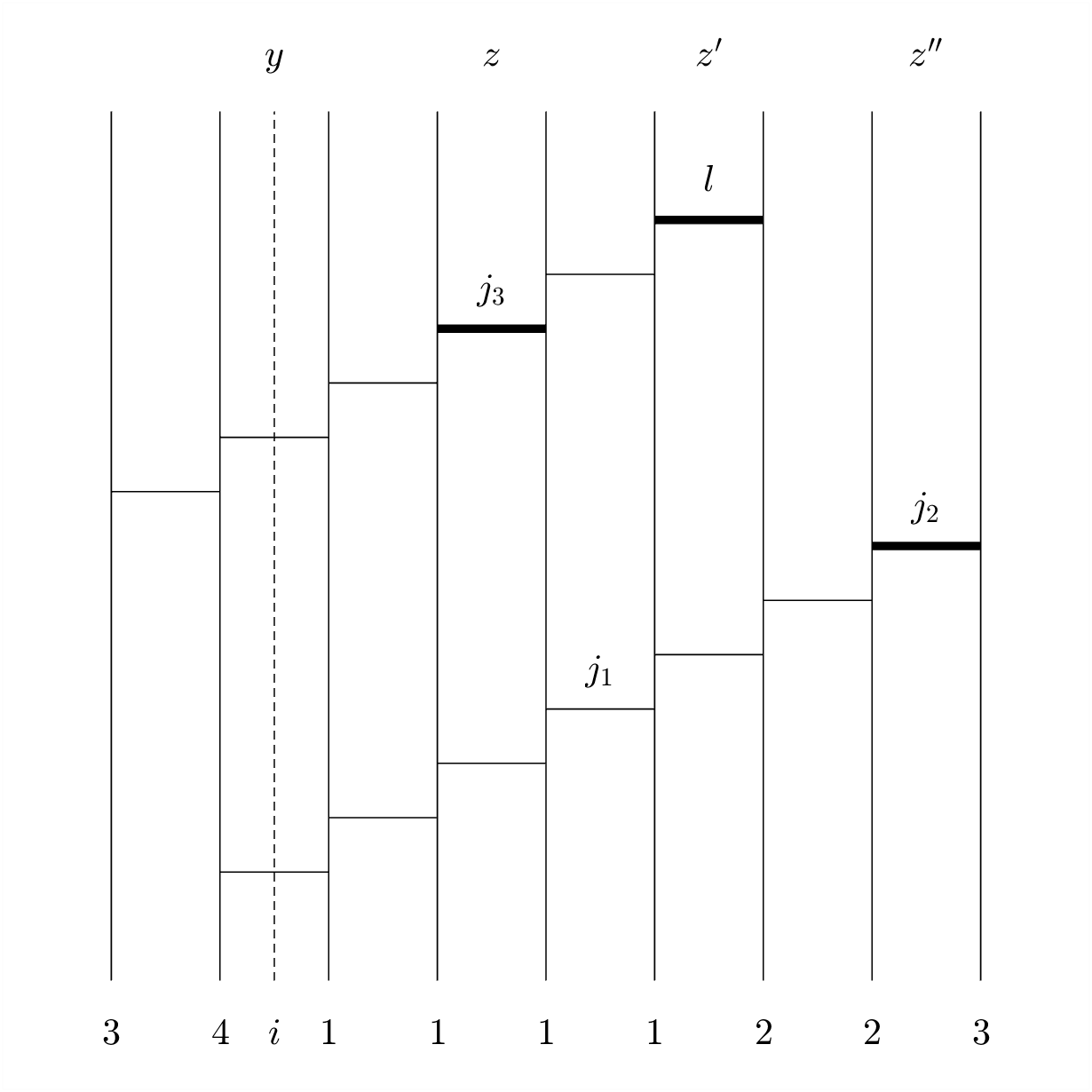}
\caption{(Case 1) $y \leq z < z' < z''$}
\label{fig1}
\end{subfigure}
\hfill
\begin{subfigure}[t]{0.4\textwidth}
\centering
\includegraphics[width=\textwidth]{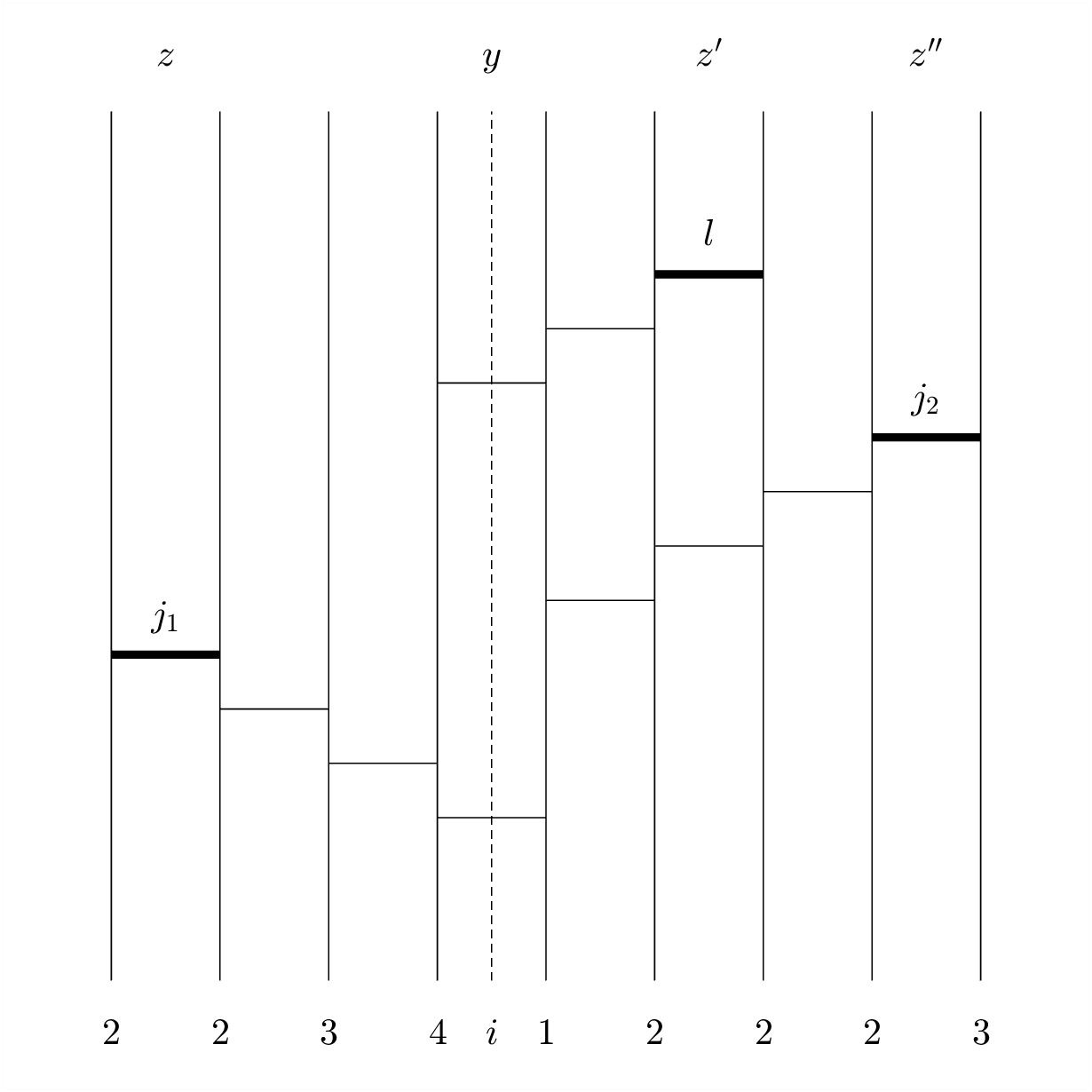}
\caption{(Case 2) $z < y \leq z' < z''$}
\label{fig2}
\end{subfigure}
\caption{$3412$ pattern in $u'$}
\label{fig012}
\end{figure}

\begin{theorem}\label{thm5.1}
Let $G= GL_n$, $W=S_n$, $I, J \subseteq [n-1]$ and $u \in {}^I W^J$. If $[u,I,J]$ has a $3412$ pattern, then $[u,I,J]$ is non-spherical.
\end{theorem}
\begin{proof}
Write $[u,J] = x_1 \cdots x_n$. Recall that $h:=n - |J|$. Since $[u,I,J]$ has a $3412$ pattern, choose $a < b< c <d$ so that $x_c < x_d <x_a < x_b$. Use $J$-expansion Lemma \ref{lem4.5} so that $1 , \cdots , x_c$ and $x_c +1 , \cdots , x_d$ and $x_d +1 , \cdots ,x_a$ and $x_a+1 , \cdots , h$ are equalized. Then $[u,I,J']$ still has a $3412$ pattern. Therefore, we may assume that $[n-1]-J = \{ z, z', z''\}$ with $z<z'<z''$ \textit{i.e.} $x_i \in [4]$.

Due to Lemma \ref{lem2.6}, there exists a $y \in [n-1]-I$ such that $b \leq y < c$. Use $I$-expansion Lemma \ref{lem4.4} on $I' := [n-1]-\{y\}$. Then $[u',I',J]$ still has a $3412$ pattern. Therefore, we may assume that $[n-1]-I = \{y\}$. In other words,
\begin{equation}\label{eqn5.1}
[u,I,J] =  \cdots 3 4 \cdots 4 i 1\cdots 12 \cdots .
\end{equation}
Choose $\lambda \in \Parj$
\begin{equation}\label{eqn5.2}
\lambda := \varpi_z + \varpi_{z'} + \varpi_{z''}.
\end{equation}
From the relation between $y$ and $z<z'<z''$, we have four cases. Due to symmetry, it is enough to cover the cases $y\leq z < z' <z''$ and $z < y \leq z' < z''$. 

$\bullet$ (Case 1) $y\leq z < z' <z''$: Due to a $3412$ pattern,
\begin{equation}\label{eqn5.57}
u(z'+1) \leq y-1, \quad u(z''+1) \leq y.
\end{equation}
As in Figure \ref{fig1}, define
\begin{equation}\label{eqn4.13}
u':=s_y s_{y+1} \cdots s_{z''-1}s_{z''} s_{y-1} s_{y} \cdots s_{z'-1}s_{z'}.
\end{equation}
Then $u' \leq u$ by applying Lemma \ref{lem5.5} with $K := \{ z'+1 , z''+1\}$ and \eqref{eqn5.57}. Also, $u'$ is the minimal representative of $W_I u'$. Denote $l = \ell (u')$ and define $(\mathbf{i}, \mathbf{m} )$ from the reduced word expression \eqref{eqn4.13} and $\lambda$ in \eqref{eqn5.2}. We can compute that
\begin{equation}\label{eqn4.14}
m_j = 
\begin{cases}
1 & j = z''-y+1 (=: j_2) \\
1 & j = z+z''-2y + 3 (=:j_3)\\
1 & j = l\\
0 & \text{otherwise}
\end{cases}
.
\end{equation}
Write $j_1:= j_3-j_2$. Define $J_{j,k}$ and $K_{j,k}$ for $1 \leq j \leq l$ and $1 \leq k \leq m_j$
\begin{equation}\label{eqn4.15}
J_{j,k}=
\begin{cases}
[l] & j = j_2 \\
[j_1]\cup [j_2 +1 , l]& j=  j_3 \\
[j_3 +2 , l] & j=l\\
\end{cases}
,
\quad
K_{j,k}=
\begin{cases}
[l] & j=j_2 \\
[j_2 + 1 , l] & j = j_3\\
[j_2 + 2,l] & j = l \\
\end{cases}
.
\end{equation}
Then the corresponding $LS$-monomials $\{\alpha_{j,k}\}$ and $\{\beta_{j,k}\}$ are
\begin{equation}\label{eqn4.17}
\alpha_{j,k} = 
\begin{cases}
\varpi_{z''} +\epsilon_{z''+1} - \epsilon_y & j = j_2 \\
\varpi_{z} + \epsilon_{z +1} + \epsilon_{z +2} - \epsilon_{y}-\epsilon_{y-1} & j = j_3 \\
\varpi_{z'} + \epsilon_{z'+1} - \epsilon_{z +2}  & j= l \\
\end{cases}
, \quad \beta_{j,k} = 
\begin{cases}
\varpi_{z''} +\epsilon_{z''+1} - \epsilon_y & j = j_2 \\
\varpi_{z} + \epsilon_{z +1}  -\epsilon_{y-1} & j = j_3 \\
\varpi_{z'} + \epsilon_{z'+1} - \epsilon_{y}  & j= l \\
\end{cases}
.
\end{equation}

$\bullet$ (Case 2) $z<y \leq z' <z''$: Due to a $3412$ pattern,
\begin{equation}\label{eqn5.58}
u(z) \geq y+1, \quad u(z' +1) \leq y-1, \quad u(z''+1) \leq y.
\end{equation}
As in Figure \ref{fig2}, define
\begin{equation}\label{eqn4.19}
u':= s_y s_{y-1} \cdots s_{z+1} s_z s_{y+1} s_{y+2} \cdots s_{z''-1} s_{z''} s_{y} s_{y+1} \cdots s_{z' -1} s_{z'}.
\end{equation}
Then $u' \leq u$ by applying Lemma \ref{lem5.5} with $K:= \{ z , z'+1 , z''+1 \}$ with \eqref{eqn5.58}. Also, $u'$ is the minimal representative of $W_I u'$. Denote $l = \ell (u')$ and define $(\mathbf{i}, \mathbf{m})$ from the reduced word expression \eqref{eqn4.19} and $\lambda$ in \eqref{eqn5.2}. We can compute that
\begin{equation}\label{eqn4.20}
m_j =
\begin{cases}
1 & j = y-z+1 (=: j_1) \\
1& j = z''-z+1 (=:j_2) \\
1 & j = l \\
0 & \text{otherwise}
\end{cases}
.
\end{equation}
Define $J_{j,k}$ and $K_{j,k}$ for $1 \leq j \leq l$ and $1 \leq k \leq m_j$
\begin{equation}\label{eqn4.21}
J_{j,k} =
\begin{cases}
[l] & j=j_1\\
[l] & j = j_2  \\
[j_2+2 , l] & j=l
\end{cases}
, \quad K_{j,k} = 
\begin{cases}
[l] & j=j_1\\
[j_1+1 , l] & j = j_2 \\
[j_2+1,l] & j = l 
\end{cases}
.
\end{equation}
Then the corresponding $LS$-monomials $\{\alpha_{j,k}\}$ and $\{\beta_{j,k}\}$ are
\begin{subequations}\label{eqn4.24}
\begin{equation}\label{eqn4.25}
\alpha_{j,k}=
\begin{cases}
\varpi_z - \epsilon_z + \epsilon_{y+1} & j=j_1 \\
\varpi_{z''} - \epsilon_{y} + \epsilon_{z''+1}& j=j_2 \\
\varpi_{z'}-\epsilon_{y+1}+\epsilon_{z'+1} & j=l
\end{cases}
, \quad \beta_{j,k} = 
\begin{cases}
\varpi_z - \epsilon_z + \epsilon_{y+1} & j=j_1 \\
\varpi_{z''} - \epsilon_{y+1} + \epsilon_{z''+1}& j=j_2 \\
\varpi_{z'}-\epsilon_{y}+\epsilon_{z'+1} & j=l
\end{cases}
.
\end{equation}
\end{subequations}

From (Case 1) to (Case 2), we constructed two distinct $I$-dominant LS-monomials $\{\alpha_{j,k}\}$ and $\{\beta_{j,k}\}$ with the same weight based on $u'$ and $\lambda$. Let $w'$ be the maximal representative of $W_I u'$. Due to Theorem \ref{thm4.1}, $H^0(X(w'W_J), L_\lambda)^*$ is not $L_I$-multiplicity-free. From $u'\leq u$, we have $w' \leq w$. Using Lemma \ref{lem2.5}, $H^0(X(wW_J), L_\lambda)^*$ is also not $L_I$-multiplicity-free. Due to Theorem \ref{lem2.4}, $X(wW_J)$ is not $L_I$-spherical.
\end{proof}

\begin{theorem}\label{thm5.2}
Let $G= GL_n$, $W=S_n$, $I, J \subseteq [n-1]$ and $u \in {}^I W^J$. If $[u,I,J]$ has a $321$ pattern, then $[u,I,J]$ is non-spherical.
\end{theorem}
\begin{proof}
Write $[u,J] = x_1 \cdots x_n$. Recall that $h:= n - |J|$. Since $[u,I,J]$ has a $321$ pattern, choose $a<b<c$ so that $x_c < x_b < x_a$. Use $J$-expansion Lemma \ref{lem4.5} so that $1 , \cdots , x_c$ and $x_c + 1 , \cdots , x_b$ and $x_b +1 , \cdots h$ are equalized. Then $[u, I , J']$ still has a $321$ pattern. Therefore, we may assume $[n-1] - J = \{z, z'\}$ with $z<z'$ \textit{i.e.} $x_j \in [3]$.

Due to Lemma \ref{lem2.6}, there exist $y , y' \in [n-1]-I$ such that $a \leq y < b \leq y' <c$. Use $I$-expansion Lemma \ref{lem4.4} on $I':= [n-1]-\{y,y'\}$. Then $[u',I',J]$ still has a $321$ pattern. Therefore, we may assume that $[n-1] - I = \{y , y'\}$. In other words,
\begin{equation}\label{eqn5.3}
[u,I,J] = \cdots 3i\cdots 2\cdots i1\cdots  .
\end{equation}

Choose $\lambda \in \Parj$
\begin{equation}\label{eqn5.4}
\lambda:= \varpi_z + \varpi_{z'}.
\end{equation}
From the relation between $y<y'$ and $z< z'$, we have six cases to consider. Due to symmetry, it is enough to cover the cases $z<z' \leq y<y'$, $y \leq z < z' \leq y'$, $z \leq y < z' \leq y'$ and $z\leq y<y' \leq z'$.

$\bullet$ (Case 1) $z<z' \leq y<y'$: Due to a $321$ pattern,
\begin{equation}\label{eqn5.59}
u(z) \geq y'+1, \quad u(z') \geq y+1.
\end{equation} 
As in Figure \ref{fig3}, define
\begin{equation}\label{eqn5.5}
u':=s_y s_{y-1} \cdots s_{z'}s_{z'-1} s_{y'} s_{y'-1} \cdots s_{z+1}s_z.
\end{equation}
Then $u' \leq u$ by applying Lemma \ref{lem5.5} with $K:=\{z,z'\}$ and \eqref{eqn5.59}. Also, $u'$ is the minimal representative of $W_Iu'$. Denote $l:= \ell (u')$ and define $(\mathbf{i}, \mathbf{m})$ from the reduced word expression in \eqref{eqn5.5} and $\lambda$ in \eqref{eqn5.4}. We can compute
\begin{equation}\label{eqn5.6}
m_j = 
\begin{cases}
1 & j = y'+y-2z'+3 (=:j_3) \\
1 & j = l \\
0 & \text{otherwise}
\end{cases}
\end{equation}

Write $j_1:=y-z'+3$ and $j_2 = y'-z'+3$. Define $J_{j,k}$ and $K_{j,k}$ for $1\leq j \leq l$ and $1 \leq k \leq m_j$
\begin{equation}\label{eqn5.7}
J_{j,k} = 
\begin{cases}
[l] & j=j_3 \\
[j_3+2,l] & j = l
\end{cases}
, \quad K_{j,k}=
\begin{cases}
[j_1,l] & j =j_3 \\
[j_2,l] & j =l
\end{cases}
.
\end{equation}
Then the corresponding LS-monomials $\{\alpha_{j,k}\}$ and $\{\beta_{j,k}\}$ are
\begin{equation}\label{eqn5.8}
\alpha_{j,k} = 
\begin{cases}
\varpi_{z'} - \epsilon_{z'-1}-\epsilon_{z'} + \epsilon_{y+1} + \epsilon_{y'+1} & j = j_3 \\
\varpi_z-\epsilon_{z} + \epsilon_{z'-1} & j=l 
\end{cases}
, \quad \beta_{j,k}=
\begin{cases}
\varpi_{z'}-\epsilon_{z'}+\epsilon_{y'+1} & j=j_3 \\
\varpi_{z}-\epsilon_z + \epsilon_{y+1} & j=l
\end{cases}
.
\end{equation}

$\bullet$ (Case 2) $y \leq z <z' \leq y'$: Due to a $321$ pattern,
\begin{equation}\label{eqn5.60}
u(z) \geq y'+1, \quad u(z' +1 ) \leq y.
\end{equation}
As in Figure \ref{fig4}, define
\begin{equation}\label{eqn5.9}
u':=s_y s_{y+1} \cdots s_{z'-2}s_{z'-1} s_{y'}s_{y'-1} \cdots s_{z+1}s_z.
\end{equation}
Then $u' \leq u$ by applying Lemma \ref{lem5.5} with $K:=\{z, z' +1 \}$ and \eqref{eqn5.60}. Also, $u'$ is the minimal representative of $W_Iu'$. Denote $l:= \ell (u')$ and define $(\mathbf{i}, \mathbf{m})$ from the reduced word expression in \eqref{eqn5.9} and $\lambda$ in \eqref{eqn5.4}. We can compute
\begin{equation}\label{eqn5.10}
m_j = 
\begin{cases}
1 & j = -y+y'+1 (=:j_3) \\
1 & j = l \\
0 & \text{otherwise}
\end{cases}
\end{equation}

Write $j_1:= z-y+1$ and $j_2 = z'-y+1$. Define $J_{j,k}$ and $K_{j,k}$ for $1\leq j \leq l$ and $1 \leq k \leq m_j$
\begin{equation}\label{eqn5.11}
J_{j,k} = 
\begin{cases}
[l] & j=j_3 \\
\emptyset& j = l
\end{cases}
, \quad K_{j,k}
\begin{cases}
[j_1-1] \cup [j_2, l] & j =j_3 \\
[j_1-1] \cup [j_3+1, l] & j =l
\end{cases}
.
\end{equation}
Then the corresponding LS-monomials $\{\alpha_{j,k}\}$ and $\{\beta_{j,k}\}$ are
\begin{equation}\label{eqn5.12}
\alpha_{j,k} = 
\begin{cases}
\varpi_{z'} - \epsilon_{y}+\epsilon_{y'+1} & j = j_3 \\
\varpi_z & j=l 
\end{cases}
, \quad \beta_{j,k}=
\begin{cases}
\varpi_{z'}-\epsilon_{z'}+\epsilon_{y'+1} & j=j_3 \\
\varpi_{z}-\epsilon_y + \epsilon_{z'} & j=l
\end{cases}
.
\end{equation}

\begin{figure}
\centering
\begin{subfigure}[t]{0.45\textwidth}
\centering
\includegraphics[width=\textwidth]{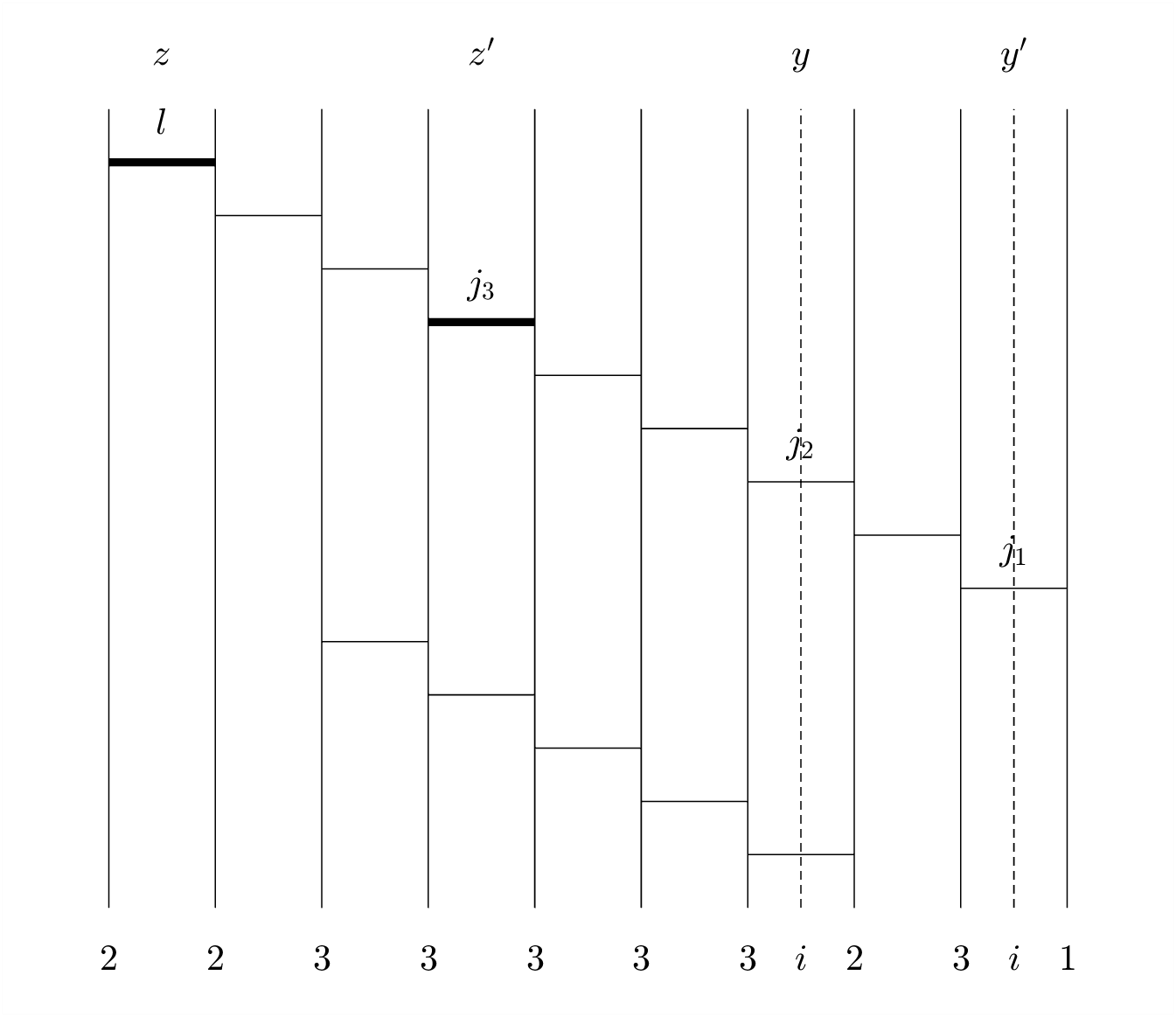}
\caption{(Case 1) $z<z' \leq y < y'$}
\label{fig3}
\end{subfigure}
\hfill
\begin{subfigure}[t]{0.45\textwidth}
\centering
\includegraphics[width=\textwidth]{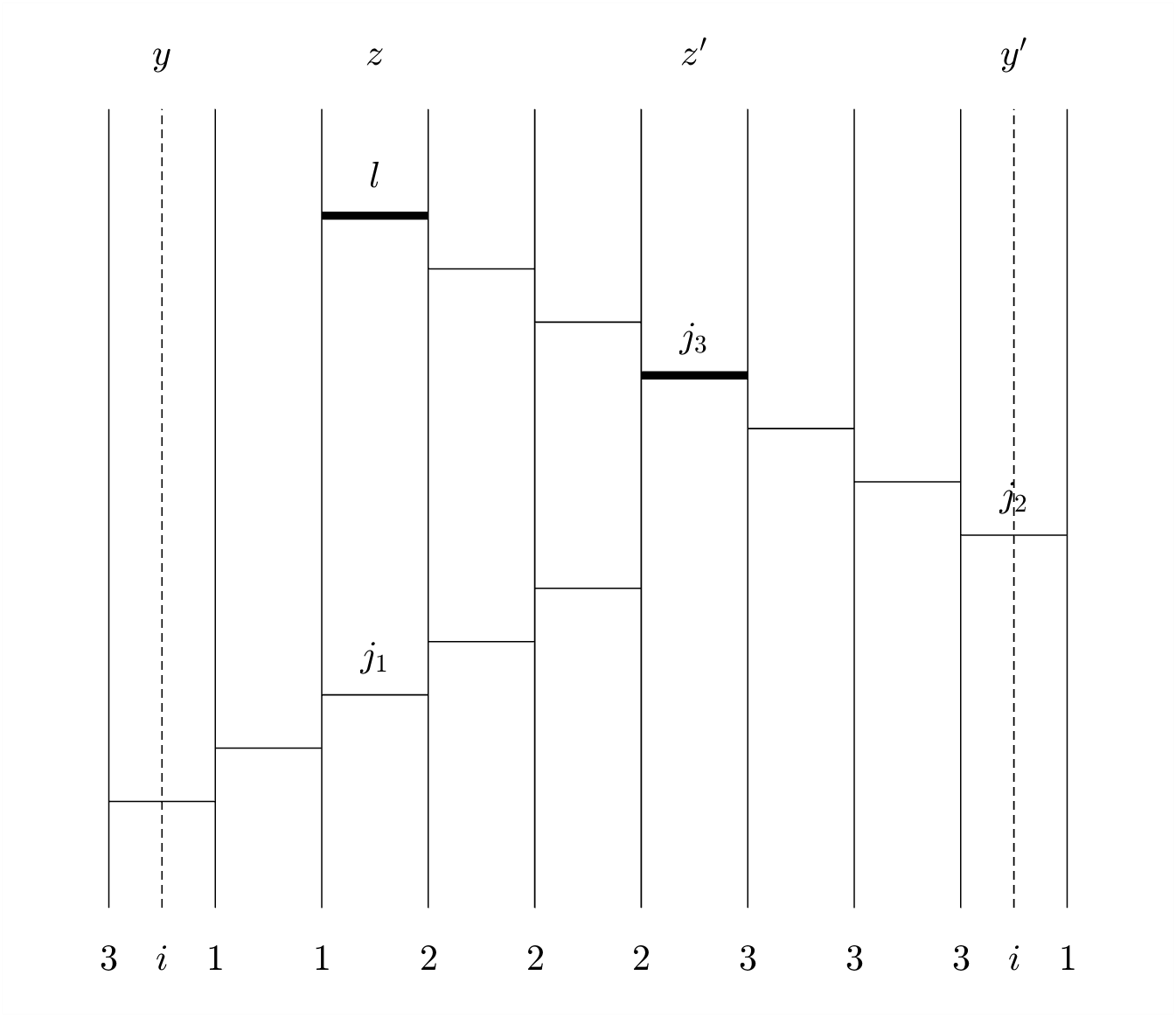}
\caption{(Case 2) $y \leq z < z' \leq y'$}
\label{fig4}
\end{subfigure}

\begin{subfigure}[t]{0.45\textwidth}
\centering
\includegraphics[width=\textwidth]{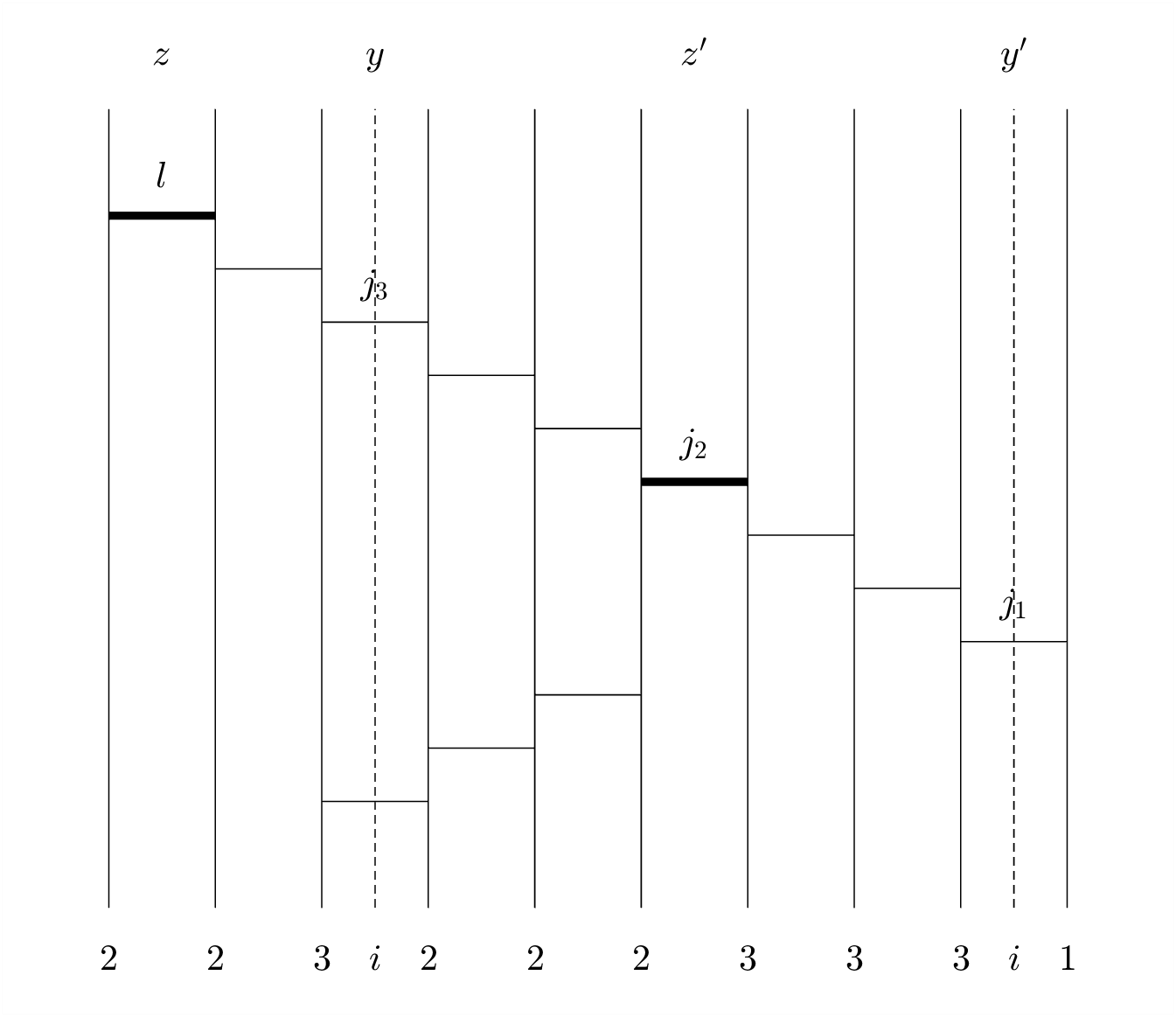}
\caption{(Case 3) $z \leq y < z' \leq y'$}
\label{fig5}
\end{subfigure}
\hfill
\begin{subfigure}[t]{0.45\textwidth}
\centering
\includegraphics[width=\textwidth]{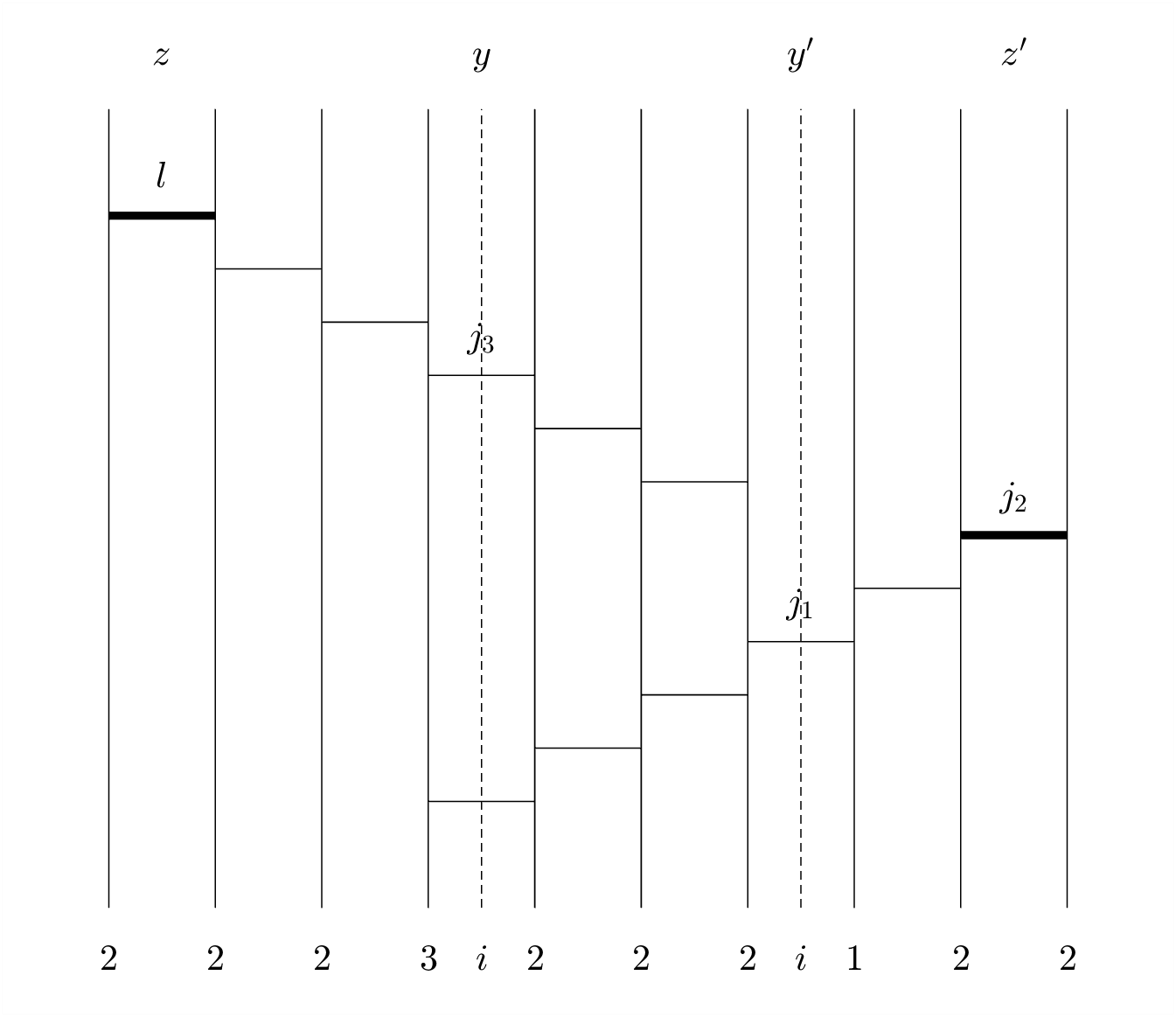}
\caption{(Case 4) $z \leq y < y' \leq z'$}
\label{fig6}
\end{subfigure}
\caption{$321$ patterns in $u'$}
\label{fig03456}
\end{figure}

$\bullet$ (Case 3) $z \leq y  < z' \leq y'$: Due to a $321$ pattern,
\begin{equation}\label{eqn5.61}
u(z) \geq y'+1, \quad u(z'+1)\leq y.
\end{equation}
As in Figure \ref{fig5}, define
\begin{equation}\label{eqn5.13}
u':=s_y s_{y+1} \cdots s_{z'-2} s_{z'-1} s_{y'}s_{y'-1} \cdots s_{z+1}s_z.
\end{equation}
Then $u' \leq u$ by applying Lemma \ref{lem5.5} with $K:=\{z,z'+1 \}$ and \eqref{eqn5.61}. Also, $u'$ is the minimal representative of $W_Iu'$. Denote $l:= \ell (u')$ and define $(\mathbf{i}, \mathbf{m})$ from the reduced word expression in \eqref{eqn5.13} and $\lambda$ in \eqref{eqn5.4}. We can compute
\begin{equation}\label{eqn5.14}
m_j = 
\begin{cases}
1 & j = -y+y'+1 (=:j_2) \\
1 & j = l \\
0 & \text{otherwise}
\end{cases}
\end{equation}

Write $j_1:= z'-y+1$ and $j_3:=z'+y'-2y+1$. Define $J_{j,k}$ and $K_{j,k}$ for $1\leq j \leq l$ and $1 \leq k \leq m_j$
\begin{equation}\label{eqn5.15}
J_{j,k} = 
\begin{cases}
[l] & j=j_2 \\
[j_3+1,l] & j = l
\end{cases}
, \quad K_{j,k}
\begin{cases}
[j_1,l] & j =j_2 \\
[j_2+1,l] & j =l
\end{cases}
.
\end{equation}
Then the corresponding LS-monomials $\{\alpha_{j,k}\}$ and $\{\beta_{j,k}\}$ are
\begin{equation}\label{eqn5.16}
\alpha_{j,k} = 
\begin{cases}
\varpi_{z'} - \epsilon_{y}+\epsilon_{y'+1} & j = j_2 \\
\varpi_z-\epsilon_z + \epsilon_y & j=l 
\end{cases}
, \quad \beta_{j,k}=
\begin{cases}
\varpi_{z'}-\epsilon_{z'}+\epsilon_{y'+1} & j=j_2 \\
\varpi_{z}-\epsilon_z + \epsilon_{z'} & j=l
\end{cases}
.
\end{equation}

$\bullet$ (Case 4) $z \leq y  < y' \leq  z'$: Due to a $321$ pattern,
\begin{equation}\label{eqn5.62}
u(z) \geq y'+1, \quad u(z'+1) \leq y.
\end{equation} 
As in Figure \ref{fig6}, define
\begin{equation}\label{eqn5.17}
u':=s_y s_{y+1} \cdots s_{z'-1} s_{z'} s_{y'-1}s_{y'-2} \cdots s_{z+1}s_z.
\end{equation}
Then $u' \leq u$ by applying Lemma \ref{lem5.5} with $K:=\{z, z'+1 \}$ and \eqref{eqn5.62}. Also, $u'$ is the minimal representative of $W_Iu'$. Denote $l:= \ell (u')$ and define $(\mathbf{i}, \mathbf{m})$ from the reduced word expression in \eqref{eqn5.17} and $\lambda$ in \eqref{eqn5.4}. We can compute
\begin{equation}\label{eqn5.18}
m_j = 
\begin{cases}
1 & j = z'-y+1 (=:j_2) \\
1 & j = l \\
0 & \text{otherwise}
\end{cases}
\end{equation}

Write $j_1:= y'-y+1$ and $j_3:= z'+y'-2y+1$. Define $J_{j,k}$ and $K_{j,k}$ for $1\leq j \leq l$ and $1 \leq k \leq m_j$
\begin{equation}\label{eqn5.19}
J_{j,k} = 
\begin{cases}
[l] & j=j_2 \\
[j_3+1,l] & j = l
\end{cases}
, \quad K_{j,k}
\begin{cases}
[j_1,l] & j =j_2 \\
[j_2+1,l] & j =l
\end{cases}
.
\end{equation}
Then the corresponding LS-monomials $\{\alpha_{j,k}\}$ and $\{\beta_{j,k}\}$ are
\begin{equation}\label{eqn5.20}
\alpha_{j,k} = 
\begin{cases}
\varpi_{z'} - \epsilon_{y}+\epsilon_{z'+1} & j = j_2 \\
\varpi_z-\epsilon_z + \epsilon_y & j=l 
\end{cases}
, \quad \beta_{j,k}=
\begin{cases}
\varpi_{z'}-\epsilon_{y'}+\epsilon_{z'+1} & j=j_2 \\
\varpi_{z}-\epsilon_z + \epsilon_{y'} & j=l
\end{cases}
.
\end{equation}

From (Case 1) to (Case 4), we constructed two distinct $I$-dominant LS-monomials $\{\alpha_{j,k}\}$ and $\{\beta_{j,k}\}$ with the same weight based on $u'$ and $\lambda$. Let $w'$ be the maximal representative of $W_I u'$. Due to Theorem \ref{thm4.1}, $H^0(X(w'W_J), L_\lambda)^*$ is not $L_I$-multiplicity-free. From $u'\leq u$, we have $w' \leq w$. Using Lemma \ref{lem2.5}, $H^0(X(wW_J), L_\lambda)^*$ also is not $L_I$-multiplicity-free. Due to Theorem \ref{lem2.4}, $X(wW_J)$ is not $L_I$-spherical.
\end{proof}

\subsection{Toricity from partial flags}

In this subsection, we discuss $3i312$, $231i1$ and $2i21i1$ patterns, observed in Theorem \ref{thm2.1}. In fact, according to \cite[Theorem 1.1]{MR4599567}, the key polynomial indexed by $(0,0,1,1)$ is multiplicity-free, whereas the one indexed by $(0,0,2,2)$ is not. This is reflected in \eqref{eqn5.48} where we choose $2 \varpi_z$ instead of $\varpi_z$.

\begin{theorem}\label{thm5.3}
Let $G= GL_n$, $W=S_n$, $I, J \subseteq [n-1]$ and $u \in {}^I W^J$. If $[u,I,J]$ has $3i312$ or $231i1$ patterns, then $[u,I,J]$ is non-spherical.
\end{theorem}
\begin{proof}
Without loss of generality, assume that $[u,I,J]$ has a $3i312$ pattern. Write $[u,J] = x_1 \cdots x_n$. Recall that $h:= n - |J|$. Choose $a\leq y < b < c < d$ so that $x_c < x_d < x_a = x_b$ and $y \in [n-1]-I$. Use $J$-expansion Lemma \ref{lem4.5} so that $1 , \cdots , x_c$ and $x_c + 1 , \cdots , x_d$ and $x_d +1 , \cdots , h$ are equalized. Then $[u, I , J']$ still has a $3i312$ pattern. Therefore, we may assume $[n-1] - J = \{z, z'\}$ with $z<z'$ \textit{i.e.} $x_j \in [3]$.

Due to Lemma \ref{lem2.6}, there exists a $y' \in [n-1]-I$ such that $b \leq y' < c$. Use $I$-expansion Lemma \ref{lem4.4} on $I':= [n-1]-\{y,y'\}$. Then $[u',I',J]$ still has a $3i312$ pattern. Therefore, we may assume that $[n-1] - I = \{y , y'\}$. In other words,
\begin{equation}\label{eqn5.21}
[u,I,J] = \cdots 3 i \cdots 3i 1 \cdots 1 2 \cdots .
\end{equation}

Choose $\lambda \in \Parj$
\begin{equation}\label{eqn5.22}
\lambda:= \varpi_z + \varpi_{z'}.
\end{equation}
From the relation between $y<y'$ and $z< z'$, we have six cases to consider: $z<z' \leq y<y'$, $z\leq y<y' \leq z'$, $y \leq z < z' \leq y'$, $z \leq y < z' \leq y'$, $y<y' \leq z < z'$ and $y \leq z < y' \leq z'$.

\begin{figure}
\centering
\begin{subfigure}[t]{0.45\textwidth}
\centering
\includegraphics[width=\textwidth]{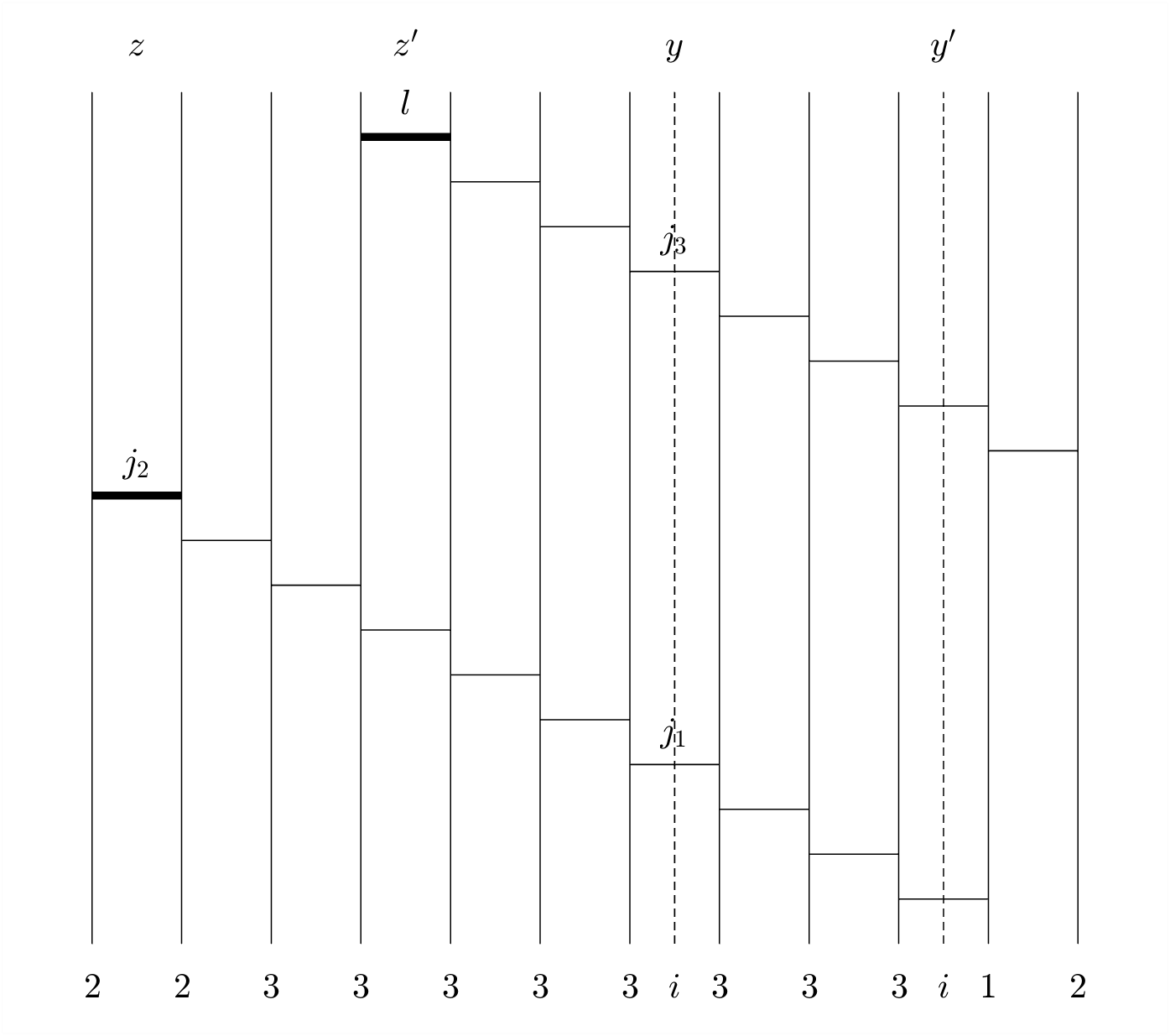}
\caption{(Case 1) $z<z' \leq y < y'$}
\label{fig7}
\end{subfigure}
\hfill
\begin{subfigure}[t]{0.45\textwidth}
\centering
\includegraphics[width=\textwidth]{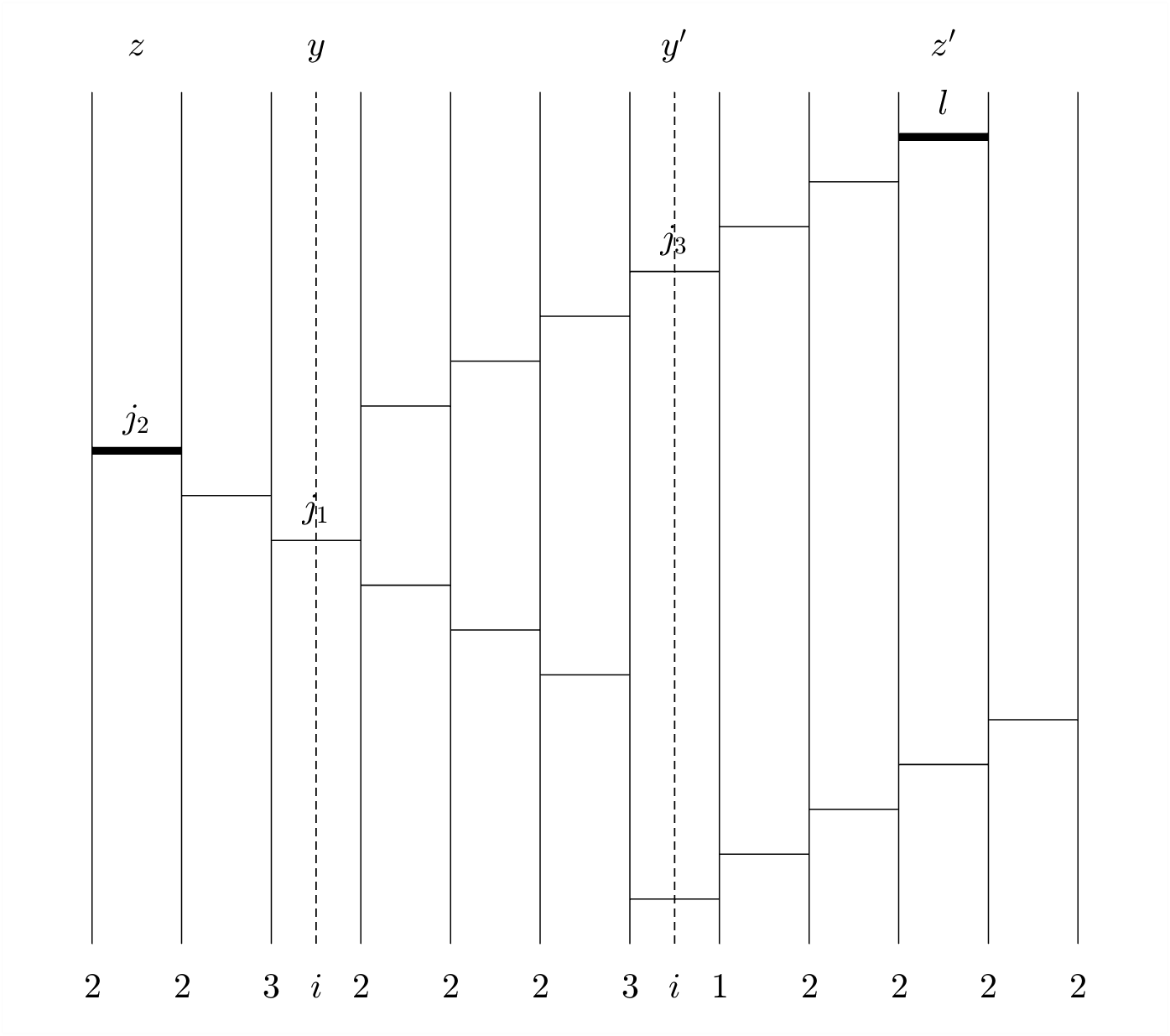}
\caption{(Case 2) $z \leq y < y' \leq z'$}
\label{fig8}
\end{subfigure}

\begin{subfigure}[t]{0.45\textwidth}
\centering
\includegraphics[width=\textwidth]{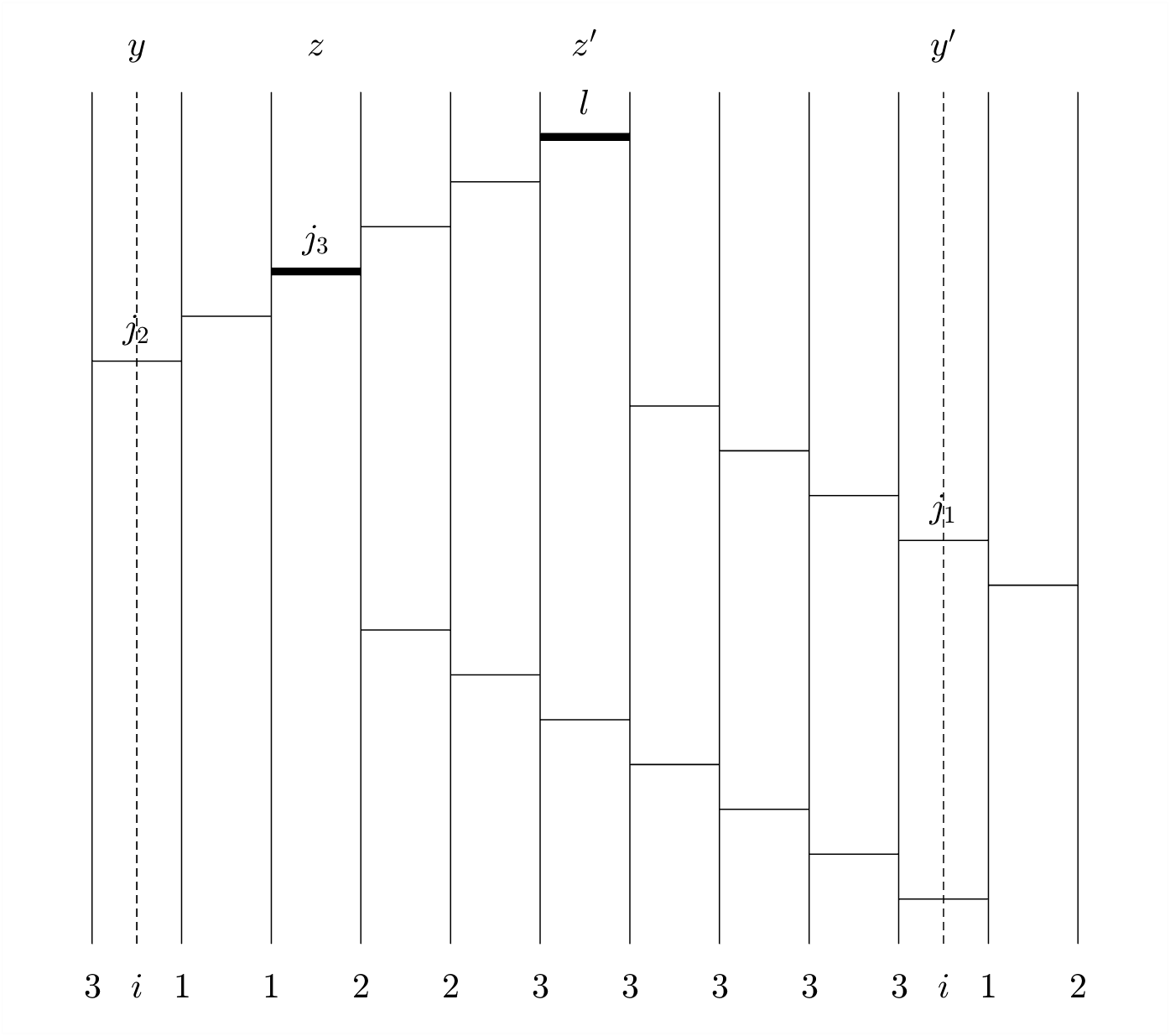}
\caption{(Case 3) $y \leq z < z' \leq y'$}
\label{fig9}
\end{subfigure}
\hfill
\begin{subfigure}[t]{0.45\textwidth}
\centering
\includegraphics[width=\textwidth]{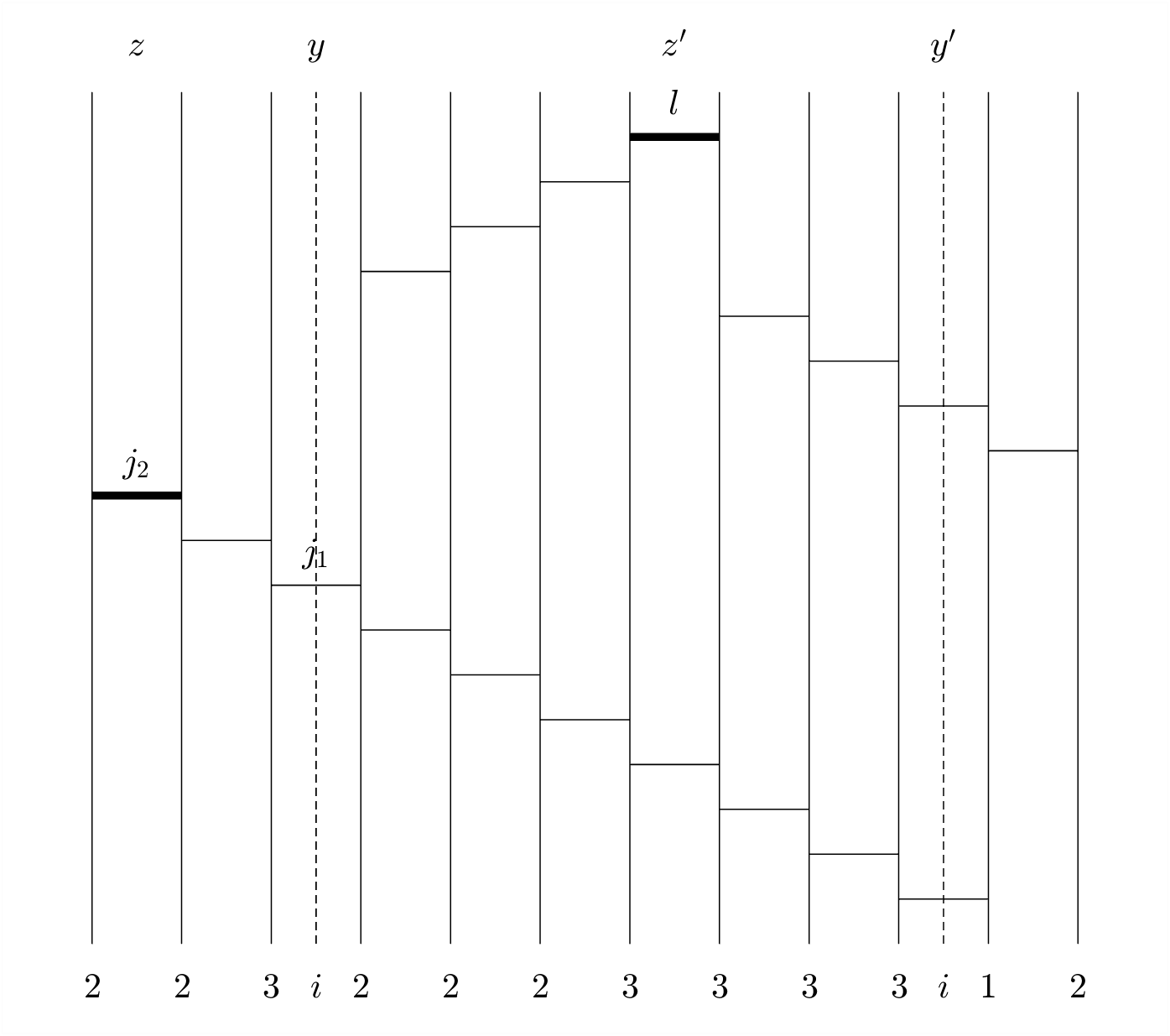}
\caption{(Case 4) $z \leq y < z' \leq y'$}
\label{fig10}
\end{subfigure}

\begin{subfigure}[t]{0.45\textwidth}
\centering
\includegraphics[width=\textwidth]{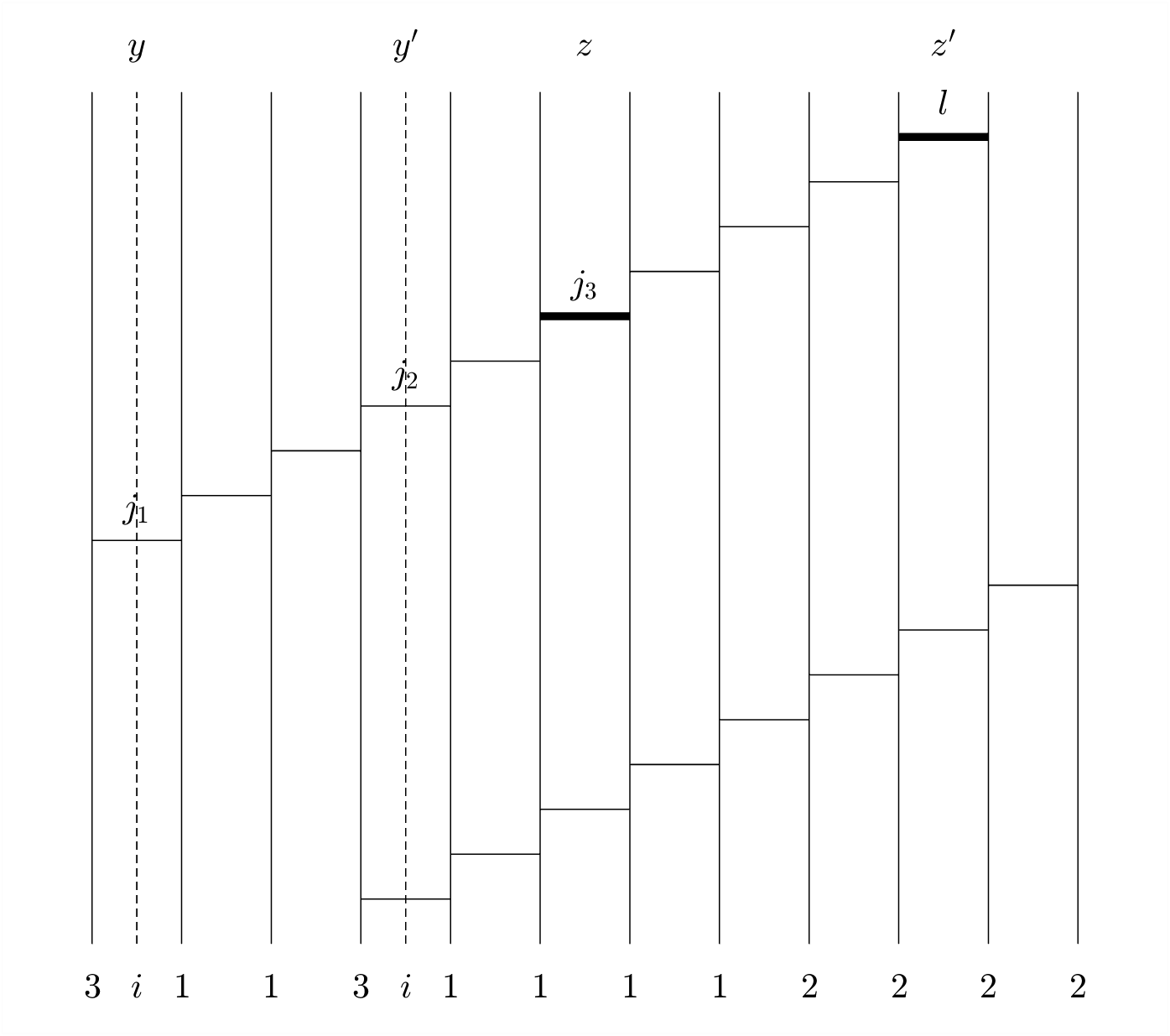}
\caption{(Case 5) $y < y' \leq z < z'$}
\label{fig11}
\end{subfigure}
\hfill
\begin{subfigure}[t]{0.45\textwidth}
\centering
\includegraphics[width=\textwidth]{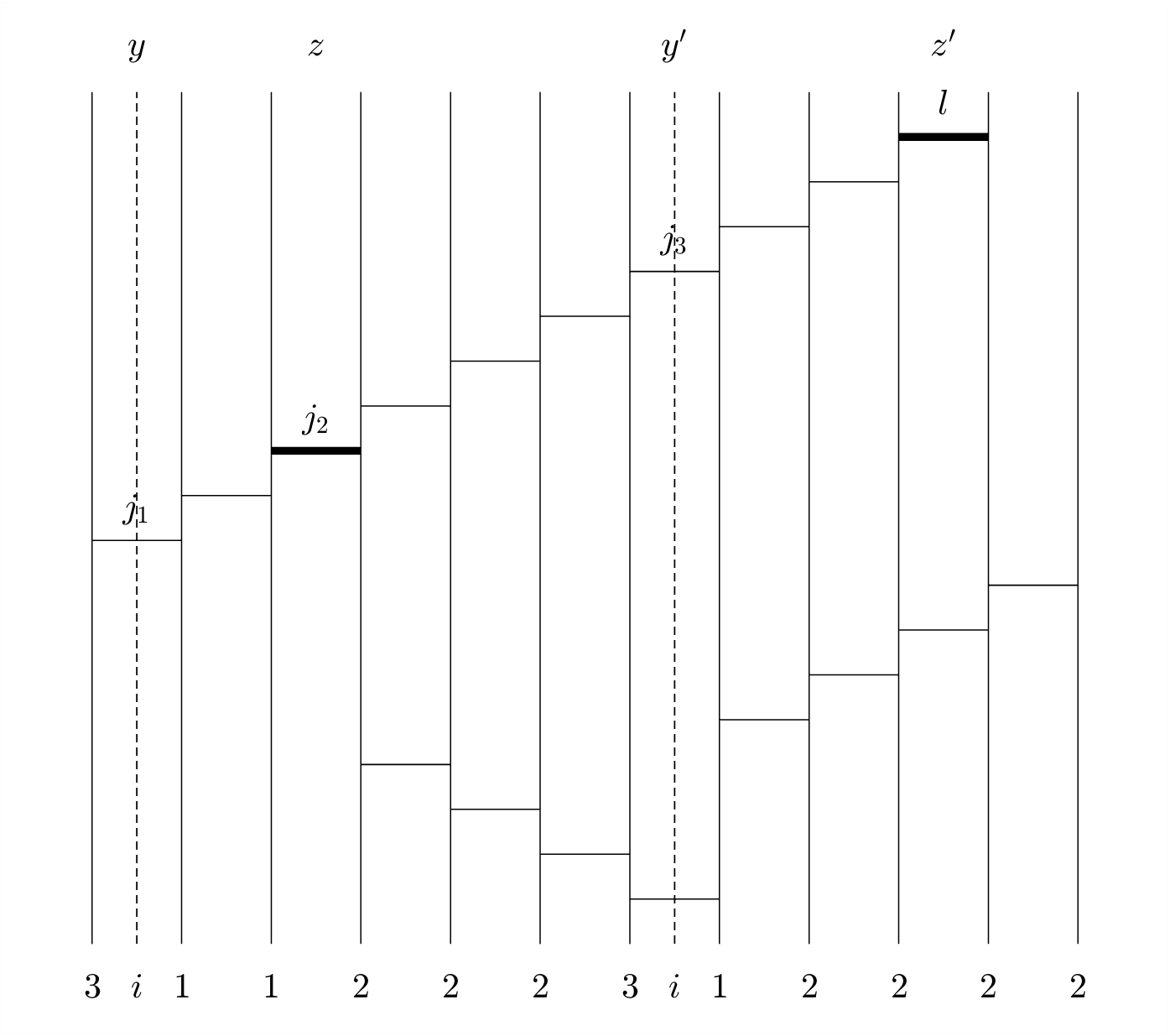}
\caption{(Case 6) $y \leq z < y' \leq z'$}
\label{fig12}
\end{subfigure}
\caption{$3i312$ patterns in $u'$}
\label{fig0789101112}
\end{figure}

$\bullet$ (Case 1) $z<z' \leq y<y'$: Due to a $3i312$ pattern,
\begin{equation}\label{eqn5.63}
u(z) \geq y'+1, \quad u(z') \geq y'+2.
\end{equation}
As in Figure \ref{fig7}, define
\begin{equation}\label{eqn5.23}
u':=s_{y'} s_{y'-1} \cdots s_{z+1}s_z s_{y'+1}s_{y'} \cdots s_{z'+1}s_{z'}.
\end{equation}
Then $u' \leq u$ by applying Lemma \ref{lem5.5} with $K:=\{z,z'\}$ and \eqref{eqn5.63}. Also, $u'$ is the minimal representative of $W_Iu'$. Denote $l:= \ell (u')$ and define $(\mathbf{i}, \mathbf{m})$ from the reduced word expression in \eqref{eqn5.23} and $\lambda$ in \eqref{eqn5.22}. We can compute
\begin{equation}\label{eqn5.24}
m_j = 
\begin{cases}
1 & j = y'-z+1 (=:j_2) \\
1 & j = l \\
0 & \text{otherwise}
\end{cases}
.
\end{equation}

Write $j_1:= y'-y+1$ and $j_3:= 2y'-y-z+3$. Define $J_{j,k}$ and $K_{j,k}$ for $1\leq j \leq l$ and $1 \leq k \leq m_j$
\begin{equation}\label{eqn5.25}
J_{j,k} = 
\begin{cases}
[l] & j=j_2 \\
[j_3,l] & j = l
\end{cases}
, \quad K_{j,k}
\begin{cases}
[j_1,l] & j =j_2 \\
[j_2+2,l] & j =l
\end{cases}
.
\end{equation}
Then the corresponding LS-monomials $\{\alpha_{j,k}\}$ and $\{\beta_{j,k}\}$ are
\begin{equation}\label{eqn5.26}
\alpha_{j,k} = 
\begin{cases}
\varpi_z-\epsilon_z + \epsilon_{y'+1} & j = j_2 \\
\varpi_{z'}-\epsilon_{z'} + \epsilon_{y+1} & j=l 
\end{cases}
, \quad \beta_{j,k}=
\begin{cases}
\varpi_z-\epsilon_z + \epsilon_{y+1} & j = j_2 \\
\varpi_{z'}-\epsilon_{z'} + \epsilon_{y'+1} & j=l 
\end{cases}
.
\end{equation}

$\bullet$ (Case 2) $z \leq y <y' \leq z'$: Due to a $3i312$ pattern,
\begin{equation}\label{eqn5.64}
u(z)\geq y'+1, \quad u(z' +1) \leq y, \quad u(z'+2) \leq y'.
\end{equation}
As in Figure \ref{fig8}, define
\begin{equation}\label{eqn5.27}
u':=s_{y'} s_{y'+1} \cdots s_{z'}s_{z'+1} s_{y'-1} s_{y'-2} \cdots s_{z+1} s_z s_{y+1} s_{y+2} \cdots s_{z'-1} s_{z'}. 
\end{equation}
Then $u' \leq u$ by applying Lemma \ref{lem5.5} with $K:=\{z , z'+1 , z'+2\}$ and \eqref{eqn5.64}. Also, $u'$ is the minimal representative of $W_Iu'$. Denote $l:= \ell (u')$ and define $(\mathbf{i}, \mathbf{m})$ from the reduced word expression in \eqref{eqn5.27} and $\lambda$ in \eqref{eqn5.22}. We can compute
\begin{equation}\label{eqn5.28}
m_j = 
\begin{cases}
1 & j = z'-z+2 (=:j_2) \\
1 & j = l \\
0 & \text{otherwise}
\end{cases}
.
\end{equation}

Write $j_1:= z'-y+2$, $j_3:= z'-z+y'-y+2$. Define $J_{j,k}$ and $K_{j,k}$ for $1\leq j \leq l$ and $1 \leq k \leq m_j$
\begin{equation}\label{eqn5.29}
J_{j,k} = 
\begin{cases}
[l] & j=j_2 \\
[j_3+1,l] & j = l
\end{cases}
, \quad K_{j,k}
\begin{cases}
[j_1,l] & j =j_2 \\
[j_2+1,l] & j =l
\end{cases}
.
\end{equation}
Then the corresponding LS-monomials $\{\alpha_{j,k}\}$ and $\{\beta_{j,k}\}$ are
\begin{equation}\label{eqn5.30}
\alpha_{j,k} = 
\begin{cases}
\varpi_{z}-\epsilon_z+\epsilon_{y'+1} & j = j_2 \\
\varpi_{z'}-\epsilon_{y'+1} + \epsilon_{z'+1} & j=l 
\end{cases}
, \quad \beta_{j,k}=
\begin{cases}
\varpi_{z}-\epsilon_{z} + \epsilon_{y+1} & j = j_2 \\
\varpi_{z'}-\epsilon_{y+1} + \epsilon_{z'+1} & j=l 
\end{cases}
.
\end{equation}

$\bullet$ (Case 3) $y \leq z<z' \leq y'$: Due to a $3i312$ pattern,
\begin{equation}\label{eqn5.65}
u(z) \geq y'+1, \quad u(z') \geq y'+2, \quad u(z'+1) \leq y.
\end{equation}
As in Figure \ref{fig9}, define
\begin{equation}\label{eqn5.31}
u':= s_{y'} s_{y'-1} \cdots s_{z+2}s_{z+1} s_{y'+1} s_{y'} \cdots s_{z'+2}s_{z'+1} s_y s_{y+1} \cdots s_{z'-1}s_{z'}.
\end{equation}
Then $u' \leq u$ by applying Lemma \ref{lem5.5} with $K:=\{ z , z' , z'+1 \}$ and \eqref{eqn5.65}. Also, $u'$ is the minimal representative of $W_Iu'$. Denote $l:= \ell (u')$ and define $(\mathbf{i}, \mathbf{m})$ from the reduced word expression in \eqref{eqn5.31} and $\lambda$ in \eqref{eqn5.22}. We can compute
\begin{equation}\label{eqn5.32}
m_j = 
\begin{cases}
1 & j = 2y'-y-z'+2(=:j_3) \\
1 & j = l \\
0 & \text{otherwise}
\end{cases}
.
\end{equation}

Write $j_1:= y'-z+2$ and $j_2:= 2y' - z - z' +2$. Define $J_{j,k}$ and $K_{j,k}$ for $1\leq j \leq l$ and $1 \leq k \leq m_j$
\begin{equation}\label{eqn5.33}
J_{j,k} = 
\begin{cases}
[l] & j=j_3 \\
\emptyset & j = l
\end{cases}
, \quad K_{j,k}
\begin{cases}
[j_1,l] & j =j_3 \\
[j_1, j_2-1] \cup [j_3+1,l]& j =l
\end{cases}
.
\end{equation}
Then the corresponding LS-monomials $\{\alpha_{j,k}\}$ and $\{\beta_{j,k}\}$ are
\begin{equation}\label{eqn5.34}
\alpha_{j,k} = 
\begin{cases}
\varpi_{z}-\epsilon_y+\epsilon_{y'+1} & j = j_3 \\
\varpi_{z'} & j=l 
\end{cases}
, \quad \beta_{j,k}=
\begin{cases}
\varpi_{z}-\epsilon_{y} + \epsilon_{z+1} & j = j_3 \\
\varpi_{z'}-\epsilon_{z+1} + \epsilon_{y'+1} & j=l 
\end{cases}
.
\end{equation}

$\bullet$ (Case 4) $z \leq y < z' \leq y'$: Due to a $3i312$ pattern,
\begin{equation}\label{eqn5.66}
u(z) \geq y'+1, \quad u(z') \geq y'+2, \quad u(z'+1) \leq y.
\end{equation}
As in Figure \ref{fig10}, define
\begin{equation}\label{eqn5.35}
u':=s_{y'}s_{y'-1} \cdots s_{z+1} s_z s_{y'+1} s_{y'} \cdots s_{z'+2} s_{z'+1} s_{y+1} s_{y+2} \cdots s_{z'-1} s_{z'}.
\end{equation}
Then $u' \leq u$ by applying Lemma \ref{lem5.5} with $K:=\{z , z' , z'+1 \}$ and \eqref{eqn5.66}. Also, $u'$ is the minimal representative of $W_Iu'$. Denote $l:= \ell (u')$ and define $(\mathbf{i}, \mathbf{m})$ from the reduced word expression in \eqref{eqn5.35} and $\lambda$ in \eqref{eqn5.22}. We can compute
\begin{equation}\label{eqn5.36}
m_j = 
\begin{cases}
1 & j = y'-z+1 (=:j_2) \\
1 & j = l \\
0 & \text{otherwise}
\end{cases}
.
\end{equation}

Write $j_1:= y'-y+1$. Define $J_{j,k}$ and $K_{j,k}$ for $1\leq j \leq l$ and $1 \leq k \leq m_j$
\begin{equation}\label{eqn5.37}
J_{j,k} = 
\begin{cases}
[l] & j=j_2 \\
\emptyset & j = l
\end{cases}
, \quad K_{j,k}
\begin{cases}
[j_1,l] & j =j_2 \\
[j_2+2,l] & j =l
\end{cases}
.
\end{equation}
Then the corresponding LS-monomials $\{\alpha_{j,k}\}$ and $\{\beta_{j,k}\}$ are
\begin{equation}\label{eqn5.38}
\alpha_{j,k} = 
\begin{cases}
\varpi_{z}-\epsilon_z+\epsilon_{y'+1} & j = j_2 \\
\varpi_{z'} & j=l 
\end{cases}
, \quad \beta_{j,k}=
\begin{cases}
\varpi_{z}-\epsilon_{z} + \epsilon_{y+1} & j = j_2 \\
\varpi_{z'}-\epsilon_{y+1} + \epsilon_{y'+1} & j=l 
\end{cases}
.
\end{equation}

$\bullet$ (Case 5) $y < y' \leq z < z'$: Due to a $3i312$ pattern,
\begin{equation}\label{eqn5.67} 
u(z'+1) \leq y, \quad u(z'+2) \leq y'.
\end{equation}
As in Figure \ref{fig11}, define
\begin{equation}\label{eqn5.39}
u':=s_{y'}s_{y'+1} \cdots s_{z'}s_{z'+1} s_{y}s_{y+1} \cdots s_{z'-1}s_{z'}.
\end{equation}
Then $u' \leq u$ by applying Lemma \ref{lem5.5} with $K:=\{ z'+1 , z'+2 \}$ and \eqref{eqn5.67}. Also, $u'$ is the minimal representative of $W_Iu'$. Denote $l:= \ell (u')$ and define $(\mathbf{i}, \mathbf{m})$ from the reduced word expression in \eqref{eqn5.39} and $\lambda$ in \eqref{eqn5.22}. We can compute
\begin{equation}\label{eqn5.40}
m_j = 
\begin{cases}
1 & j = z'-y'+3 + z -y(=:j_3) \\
1 & j = l \\
0 & \text{otherwise}
\end{cases}
.
\end{equation}

Write $j_1:= z'-y'+3$ and $j_2:=z'-y+3$. Define $J_{j,k}$ and $K_{j,k}$ for $1\leq j \leq l$ and $1 \leq k \leq m_j$
\begin{equation}\label{eqn5.41}
J_{j,k} = 
\begin{cases}
[j_1,l] & j=j_3 \\
[j_2,l] & j = l
\end{cases}
, \quad K_{j,k}
\begin{cases}
[l]& j =j_3 \\
[j_3+2,l]& j =l
\end{cases}
.
\end{equation}
Then the corresponding LS-monomials $\{\alpha_{j,k}\}$ and $\{\beta_{j,k}\}$ are
\begin{equation}\label{eqn5.42}
\alpha_{j,k} = 
\begin{cases}
\varpi_{z}-\epsilon_y+\epsilon_{z+1} & j = j_3 \\
\varpi_{z'}-\epsilon_{y'}+\epsilon_{z'+1} & j=l 
\end{cases}
, \quad \beta_{j,k}=
\begin{cases}
\varpi_{z}-\epsilon_{y} - \epsilon_{y'} + \epsilon_{z+1}+\epsilon_{z+2} & j = j_3 \\
\varpi_{z'}-\epsilon_{z+2} + \epsilon_{z'+1} & j=l 
\end{cases}
.
\end{equation}

$\bullet$ (Case 6) $y \leq z <  y' \leq z'$: Due to a $3i312$ pattern,
\begin{equation}\label{eqn5.68}
u(z) \geq y'+1, \quad u(z'+1) \leq y, \quad u(z'+2) \leq y'.
\end{equation}
As in Figure \ref{fig12}, define
\begin{equation}\label{eqn5.43}
u':=s_{y'} s_{y'-1} \cdots s_{z+2} s_{z+1} s_{y'+1} s_{y'+2} \cdots s_{z'} s_{z'+1} s_{y}s_{y+1} \cdots s_{z'-1} s_{z'}.
\end{equation}
Then $u' \leq u$ by applying Lemma \ref{lem5.5} with $K:=\{z , z'+1 , z'+2\}$ and \eqref{eqn5.68}. Also, $u'$ is the minimal representative of $W_Iu'$. Denote $l:= \ell (u')$ and define $(\mathbf{i}, \mathbf{m})$ from the reduced word expression in \eqref{eqn5.43} and $\lambda$ in \eqref{eqn5.22}. We can compute
\begin{equation}\label{eqn5.44}
m_j = 
\begin{cases}
1 & j = z'-y+2(=:j_2) \\
1 & j = l \\
0 & \text{otherwise}
\end{cases}
.
\end{equation}

Write $j_1:= z'-z+2$ and $j_3:= y'-y+z'-z+2$. Define $J_{j,k}$ and $K_{j,k}$ for $1\leq j \leq l$ and $1 \leq k \leq m_j$
\begin{equation}\label{eqn5.45}
J_{j,k} = 
\begin{cases}
[l] & j=j_2 \\
[j_3 +1 , l] & j = l
\end{cases}
, \quad K_{j,k}
\begin{cases}
[j_1,l]& j =j_2 \\
[j_2+1,l]& j =l
\end{cases}
.
\end{equation}
Then the corresponding LS-monomials $\{\alpha_{j,k}\}$ and $\{\beta_{j,k}\}$ are
\begin{equation}\label{eqn5.46}
\alpha_{j,k} = 
\begin{cases}
\varpi_{z}-\epsilon_y+\epsilon_{y'+1} & j = j_2 \\
\varpi_{z'}-\epsilon_{y'+1}+\epsilon_{z'+1} & j=l 
\end{cases}
, \quad \beta_{j,k}=
\begin{cases}
\varpi_{z}-\epsilon_{y} + \epsilon_{z+1} & j = j_2 \\
\varpi_{z'}-\epsilon_{z+1} + \epsilon_{z'+1} & j=l 
\end{cases}
.
\end{equation}

From (Case 1) to (Case 6), we constructed two distinct $I$-dominant LS-monomials $\{\alpha_{j,k}\}$ and $\{\beta_{j,k}\}$ with the same weight based on $u'$ and $\lambda$. Let $w'$ be the maximal representative of $W_I u'$. Due to Theorem \ref{thm4.1}, $H^0(X(w'W_J), L_\lambda)^*$ is not $L_I$-multiplicity-free. From $u'\leq u$, we have $w' \leq w$. Using Lemma \ref{lem2.5}, $H^0(X(wW_J), L_\lambda)^*$ is also not $L_I$-multiplicity-free. Due to Theorem \ref{lem2.4}, $X(wW_J)$ is not $L_I$-spherical.
\end{proof}

\begin{figure}
\centering
\begin{subfigure}[t]{0.45\textwidth}
\centering
	\includegraphics[width=\textwidth]{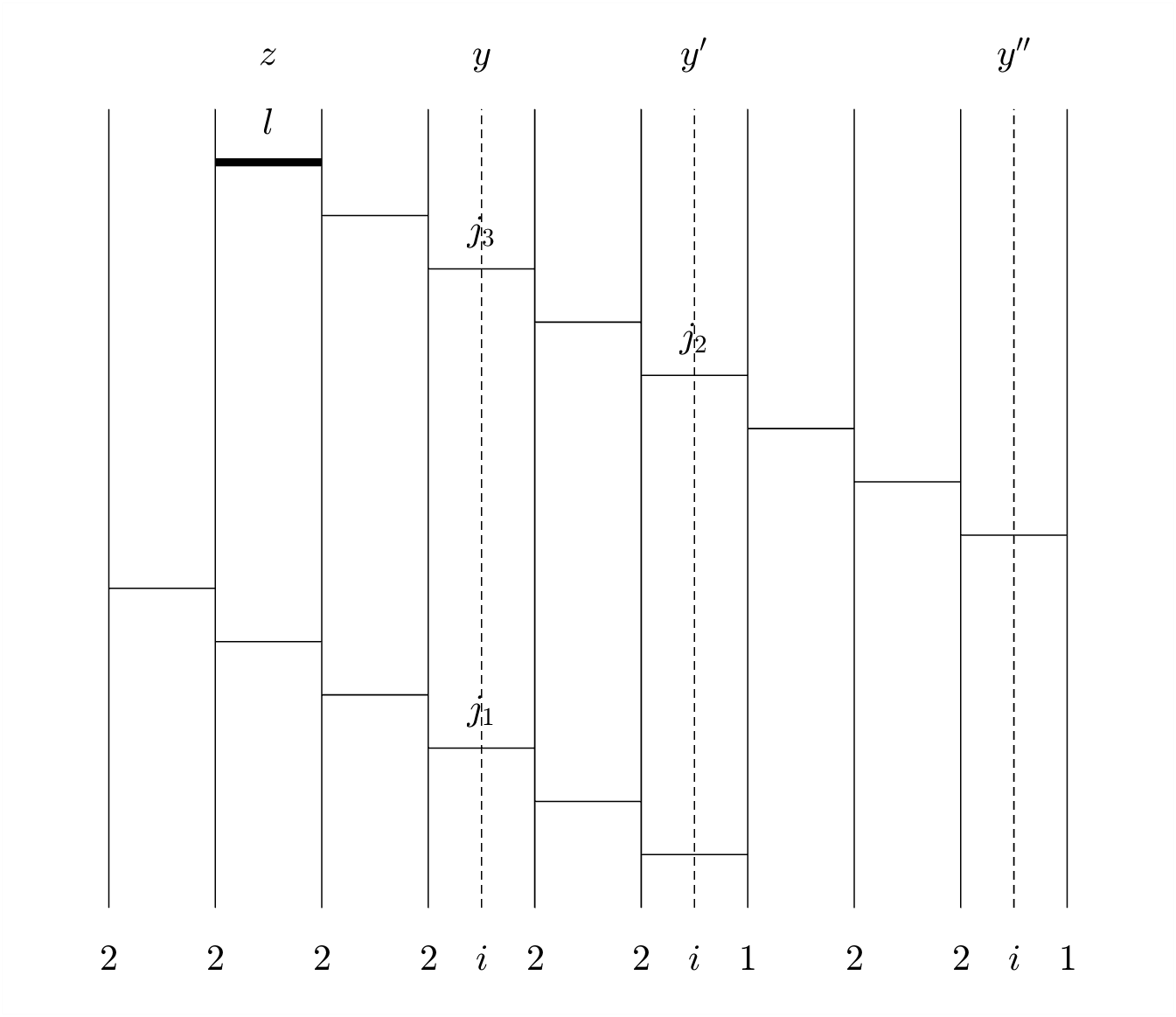}
	\caption{(Case 1) $z \leq y < y' < y''$}
	\label{fig13}
\end{subfigure}
\hfill
\begin{subfigure}[t]{0.45\textwidth}
\centering
\includegraphics[width=\textwidth]{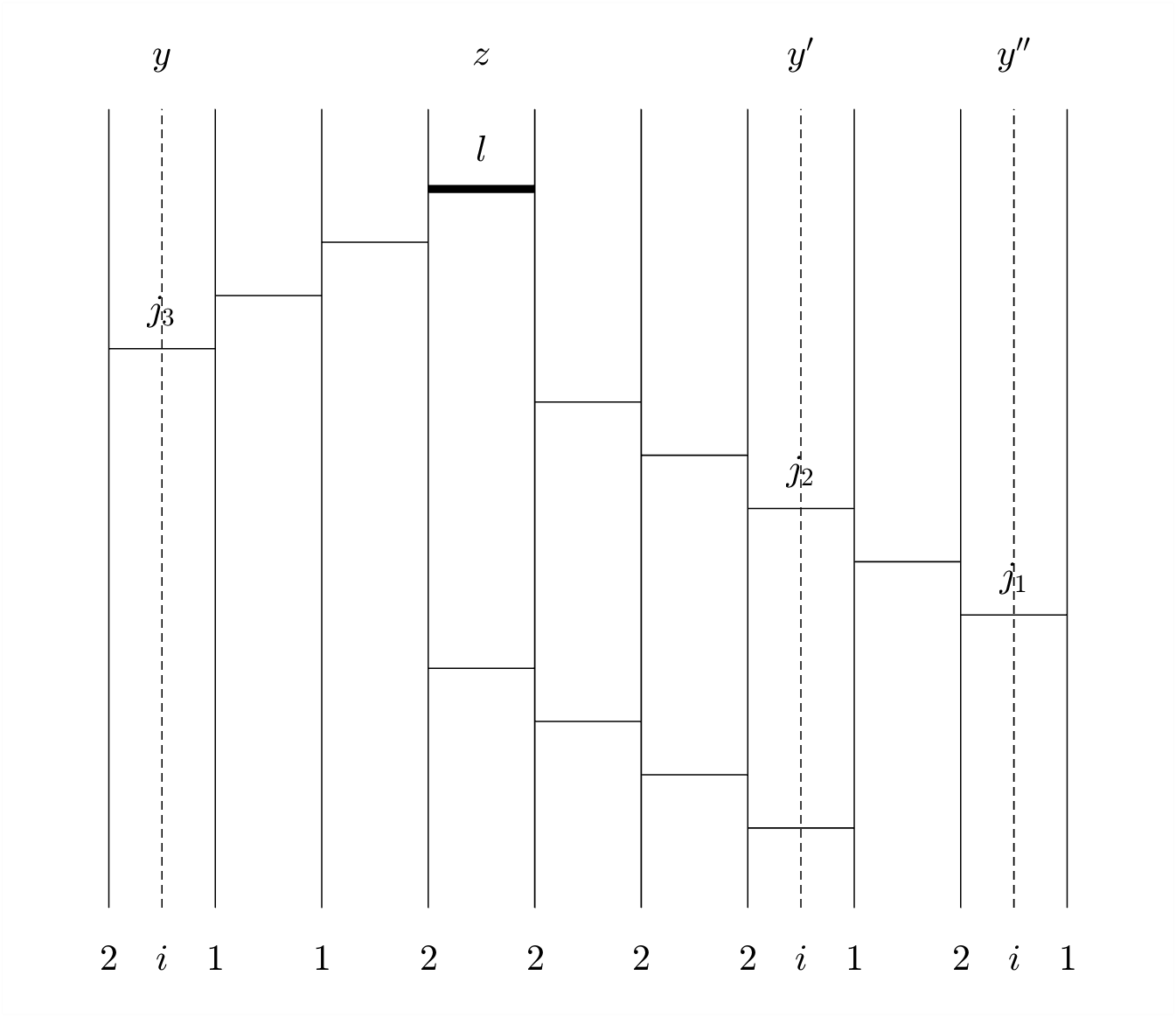}
	\caption{(Case 2) $y < z \leq y' < y''$}
	\label{fig14}
\end{subfigure}
\caption{$2i21i1$ patterns in $u'$}
\label{fig01314}
\end{figure}

\begin{theorem}\label{thm5.4}
Let $G= GL_n$, $W=S_n$, $I, J \subseteq [n-1]$ and $u \in {}^I W^J$. If $[u,I,J]$ has a $2i21i1$ pattern, then $[u,I,J]$ is non-spherical.
\end{theorem}
\begin{proof}
Write $[u,J] = x_1  \cdots x_n$. Recall that $h:=n - |J|$. Since $[u,I,J]$ has a $2i21i1$ pattern, choose $a \leq y < b< c \leq y'' <d$ so that $x_c = x_d <x_a = x_b$ and $y,y'' \in [n-1]-I$. Use $J$-expansion Lemma \ref{lem4.5} so that $1 , \cdots , x_c$ and $x_c +1 , \cdots ,  h$ are equalized. Then $[u,I,J']$ still has a $2i21i1$ pattern. Therefore, we may assume that $[n-1]-J = \{ z\}$ \textit{i.e.} $x_i \in [2]$.

Due to Lemma \ref{lem2.6}, there exists a $y' \in [n-1]-I$ such that $b \leq y' < c$. Use $I$-expansion Lemma \ref{lem4.4} on $I' := [n-1]-\{y,y',y''\}$. Then $[u',I',J]$ still has a $2i21i1$ pattern. Therefore, we may assume that $[n-1]-I = \{y,y',y''\}$. In other words,
\begin{equation}\label{eqn5.47}
[u,I,J] =  \cdots 2i \cdots 2i1\cdots i1\cdots .
\end{equation}
Choose $\lambda \in \Parj$
\begin{equation}\label{eqn5.48}
\lambda := 2 \varpi_z.
\end{equation}
From the relation between $y<y'<y''$ and $z$, we have four cases. Due to symmetry, it is enough to cover the cases $z\leq y < y' <y''$ and $y < z \leq y'  < y''$. 

$\bullet$ (Case 1) $z\leq y < y' <y''$: Due to a $2i21i1$ pattern,
\begin{equation}\label{eqn5.69}
u(z-1)\geq y'+1, \quad u(z) \geq y''+1.
\end{equation}
As in Figure \ref{fig13}, define
\begin{equation}\label{eqn5.49}
u':= s_{y'}s_{y'-1} \cdots s_{z}s_{z-1} s_{y''} s_{y''-1} \cdots s_{z+1}s_z.
\end{equation}
Then $u' \leq u$ by applying Lemma \ref{lem5.5} with $K:=\{z-1, z\}$ and \eqref{eqn5.69}. Also, $u'$ is the minimal representative of $W_Iu'$. Denote $l:= \ell (u')$ and define $(\mathbf{i}, \mathbf{m})$ from the reduced word expression in \eqref{eqn5.49} and $\lambda$ in \eqref{eqn5.48}. We can compute
\begin{equation}\label{eqn5.50}
m_j = 
\begin{cases}
2 & j = l \\
0 & \text{otherwise}
\end{cases}
.
\end{equation}

Write $j_1 := y'-y+1$, $j_2:= y''-z+3$ and $j_3:=y'-y+y''-z+3$. Define $J_{j,k}$ and $K_{j,k}$ for $1\leq j \leq l$ and $1 \leq k \leq m_j$
\begin{equation}\label{eqn5.51}
J_{j,k} = 
\begin{cases}
[l] & j=l , \hsp k=1 \\
[j_3,l] & j = l , \hsp k=2
\end{cases}
, \quad K_{j,k}
\begin{cases}
[j_1,l]& j =l, \hsp k=1 \\
[j_2,l]& j =l, \hsp k=2
\end{cases}
.
\end{equation}
Then the corresponding LS-monomials $\{\alpha_{j,k}\}$ and $\{\beta_{j,k}\}$ are
\begin{subequations}\label{eqn5.52}
\begin{equation}
\alpha_{j,k} = 
\begin{cases}
\varpi_{z}-\epsilon_z - \epsilon_{z-1}+\epsilon_{y'+1}+\epsilon_{y''+1} & j = l , \hsp k=1 \\
\varpi_{z}-\epsilon_{z}+\epsilon_{y+1} & j=l ,\hsp k=2
\end{cases}
, 
\end{equation}
\begin{equation}
\beta_{j,k}=
\begin{cases}
\varpi_{z}-\epsilon_z - \epsilon_{z-1}+\epsilon_{y+1}+\epsilon_{y''+1} & j = l , \hsp k=1 \\
\varpi_{z}-\epsilon_{z}+\epsilon_{y'+1} & j=l ,\hsp k=2
\end{cases}
.
\end{equation}
\end{subequations}

$\bullet$ (Case 2) $y < z \leq y' < y''$: Due to a $2i21i1$ pattern,
\begin{equation}\label{eqn5.70}
u(z-1) \geq y'+1, \quad u(z) \geq y''+1, \quad u(z+1) \leq y.
\end{equation}
As in Figure \ref{fig14}, define
\begin{equation}\label{eqn5.53}
u':= s_{y'} s_{y'-1} \cdots s_{z+1} s_z s_{y''}s_{y''-1} \cdots s_{z+2}s_{z+1}s_y s_{y+1} \cdots s_{z-1}s_z.
\end{equation}
Then $u' \leq u$ by applying Lemma \ref{lem5.5} with $K:=\{z-1, z , z+1\}$ and \eqref{eqn5.70}. Also, $u'$ is the minimal representative of $W_Iu'$. Denote $l:= \ell (u')$ and define $(\mathbf{i}, \mathbf{m})$ from the reduced word expression in \eqref{eqn5.53} and $\lambda$ in \eqref{eqn5.48}. We can compute
\begin{equation}\label{eqn5.54}
m_j = 
\begin{cases}
2 & j = l \\
0 & \text{otherwise}
\end{cases}
.
\end{equation}

Write $j_1:= y'-z+2$, $j_2:= y''-z+2$ and $j_3:= y''+y'-2z+2$. Define $J_{j,k}$ and $K_{j,k}$ for $1\leq j \leq l$ and $1 \leq k \leq m_j$
\begin{equation}\label{eqn5.55}
J_{j,k} = 
\begin{cases}
[l] & j=l , \hsp k=1 \\
\emptyset & j = l , \hsp k=2
\end{cases}
, \quad K_{j,k}
\begin{cases}
[j_1,l]& j =l, \hsp k=1 \\
[j_2,j_3-1] \cup \{l\}& j =l, \hsp k=2
\end{cases}
.
\end{equation}
Then the corresponding LS-monomials $\{\alpha_{j,k}\}$ and $\{\beta_{j,k}\}$ are
\begin{subequations}\label{eqn5.56}
\begin{equation}
\alpha_{j,k} = 
\begin{cases}
\varpi_{z}-\epsilon_y - \epsilon_{z}+\epsilon_{y'+1}+\epsilon_{y''+1} & j = l , \hsp k=1 \\
\varpi_{z} & j=l ,\hsp k=2
\end{cases}
, 
\end{equation}
\begin{equation}
\beta_{j,k}=
\begin{cases}
\varpi_{z}-\epsilon_y +\epsilon_{y''+1} & j = l , \hsp k=1 \\
\varpi_{z}-\epsilon_{z}+\epsilon_{y'+1} & j=l ,\hsp k=2
\end{cases}
.
\end{equation}
\end{subequations}

From (Case 1) to (Case 2), we constructed two distinct $I$-dominant LS-monomials $\{\alpha_{j,k}\}$ and $\{\beta_{j,k}\}$ with the same weight based on $u'$ and $\lambda$. Let $w'$ be the maximal representative of $W_I u'$. Due to Theorem \ref{thm4.1}, $H^0(X(w'W_J), L_\lambda)^*$ is not $L_I$-multiplicity-free. From $u'\leq u$, we have $w' \leq w$. Using Lemma \ref{lem2.5}, $H^0(X(wW_J), L_\lambda)^*$ is also not $L_I$-multiplicity-free. Due to Theorem \ref{lem2.4}, $X(wW_J)$ is not $L_I$-spherical.
\end{proof}

\section{Sphericality of Grassmannians}\label{sec6}

In \cite{MR1395667}, a permutation $w \in S_n$ is \textbf{bigrassmannian}\footnote{By convention, we include the identity as a bigrassmannian permutation.} if $w$ has a unique descent and so does $w^{-1}$. Then there exist $p , q \in [n-1]$ and $r \in \Z_{>0}$ such that
\begin{equation}
w = (s_{p} \cdots s_{q+r-1})(s_{p-1}\cdots s_{q+r-2} )\cdots (s_{p-r+1} \cdots s_q).
\end{equation}
$w$ has a wiring diagram of rectangle shape; see Figure \ref{fig02021}. Bigrassmannian permutations can be characterized as permutations avoiding the patterns $321$, $2143$, $3142$ and $2413$ so that there are no multiple descents.

A permutation $w \in S_n$ is \textbf{standard Coxeter} if there exists a reduced word expression of $w$ using each $s_i$ at most once. Equivalently, $w$ avoids $321$ and $3412$ patterns due to \cite[Theorem 4.3]{MR2333139}, where $w$ is called \textit{Boolean}.

In this section, we discuss how to reduce the classification problem of spherical $[u,I,J]$ to a special case when $u$ is bigrassmannian. To do this, we generalize the concept of bigrassmannian permutations and standard Coxeter permutations.

\subsection{Generalization of Coxeter elements} 

Let $G = GL_n$ and $W = S_n$. Consider a minimal sequence of $n$. Recall that it can be uniquely expressed by $[u,I,J]$ where $I, J \subseteq [n-1]$ and $u \in {}^I W^J$. To compute $u$, we first eliminate $i$ inserted in the sequence. Read off the numbers of the sequence in weakly increasing order. If there are entries with the same number, read from the left to the right. As we read the numbers, replace them with $1,2, \cdots , n$ and we have the one-line notation of $u^{-1}$. For instance, consider $13i4i223i4$ which is a minimal sequence of $7$. Remove $i$ and we have $1342234$. Replace them with $1,2, \cdots ,7$ and we have $u^{-1} = 1462357$.

\begin{lemma}\label{lem6.0}
Let $G = GL_n$, $W= S_n$, $I,J \subseteq [n-1]$ and $u \in {}^I W^J$. Suppose $[u,I,J]$ avoids $321$, $3412$, $3i312$, $231i1$ and $2i21i1$ patterns. Then $u$ avoids permutation patterns as follows.
\begin{equation}\label{eqn6.0}
321, \quad 246135, \quad 351624, \quad 415263.
\end{equation}
\end{lemma}
\begin{proof}
The inverses of \eqref{eqn6.0} are $321$, $415263$, $351624$ and $246135$, respectively. Therefore, it is sufficient to prove that $u^{-1}$ avoids \eqref{eqn6.0}.

Suppose $u^{-1}$ has a $321$ pattern. After replacement, $[u,I,J]$ still has a $321$ pattern, contradicting the assumption.

Suppose $u^{-1}$ has a $246135$ pattern. From the underlined numbers of $2\underline{4613}5$, it has a $3412$ pattern. Since $[u,I,J]$ does not have a $3412$ pattern, at least two numbers among $4613$ should be equalized in $[u,I,J]$. This leads to a contradiction. For instance, if $1$ and $3$ are equalized, then $\underline{2}46\underline{13}5$ should be all equalized, contradicting $u \in {}^I W^J$. 

Suppose $u^{-1}$ has a $351624$ pattern. Again, it has a $3412$ pattern in $\underline{351}6\underline{2}4$. Then at least two numbers among $3512$ should be equalized in $[u,I,J]$. $3$ and $5$ cannot be equalized because $4$ is on the right. Therefore, $1$ and $2$ should be equalized in $35\underline{1}6 \underline{2}4$. Due to $6$ between $1$ and $2$, there is an $i$ inserted between equalized $1$ and $2$ in $[u,I,J]$. In other words, there is a $231i1$ pattern in $[u,I,J]$, leading to a contradiction.

Suppose $u^{-1}$ has a $415263$ pattern. Again, it has a $3412$ pattern in $\underline{4}1\underline{52}6\underline{3}$. If $4$ and $5$ are equalized in $[u,I,J]$, there is an $i$ inserted between them. There is also an $i$ between $2$ and $3$. Therefore, $[u,I,J]$ has $231i1$, $3i312$ or $2i21i1$ patterns, leading to a contradiction.
\end{proof}

\begin{figure}
\centering
\begin{subfigure}[t]{\textwidth}
\centering
\includegraphics[width=\textwidth]{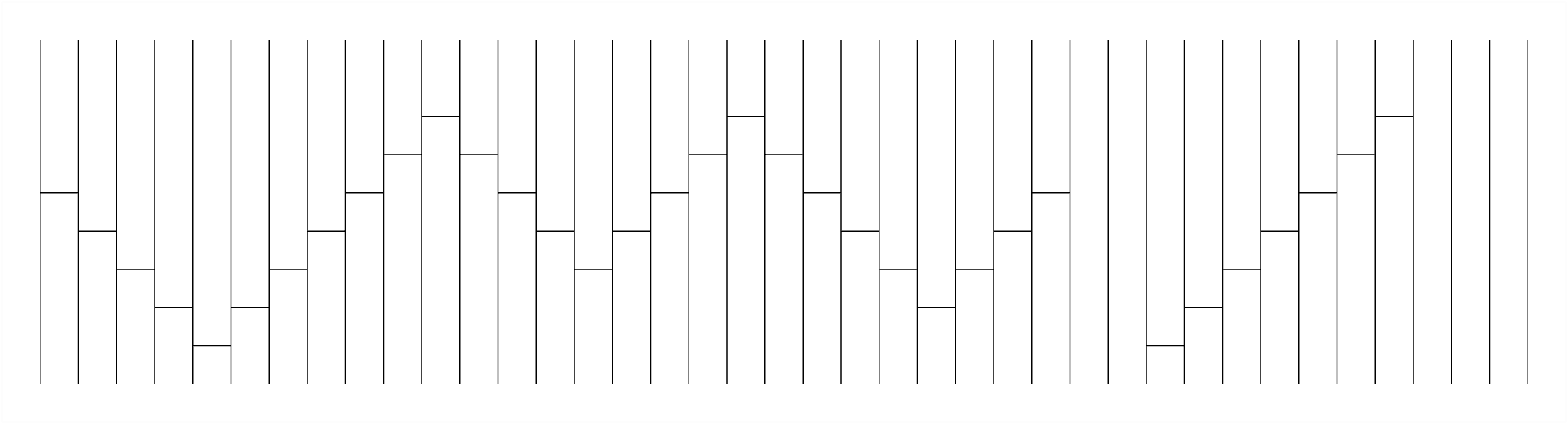}
\caption{A $321$ and $3412$ avoiding permutation}
\label{fig27}
\end{subfigure}
\hfill
\begin{subfigure}[t]{\textwidth}
\centering
\includegraphics[width=\textwidth]{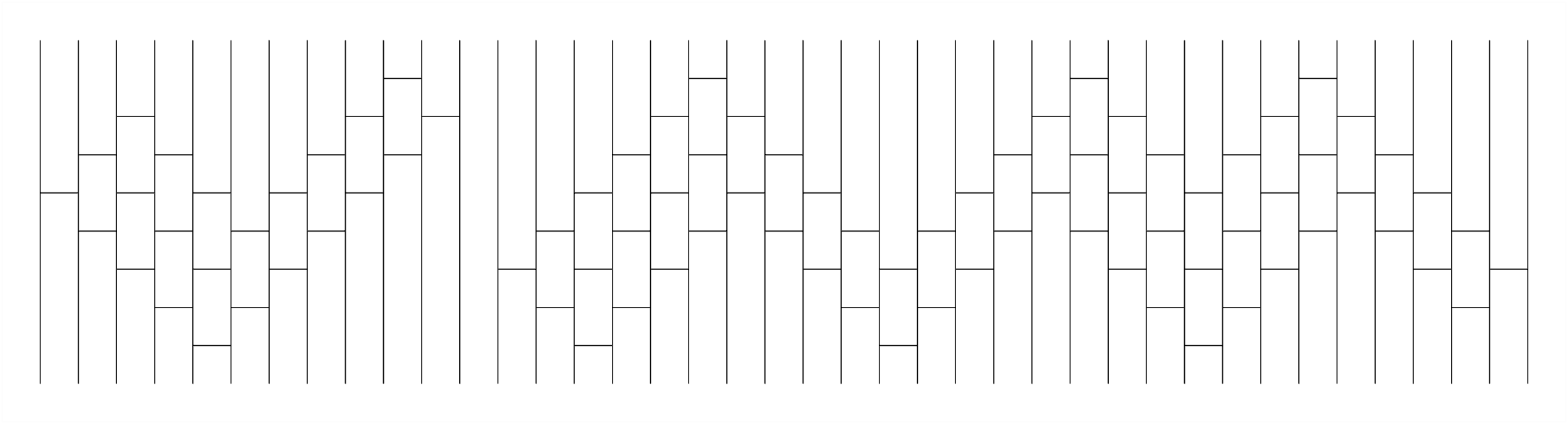}
\caption{A $321$, $246135$, $351624$ and $415263$ avoiding permutation}
\label{fig28}
\end{subfigure}
\hfill
\begin{subfigure}[t]{\textwidth}
\centering
\includegraphics[width=\textwidth]{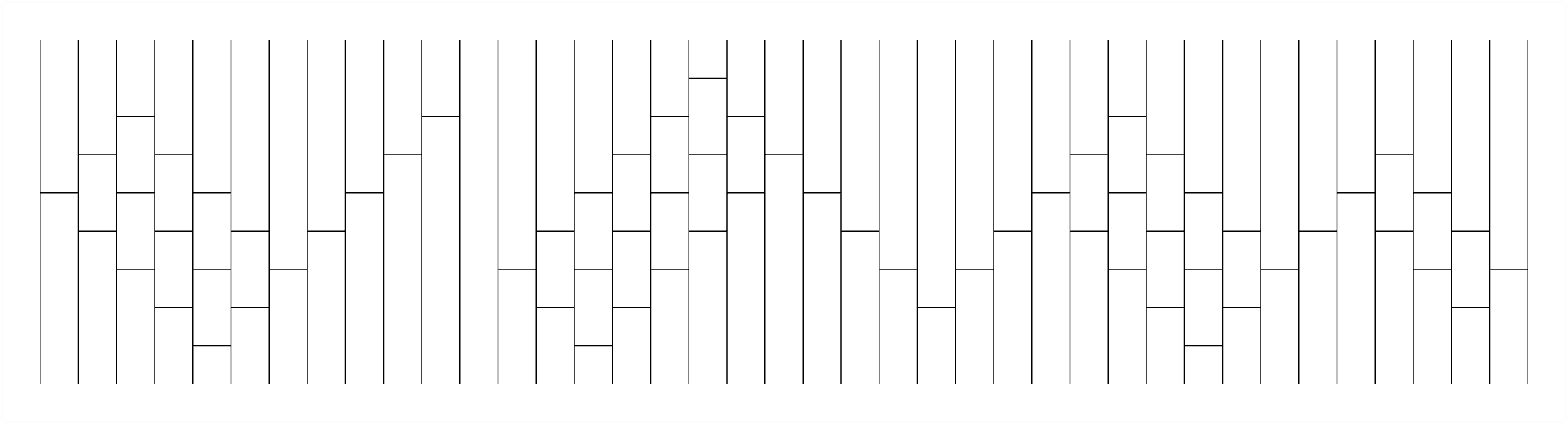}
\caption{A $321$, $246135$, $351624$, $356124$, $415263$ and $451623$ avoiding permutation}
\label{fig29}
\end{subfigure}
\caption{Generalizing bigrassmannian and standard Coxeter permutations}
\label{fig0272829}
\end{figure}

For each $w \in S_n$ and from the one-line notation of $w = w_1 w_2 \cdots w_n$, we do the following reduction. Use parentheses to wrap around some consecutive entries with consecutive numbers. For instance, from the one-line notation of the permutation in Figure \ref{fig28}, let $w$ be defined
\begin{multline}\label{eqn7.1}
(6,7,8)(1,2,3)(9,10,11,12)(4,5)(16,17,18,19)(24,25)(13,14,15) \\
(20,21)(26,27)(32,33,34)(22,23)(28)(35)(39,40)(29,30,31)(36,37,38). 
\end{multline} 
Read off parentheses in increasing order while replacing them with $1,2,3, \cdots$ to construct another permutation $\overline{w}$. For instance, from \eqref{eqn7.1}, replace $(1,2,3)$ with $1$, $(4,5)$ with $2$ and so on. Then we have
\begin{equation}\label{eqn7.2}
\overline{w} = 3 \hsp \hsp 1 \hsp \hsp 4 \hsp \hsp 2 \hsp \hsp 6 \hsp \hsp 9 \hsp  \hsp 5 \hsp \hsp 7 \hsp \hsp 10 \hsp \hsp 13 \hsp \hsp 8 \hsp \hsp 11 \hsp \hsp 14 \hsp \hsp 16 \hsp \hsp 12 \hsp \hsp 15.
\end{equation} 
If such reduction is no longer possible, $\overline{w}$ is called the \textbf{minimal reduction} of $w$. For instance, \eqref{eqn7.2} is the minimal reduction of \eqref{eqn7.1}.

Suppose $\overline{w} \in S_m$. In the process of reduction, recall that for each $i \in [m]$, the $i$th parentheses were replaced with $\overline{w} (i)$. Record the number of entries between parentheses by $r_{w , \overline{w}}(i)$. For instance, if $w$ is given as \eqref{eqn7.1}, the first parentheses $(6,7,8)$ have three entries. Therefore, $r_{w, \overline{w}} (1) = 3$. Finishing this process, we get
\begin{equation}
r_{w, \overline{w}} = 3 \hsp \hsp 3 \hsp \hsp 4 \hsp \hsp 2 \hsp \hsp 4 \hsp \hsp 2 \hsp \hsp 3 \hsp \hsp 2 \hsp \hsp 2 \hsp \hsp 3 \hsp \hsp 2 \hsp \hsp 1 \hsp \hsp 1 \hsp \hsp 2 \hsp \hsp 3 \hsp \hsp 3. 
\end{equation}

\begin{lemma}\label{lem4.9}
Let $w \in S_n$. If $w$ avoids permutation patterns
\begin{equation}
321, \quad 246135, \quad 351624, \quad 415263,
\end{equation}
then the minimal reduction $\overline{w}$ avoids $321$ and $3412$ permutation patterns.
\end{lemma}
\begin{proof}
Write the one-line notation $\overline{w} = \overline{w}_1  \overline{w}_2  \cdots  \overline{w}_m$. First of all, $\overline{w}$ avoids $321$ since three numbers corresponding to $321$ are not in the same parentheses in $w$. Suppose $\overline{w}$ has a $3412$ pattern. In other words, there exist $a < b < c < d $ such that $\overline{w}_c < \overline{w}_d < \overline{w}_a < \overline{w}_b$. Assume that $(d-a)$ is minimal compared to other $3412$ patterns which appear.

Consider $a$ and $b$. Since $\overline{w}$ cannot be reduced any more, either $b-a > 1$ or $\overline{w}_b - \overline{w}_a >1$. If $b-a > 1$, there exists an integer $a < e <b$. Since there is no $321$ pattern and $(d-a)$ is minimal, $\overline{w}_e < \overline{w}_c < \overline{w}_d < \overline{w}_a < \overline{w}_b$. On the other hand, if $\overline{w}_b - \overline{w}_a >1$, then there exists an $e$ such that $\overline{w}_a < \overline{w}_e < \overline{w}_b$. Again, since there is no $321$ pattern and $(d-a)$ is minimal, $a<b<c<d <e$.

For $c$ and $d$, we again have two cases using $f$. The conclusion is that there exists a $246135$, $351624$ or $415263$ pattern in $\overline{w}$. Such a pattern also occurs in $w$, leading to a contradiction. Hence, $\overline{w}$ avoids $321$ and $3412$ patterns. 
\end{proof}

\begin{figure}
\centering
\begin{subfigure}[t]{\textwidth}
\centering
\includegraphics[width=\textwidth]{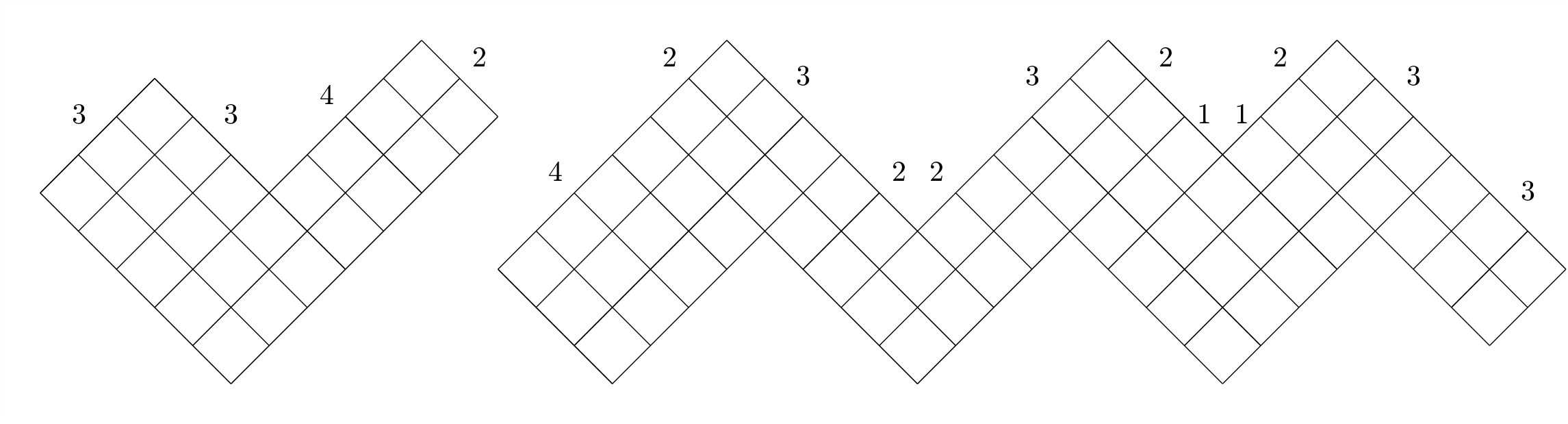}
\caption{A $321$, $246135$, $351624$ and $415263$ avoiding permutation}
\label{fig30}
\end{subfigure}
\hfill
\begin{subfigure}[t]{\textwidth}
\centering
\includegraphics[width=\textwidth]{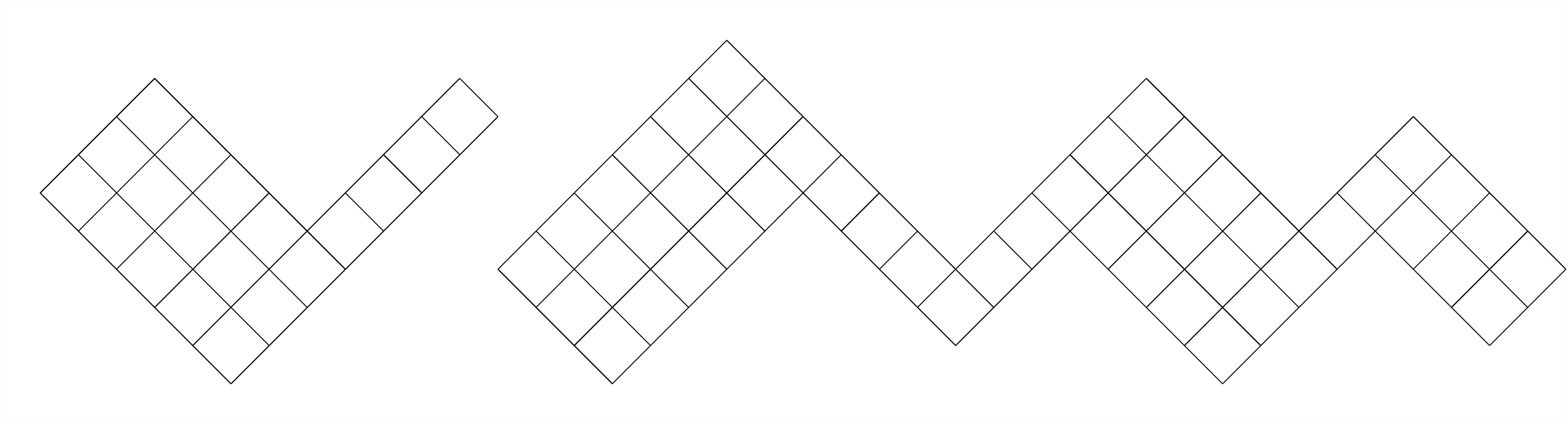}
\caption{A $321$, $246135$, $351624$, $356124$, $415263$ and $451623$ avoiding permutation}
\label{fig31}
\end{subfigure}
\caption{Deformation of wiring diagrams}
\label{fig03031}
\end{figure}

The relation between $w$ and $\overline{w}$ enables us to easily draw a picture of $w$ which avoids patterns \eqref{eqn6.0}. For instance, the wiring diagram of $w$ in Figure \ref{fig28} can be regarded as Figure \ref{fig30}. In Figure \ref{fig30}, there is a zigzag boundary curve composed of line segments of length $3,3,4,2,4,2,3,2,2,3,2,1,1,2,3,3$. These lengths match the number of entries between parentheses in \eqref{eqn7.1} \textit{i.e.} the sequence $r_{w , \overline{w}}$.

For each $w \in S_n$, the \textbf{support} of $w$ is
\begin{equation}
{\sf supp} (w) = \{ i \in [n] \mid w(i) \neq i\}.
\end{equation}
Also, if $w$ is not the identity, define
\begin{equation}
{\sf mins}(w) = \min \{ i \in [n] \mid w(i) \neq i \}, \quad {\sf maxs }(w) = \max \{ i \in [n] \mid w(i) \neq i \}.
\end{equation}

For each bigrassmannian permutation $w$ which is not the identity, define
\begin{equation}
{\sf sw}(w) = w({\sf mins}(w)) - {\sf mins}(w) , \quad {\sf se}(w) = {\sf maxs}(w) - w({\sf maxs}(w)) 
\end{equation}
Conversely, a bigrassmannian permutation $w$ which is not the identity can be uniquely determined by ${\sf mins}(w)$, ${\sf maxs}(w)$, ${\sf sw}(w)$ and ${\sf se}(w)$. Indeed, if
\begin{equation}
{\sf mins}(w) = p , \quad {\sf maxs}(w) = q , \quad {\sf sw}(w) = a , \quad {\sf se}(w) = b,
\end{equation}
then $q - p = a+ b -1$ and
\begin{equation}
w = (s_{p+a-1} s_{p+a-2} \cdots s_p)(s_{p+a} s_{p+a-1} \cdots s_{p+1}) \cdots (s_{q-1} s_{q-2} \cdots s_{q-a}).
\end{equation}

\begin{lemma}\label{lem8.1}
Let $w \in S_n$. If $w$ is not the identity and avoids permutation patterns
\begin{equation}\label{eqn8.1}
321, \quad 246135, \quad 351624, \quad 356124, \quad 415263, \quad 451623,
\end{equation}
then for some $k \in \Z_{> 0}$, there exist $w^{(1)} , \cdots , w^{(k)} \in S_n$ satisfying the following conditions.
\begin{itemize}
\item $w = w^{(1)}  \cdots  w^{(k)}$.
\item For each $1 \leq i \leq k$, $w^{(i)}$ is a bigrassmannian permutation which is not the identity.
\item For each $1 \leq i < j \leq k$,
\begin{equation}
{\sf maxs}(w^{(i)}) \leq {\sf mins}(w^{(j)}) \quad \text{or} \quad {\sf maxs}(w^{(j)}) \leq {\sf mins}(w^{(i)}).
\end{equation}
\end{itemize} 
\end{lemma}
\begin{proof}
Let ${\sf maxs}(w) = N$. We use induction on $N$. By assumption, $w(N) \neq N$ whereas $w(i)= i$ for $i > N$. We may identify $w$ with a permutation in $S_N$ to find the decomposition. Therefore, assume $N = n$.

Let $\overline{w} \in S_m$ be the minimal reduction of $w$. From Lemma \ref{lem4.9}, $\overline{w}$ is a standard Coxeter element. If we choose a reduced word $(i_1 , \cdots , i_l)$ of $\overline{w}$, then $i_j \neq i_k$ if $j \neq k$. From $w(N) \neq N$, $w(m) \neq m$. Therefore, there uniquely exists $a \in [l]$ such that $i_a = m-1$. Then we have three cases as follows. In the first case, there is no $b \in [l]$ such that $i_b = m-2$. Otherwise, there uniquely exists $b \in [l]$ such that $i_b = m-2$. The second case is $a < b$ and the last case is $a > b$.

$\bullet$ (Case 1) Assume that there is no $b \in [l]$ such that $i_b = m-2$. Then $\overline{w}(m-1) = m$ and $\overline{w} (m) = m-1$. If we write $r_{w , \overline{w}} (m-1) = p$ and $r_{w , \overline{w}} (m) = q$, then $w(i) = i+q$ for $i \in [n-p-q+1 , n-q]$ and $w (i) = i-p$ for $i \in [n-q+1 , n]$. Choose a bigrassmannian permutation $u$ such that
\begin{equation}
{\sf mins}(u) = n-p-q+1 , \quad {\sf maxs}(u) = n , \quad {\sf sw}(u) = q, \quad {\sf se}(u) = p.
\end{equation}
Let $v:= w u^{-1}$. If $v$ is the identity, we are done. Otherwise,
\begin{equation}\label{eqn8.2}
v (i) = 
\begin{cases}
w(i) & i \in [n-p-q] \\
i & i \in [n-p-q+1 , n]
\end{cases}
.
\end{equation}
Therefore, ${\sf maxs} (v) \leq n-p-q$ and $\ell (w) = \ell (u) + \ell (v)$. 

Suppose there exist $j_1 < \cdots < j_t$ such that $v$ has one of the permutation patterns in \eqref{eqn8.1} at $v(j_1) , \cdots , v(j_t)$. From \eqref{eqn8.1}, $v(j_t)$ is not the maximum among $v(j_1) , \cdots , v(j_t)$. From \eqref{eqn8.2}, $j_t \leq n-p-q$. This implies that $w$ also has one of the permutation patterns in \eqref{eqn8.1} at $w(j_1) , \cdots , w(j_t)$, leading to a contradiction. Hence, $v$ also avoids permutation patterns in \eqref{eqn8.1}. Due to the induction hypothesis on ${\sf maxs}(w) = N = n$, $v$ has a decomposition $v^{(1)} , \cdots , v^{(k)}$ as desired. So does $w$ by including $v^{(1)} , \cdots , v^{(k)} , u$, completing the induction steps.

$\bullet$ (Case 2) Assume that there uniquely exists $b \in [l]$ such that $i_b = m-2$ and $a < b$. Suppose there is no $c \in [l]$ such that $i_c = m-3$. Then $\overline{w} (m-1)= m-2$ and $\overline{w} (m) = m-1$, contradicting the assumption that $\overline{w}$ is the minimal reduction. Therefore, there uniquely exists such $c \in [l]$. If $b < c$, then again, $\overline{w} (m-1)= m-2$ and $\overline{w} (m) = m-1$. Hence, $b >c$ and we have
\begin{equation}\label{eqn8.3}
\overline{w}(m-2) = m , \quad \overline{w} (m-1 ) \leq m-3 , \quad \overline{w} (m) = m-1.
\end{equation}

Next, suppose $r_{w , \overline{w}} (m-2) \geq 2$ and $r_{w , \overline{w}} (m-1) \geq 2$. From \eqref{eqn8.3}, there exists $r < m-2$ such that $\overline{w} (r) = m-2$. Then
\begin{equation}
r < m-2 < m-1 < m , \quad \overline{w} (m-1) < \overline{w} (r) < \overline{w} (m) < \overline{w} (m-2)
\end{equation}
with $r_{w , \overline{w}} (m-2) \geq 2$ and $r_{w , \overline{w}} (m-1) \geq 2$ imply that $w$ has a $356124$ permutation pattern, leading to a contradiction. Therefore, we have two subcases as follows. We have $r_{w , \overline{w}} (m-2) = 1$ or $r_{w , \overline{w}} (m-1) = 1$.

$\bullet$ (Case 2-1) Assume that $r_{w , \overline{w}} (m-2) = 1$. If we write $r_{w , \overline{w}} (m-1) = p$ and $r_{w , \overline{w}} (m) = q$, then $w (n-p-q) = n$ and $w (i) = i-1$ for $i \in [n-q+1 , n]$ due to \eqref{eqn8.3}. Choose a bigrassmannian permutation $u$ such that
\begin{equation}
{\sf mins}(u) = n-q , \quad {\sf maxs}(u) = n , \quad {\sf sw}(u) = q, \quad {\sf se}(u) = 1.
\end{equation}
Let $v:= u^{-1} w$. If $v$ is the identity, we are done. Otherwise,
\begin{equation}\label{eqn8.4}
v (i) = 
\begin{cases}
w(i) & i \in [n-p-q-1] \\
n-q & i = n-p -q \\
w(i) & i \in [n-p -q +1 , n-q] \\
i & i \in [n-q +1 , n]
\end{cases}
.
\end{equation}
Therefore, ${\sf maxs} (v) \leq n-q$ and $\ell (w) = \ell (u) + \ell (v)$. 

Suppose there exist $j_1 < \cdots < j_t$ such that $v$ has one of the permutation patterns in \eqref{eqn8.1} at $v(j_1) , \cdots , v(j_t)$. From \eqref{eqn8.1}, $v(j_t)$ is not the maximum among $v(j_1) , \cdots , v(j_t)$. From \eqref{eqn8.4}, $j_t \leq n-q$. This implies that $w$ also has one of the permutation patterns in \eqref{eqn8.1} at $w(j_1) , \cdots , w(j_t)$, leading to a contradiction. Hence, $v$ also avoids permutation patterns in \eqref{eqn8.1}. Due to the induction hypothesis on ${\sf maxs}(w) = N = n$, $v$ has a decomposition $v^{(1)} , \cdots , v^{(k)}$ as desired. So does $w$ by including $u , v^{(1)} , \cdots , v^{(k)}$, completing the induction steps.

$\bullet$ (Case 2-2) Assume that $r_{w , \overline{w}} (m-1) = 1$. If we write $r_{w , \overline{w}} (m-2) = p$ and $r_{w , \overline{w}} (m) = q$, then $w(i) = i+1+q$ for $i \in [n-p-q , n-q-1]$, $w(n-q) <n-p-q$ and $w(i) = i-p$ for $i \in [n-q+1 , n]$ due to \eqref{eqn8.3}. Choose a bigrassmannian permutation $u$ such that
\begin{equation}
{\sf mins}(u) = n-p-q , \quad {\sf maxs}(u) = n , \quad {\sf sw}(u) = q+1, \quad {\sf se}(u) = p.
\end{equation}
Let $v:= wu^{-1}$. If $v$ is the identity, we are done. Otherwise,
\begin{equation}\label{eqn8.5}
v (i) = 
\begin{cases}
w(i) & i \in [n-p-q-1] \\
w(n-q) & i = n-p -q \\
i & i \in [n-p-q+1 , n]
\end{cases}
.
\end{equation}
Therefore, ${\sf maxs} (v) \leq n-p-q$ and $\ell (w) = \ell (u) + \ell (v)$. 

Suppose there exist $j_1 < \cdots < j_t$ such that $v$ has one of the permutation patterns in \eqref{eqn8.1} at $v(j_1) , \cdots , v(j_t)$. From \eqref{eqn8.1}, $v(j_t)$ is not the maximum among $v(j_1) , \cdots , v(j_t)$. From \eqref{eqn8.5}, $j_t \leq n-p-q$. If $j_t \neq n-p-q$, then $w$ also has one of the permutation patterns in \eqref{eqn8.1} at $w(j_1) , \cdots , w(j_t)$. If $j_t = n-p-q$, then $w$ also has one of the permutation patterns in \eqref{eqn8.1} at $w(j_1) , \cdots , w(j_{t-1}), w(n-q)$. Either way, it leads to a contradiction. Hence, $v$ also avoids permutation patterns in \eqref{eqn8.1}. Due to the induction hypothesis on ${\sf maxs}(w) = N = n$, $v$ has a decomposition $v^{(1)} , \cdots , v^{(k)}$ as desired. So does $w$ by including $v^{(1)} , \cdots , v^{(k)} , u$, completing the induction steps.

$\bullet$ (Case 3) Assume that there uniquely exists $b \in [l]$ such that $i_b = m-2$ and $a > b$. Suppose there is no $c \in [l]$ such that $i_c = m-3$. Then $\overline{w} (m-2)= m-1$ and $\overline{w} (m-1) = m$, contradicting $\overline{w}$ is the minimal reduction. Therefore, there uniquely exists such $c \in [l]$. If $b > c$, then again, $\overline{w} (m-2)= m-1$ and $\overline{w} (m-1) = m$. Hence, $b <c$ and we have
\begin{equation}\label{eqn8.6}
\overline{w}^{-1}(m-2) = m , \quad \overline{w}^{-1} (m-1 ) \leq m-3 , \quad \overline{w}^{-1} (m) = m-1.
\end{equation}
Let $\overline{w} (s) = m-1$ where $s \leq m-3$. Suppose $r_{w , \overline{w}} (s) \geq 2$ and $r_{w , \overline{w}} (m) \geq 2$. From \eqref{eqn8.6}, $\overline{w} (m-2) < m-2$. Then
\begin{equation}
s < m-2 < m-1 < m, \quad \overline{w} (m-2) < \overline{w} (m) < \overline{w} (s) < \overline{w} (m-1)
\end{equation}
with $r_{w , \overline{w}} (s) \geq 2$ and $r_{w , \overline{w}} (m) \geq 2$ imply that $w$ has a $451623$ permutation pattern, leading to a contradiction. Therefore, we have two subcases as follows. We have $r_{w , \overline{w}} (m) = 1$ or $r_{w , \overline{w}} (s) = 1$.

$\bullet$ (Case 3-1) Assume that $r_{w , \overline{w}} (m) = 1$. If we write $r_{w , \overline{w}} (m-1) = p$, then $w(i) = i+1$ for $i \in [n-p , n-1]$ and $w (n) <n-p$ due to \eqref{eqn8.6}. Choose a bigrassmannian permutation $u$ such that
\begin{equation}
{\sf mins}(u) = n-p , \quad {\sf maxs}(u) = n , \quad {\sf sw}(u) = 1, \quad {\sf se}(u) = p.
\end{equation}
Let $v:= wu^{-1}$. If $v$ is the identity, we are done. Otherwise,
\begin{equation}\label{eqn8.7}
v (i) = 
\begin{cases}
w(i) & i \in [n-p-1] \\
w(n) & i = n-p\\
i & i \in [n-p+1 , n]
\end{cases}
.
\end{equation}
Therefore, ${\sf maxs} (v) \leq n-p$ and $\ell (w) = \ell (u) + \ell (v)$. 

Suppose there exist $j_1 < \cdots < j_t$ such that $v$ has one of the permutation patterns in \eqref{eqn8.1} at $v(j_1) , \cdots , v(j_t)$. From \eqref{eqn8.1}, $v(j_t)$ is not the maximum among $v(j_1) , \cdots , v(j_t)$. From \eqref{eqn8.7}, $j_t \leq n-p$. If $j_t \neq n-p$, then $w$ also has one of the permutation patterns in \eqref{eqn8.1} at $w(j_1) , \cdots , w(j_t)$. If $j_t = n-p$, then $w$ also has one of the permutation patterns in \eqref{eqn8.1} at $w(j_1) , \cdots , w(j_{t-1}), w(n)$. Either way, it leads to a contradiction. Hence, $v$ also avoids permutation patterns in \eqref{eqn8.1}. Due to the induction hypothesis on ${\sf maxs}(w) = N = n$, $v$ has a decomposition $v^{(1)} , \cdots , v^{(k)}$ as desired. So does $w$ by including $v^{(1)} , \cdots , v^{(k)} , u$, completing the induction steps.

$\bullet$ (Case 3-2) Assume that $r_{w , \overline{w}} (s) = 1$. If we write $r_{w , \overline{w}} (m-1) = p$ and $r_{w , \overline{w}} (m) = q$, then $w(i) = i+q$ for $i \in [n-p-q+1 , n-q]$, $w(i) = i-p-1$ for $i \in [n-q+1 , n]$ and $w^{-1} (n-p) < n-p-q$ due to \eqref{eqn8.6}. Choose a bigrassmannian permutation $u$ such that
\begin{equation}
{\sf mins}(u) = n-p-q , \quad {\sf maxs}(u) = n , \quad {\sf sw}(u) = q, \quad {\sf se}(u) = p+1.
\end{equation}
Let $v:= u^{-1} w$. If $v$ is the identity, we are done. Otherwise,
\begin{equation}\label{eqn8.8}
v (i) = 
\begin{cases}
i & i \in [n-p-q+1, n] \\
n-p-q & i = w^{-1}(n-p)\\
w(i) & \text{otherwise}
\end{cases}
.
\end{equation}
Therefore, ${\sf maxs} (v) \leq n-p-q$ and $\ell (w) = \ell (u) + \ell (v)$. 

Suppose there exist $j_1 < \cdots < j_t$ such that $v$ has one of the permutation patterns in \eqref{eqn8.1} at $v(j_1) , \cdots , v(j_t)$. From \eqref{eqn8.1}, $v(j_t)$ is not the maximum among $v(j_1) , \cdots , v(j_t)$. From \eqref{eqn8.8}, $j_t \leq n-p-q$. This implies that $w$ also has one of the permutation patterns in \eqref{eqn8.1} at $w(j_1) , \cdots , w(j_t)$, leading to a contradiction. Hence, $v$ also avoids permutation patterns in \eqref{eqn8.1}. Due to the induction hypothesis on ${\sf maxs}(w) = N = n$, $v$ has a decomposition $v^{(1)} , \cdots , v^{(k)}$ as desired. So does $w$ by including $u, v^{(1)} , \cdots , v^{(k)} $, completing the induction steps.
\end{proof}

\begin{lemma}\label{lem4.11}
Let $w \in S_n$. If $w$ is not the identity and avoids permutation patterns
\begin{equation}
321, \quad 246135, \quad 351624, \quad 356124, \quad 415263, \quad 451623, \quad 2561374, \quad 4157236,
\end{equation}
then for some $k \in \Z_{> 0}$, there exist $w^{(1)} , \cdots , w^{(k)} \in S_n$ satisfying the following conditions.
\begin{itemize}
\item $w = w^{(1)}  \cdots  w^{(k)}$.
\item For each $1 \leq i \leq k$, $w^{(i)}$ is a bigrassmannian permutation which is not the identity.
\item For each $1 \leq i < j \leq k$,
\begin{equation}
{\sf maxs}(w^{(i)}) \leq {\sf mins}(w^{(j)}) \quad \text{or} \quad {\sf maxs}(w^{(j)}) \leq {\sf mins}(w^{(i)}).
\end{equation}
\item If for $1 \leq p < q < r \leq k$,
\begin{equation}
{\sf maxs}(w^{(p)}) = {\sf mins}(w^{(q)}) , \quad {\sf maxs}(w^{(q)}) = {\sf mins}(w^{(r)}),
\end{equation}
then
\begin{equation}
{\sf se}(w^{(q)}) = 1 \quad \text{or} \quad {\sf sw}(w^{(q)}) \leq 2.
\end{equation}
\item If for $1 \leq p < q < r \leq k$,
\begin{equation}
{\sf mins}(w^{(p)}) = {\sf maxs}(w^{(q)}) , \quad {\sf mins}(w^{(q)}) = {\sf maxs}(w^{(r)}),
\end{equation}
then
\begin{equation}
{\sf sw}(w^{(q)}) = 1 \quad \text{or} \quad {\sf se}(w^{(q)}) \leq 2.
\end{equation}
\end{itemize}
\end{lemma}
\begin{proof}
From Lemma \ref{lem8.1}, construct $w^{(1)} , \cdots , w^{(k)}$. Suppose there exist $1 \leq p < q < r \leq k$ such that
\begin{equation}
{\sf maxs}(w^{(p)}) = {\sf mins}(w^{(q)}) , \quad {\sf maxs}(w^{(q)}) = {\sf mins}(w^{(r)}), \quad {\sf se}(w^{(q)}) \geq 2, \quad {\sf sw}(w^{(q)}) \geq 3.
\end{equation}
Choose $1 \leq a < b < c < d < e < f < g \leq n$ as follows.
\begin{subequations}
\begin{equation}
a= w^{-1}({\sf mins} (w^{(q)})), \quad b= {\sf mins} (w^{(q)}), \quad c= {\sf mins} (w^{(q)})+1, 
\end{equation}
\begin{equation}
d= {\sf maxs}(w^{(q)}) -  {\sf sw}(w^{(q)}) +1 , \quad e= {\sf maxs}(w^{(q)}) -1 , \quad f= {\sf maxs}(w^{(q)}), 
\end{equation}
\begin{equation}
g = w^{-1}({\sf maxs}(w^{(q)})-{\sf se}(w^{(q)})).
\end{equation}
\end{subequations}
Then we have $w(d) < w(a) < w(e) < w(g) < w(b) < w(c) < w(f)$. This implies a $2561374$ permutation pattern in $w$, leading to a contradiction.

Lastly, the opposite case when there exist $1 \leq p < q < r \leq k$ such that
\begin{equation}
{\sf mins}(w^{(p)}) = {\sf maxs}(w^{(q)}) , \quad {\sf mins}(w^{(q)}) = {\sf maxs}(w^{(r)}), \quad {\sf sw}(w^{(q)}) \geq 2, \quad {\sf se}(w^{(q)}) \geq 3
\end{equation}
can be treated similarly.
\end{proof}

We now consider $244113$ and $233112$ patterns of $[u,I,J]$, which are not implied from Theorem \ref{thm2.1}, \ref{thm2.2} or \ref{thm2.3}.

\begin{theorem}\label{thm6.1}
Let $G= GL_n$, $W=S_n$, $I, J \subseteq [n-1]$ and $u \in {}^I W^J$. If $[u,I,J]$ has a $244113$ or $233112$ pattern, then $[u,I,J]$ is non-spherical.
\end{theorem}
\begin{proof}
If $[u,I,J]$ has a $244113$ pattern, we use $J$-expansion Lemma \ref{lem4.5} to convert the pattern into a $233112$ pattern. Therefore, we may assume that $[u,I,J]$ has a $233112$ pattern.

Write $[u,J] = x_1 , \cdots , x_n$. Recall that $h:=n - |J|$. Since $[u,I,J]$ has a $233112$ pattern, choose $a < b< c <d < e < f$ so that $x_d = x_e < x_a = x_f < x_b = x_c$. Use $J$-expansion Lemma \ref{lem4.5} so that $1 , \cdots , x_d$ and $x_d +1  ,\cdots , x_f$ and $x_f +1 , \cdots , h$ are equalized. Then $[u,I,J']$ still has a $233112$ pattern. Therefore, we may assume that $[n-1]-J = \{ z , z' \}$ with $z < z' $ \textit{i.e.} $x_i \in [3]$.

Due to Lemma \ref{lem2.6}, there exists a $y \in [n-1]-I$ such that $c \leq y < d$. Use $I$-expansion Lemma \ref{lem4.4} on $I' := [n-1]-\{y\}$. Then $[u',I',J]$ still has a $233112$ pattern. Therefore, we may assume that $[n-1]-I = \{y\}$. In other words,
\begin{equation}\label{eqn6.1}
[u,I,J] =  \cdots 2 \cdots 33i11 \cdots 2 \cdots .
\end{equation}
Choose $\lambda \in \Parj$
\begin{equation}\label{eqn6.2}
\lambda := 2\varpi_z + 2\varpi_{z'}.
\end{equation}
From the relation between $y$ and $z<z'$, we have three cases. Due to symmetry, it is enough to cover the cases $y\leq z < z'$ and $z < y <z'$. From a $233112$ pattern, $z'-z \geq 2$.

\begin{figure}
\centering
\begin{subfigure}[t]{0.45\textwidth}
\centering
\includegraphics[width=\textwidth]{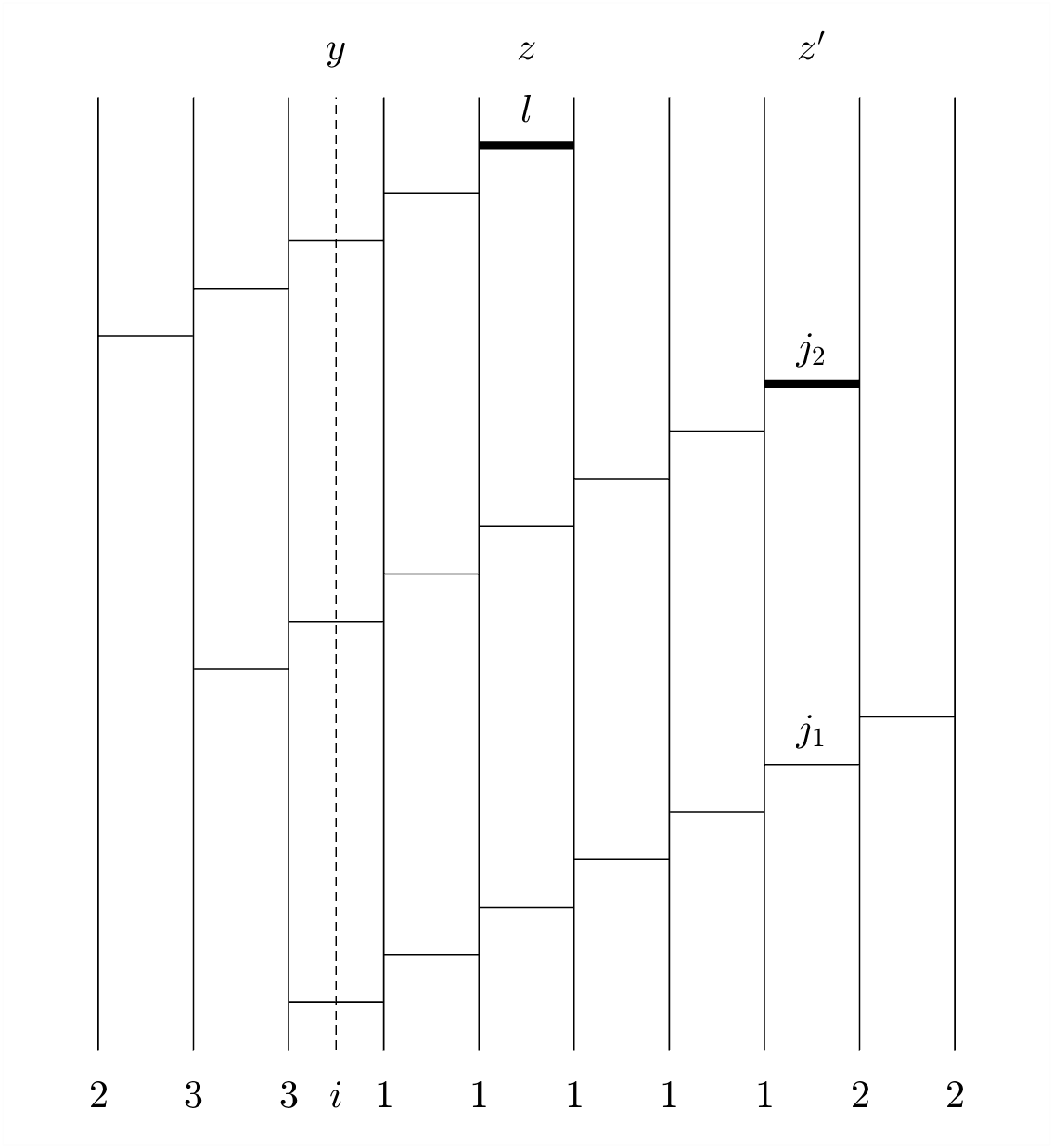}
\caption{(Case 1) $y \leq z < z'$}
\label{fig15}
\end{subfigure}
\hfill
\begin{subfigure}[t]{0.45\textwidth}
\centering
\includegraphics[width=\textwidth]{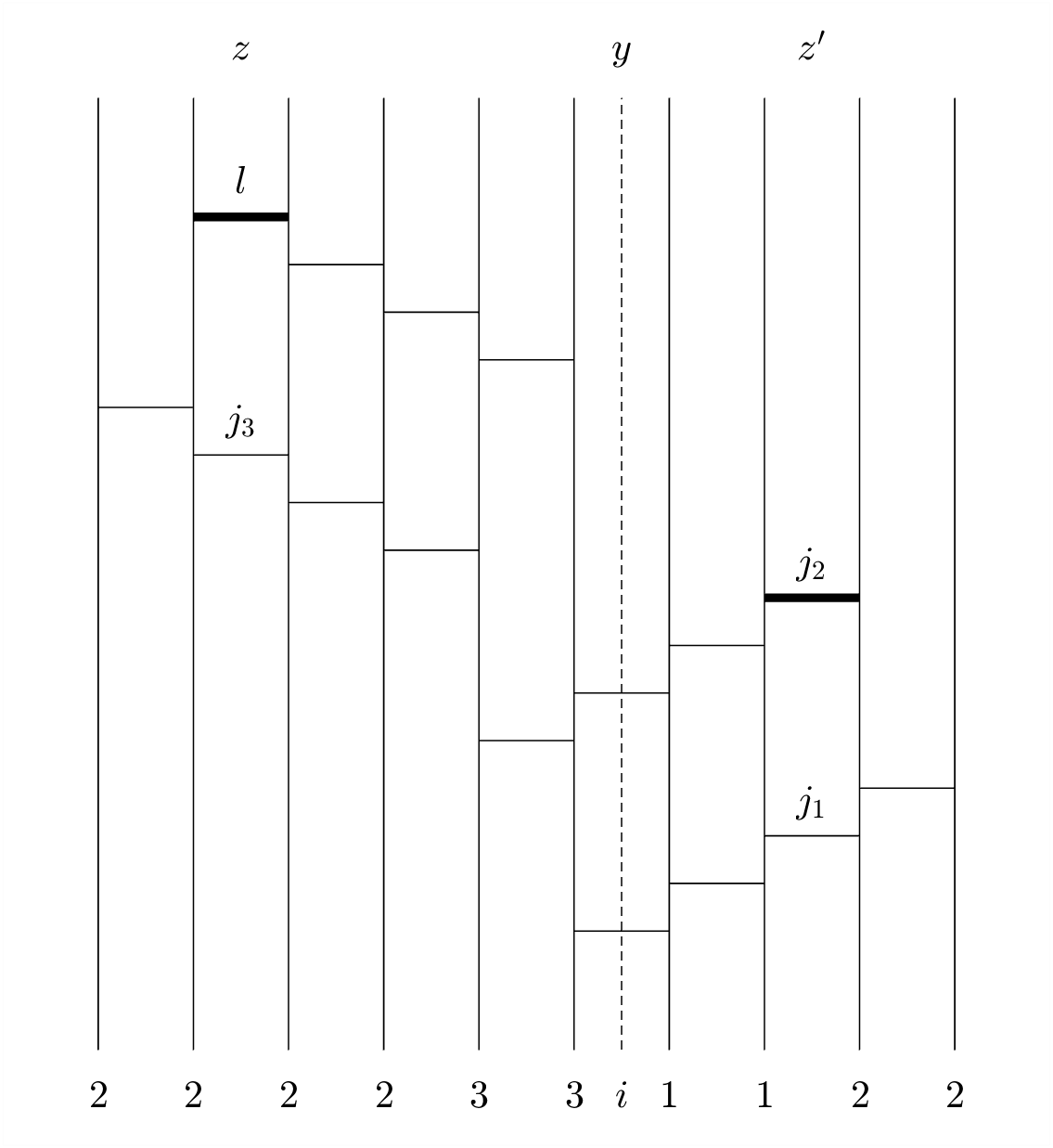}
\caption{(Case 2) $z< y < z'$}
\label{fig16}
\end{subfigure}
\caption{$233112$ patterns in $u'$}
\label{fig01516}
\end{figure}

$\bullet$ (Case 1) $y \leq z< z'$: Due to a $233112$ pattern,
\begin{equation}\label{eqn6.27}
u(z+1) \leq y-2, \quad u(z'+1) \leq y-1, \quad u(z'+2) \leq y.
\end{equation}
As in Figure \ref{fig15}, define
\begin{equation}\label{eqn6.3}
u':=s_y s_{y+1} \cdots s_{z'} s_{z'+1} s_{y-1} s_y \cdots s_{z'-1}s_{z'} s_{y-2} s_{y-1} \cdots s_{z-1}s_z.
\end{equation}
Then $u' \leq u$ by applying Lemma \ref{lem5.5} with $K:=\{z+1, z'+1 , z'+2\}$ and \eqref{eqn6.27}. Also, $u'$ is the minimal representative of $W_Iu'$. Denote $l:= \ell (u')$ and define $(\mathbf{i}, \mathbf{m})$ from the reduced word expression in \eqref{eqn6.3} and $\lambda$ in \eqref{eqn6.2}. We can compute
\begin{equation}\label{eqn6.4}
m_j = 
\begin{cases}
2 & j = 2z'-2y+4 (=:j_2) \\
2 & j = l \\
0 & \text{otherwise}
\end{cases}
.
\end{equation}

Write $j_1:=z'-y+1$. Define $J_{j,k}$ and $K_{j,k}$ for $1\leq j \leq l$ and $1 \leq k \leq m_j$
\begin{equation}\label{eqn6.5}
J_{j,k} = 
\begin{cases}
[l] & j=j_2, \hsp k = 1 \\
[l]-\{j_2\} & j=j_2 , \hsp k=2 \\
[j_1+2, l]-\{j_2\} & j=l, \hsp k=1 \\
[j_2+3,l] & j=l, \hsp k=2
\end{cases}
, \quad K_{j,k}
\begin{cases}
[l] & j =j_2, \hsp k=1 \\
[l] - \{j_2\} & j =j_2, \hsp k=2 \\
[j_1+3,l]-\{j_2\} & j=l, \hsp k=1 \\
[j_2+2,l] & j=l, \hsp k=2
\end{cases}
.
\end{equation}
Then the corresponding LS-monomials $\{\alpha_{j,k}\}$ and $\{\beta_{j,k}\}$ are
\begin{subequations}\label{eqn6.6}
\begin{equation}
\alpha_{j,k} = 
\begin{cases}
\varpi_{z'}-\epsilon_{y-1}-\epsilon_y + \epsilon_{z'+1}+\epsilon_{z'+2}  & j = j_2 , \hsp k=1 \\
\varpi_{z'}-\epsilon_y + \epsilon_{z'+1}  & j = j_2 , \hsp k=2 \\
\varpi_{z}-\epsilon_{y-2}-\epsilon_{y-1} + \epsilon_{z+1}+\epsilon_{z+2}  & j =l , \hsp k=1 \\
\varpi_z - \epsilon_{y} + \epsilon_{z+1} & j=l , \hsp k=2
\end{cases}
,
\end{equation}
\begin{equation}
\beta_{j,k}=
\begin{cases}
\varpi_{z'}-\epsilon_{y-1}-\epsilon_y + \epsilon_{z'+1}+\epsilon_{z'+2}  & j = j_2 , \hsp k=1 \\
\varpi_{z'}-\epsilon_y + \epsilon_{z'+1}  & j = j_2 , \hsp k=2 \\
\varpi_{z}-\epsilon_{y-2}-\epsilon_{y} + \epsilon_{z+1}+\epsilon_{z+2}  & j =l , \hsp k=1 \\
\varpi_z - \epsilon_{y-1} + \epsilon_{z+1} & j=l , \hsp k=2
\end{cases}
.
\end{equation}
\end{subequations}

$\bullet$ (Case 2) $z< y < z'$: Due to a $233112$ pattern,
\begin{equation}\label{eqn6.28}
u(z-1) \geq y+1, \quad u(z) \geq y+2, \quad u(z'+1) \leq y-1, \quad u(z'+2) \leq y.
\end{equation}
As in Figure \ref{fig16}, define
\begin{equation}\label{eqn6.7}
u':=s_y s_{y+1} \cdots s_{z'} s_{z'+1} s_{y-1} s_y \cdots s_{z'-1} s_{z'} s_{y-2} s_{y-3} \cdots s_{z} s_{z-1} s_{y-1} s_{y-2} \cdots s_{z+1} s_z.
\end{equation}
Then $u' \leq u$ by applying Lemma \ref{lem5.5} with $K:=\{z-1 , z , z'+1 , z'+2 \}$ and \eqref{eqn6.28}. Also, $u'$ is the minimal representative of $W_Iu'$. Denote $l:= \ell (u')$ and define $(\mathbf{i}, \mathbf{m})$ from the reduced word expression in \eqref{eqn6.7} and $\lambda$ in \eqref{eqn6.2}. We can compute
\begin{equation}\label{eqn6.8}
m_j = 
\begin{cases}
2 & j = 2z'-2y+4 (=:j_2) \\
2 & j = l \\
0 & \text{otherwise}
\end{cases}
.
\end{equation}

Write $j_1:= z'- y +1$ and $j_3:= 2z'-y-z+3$. Define $J_{j,k}$ and $K_{j,k}$ for $1\leq j \leq l$ and $1 \leq k \leq m_j$
\begin{subequations}\label{eqn6.9}
\begin{equation}
J_{j,k} = 
\begin{cases}
[l] & j=j_2, \hsp k=1 \\
[l]-\{j_2\} & j=j_2 , \hsp k=2 \\
[j_1+2,l] -\{j_2\} & j=l, \hsp k=1 \\
[j_3+3,l] & j=l, \hsp k=2
\end{cases}
,
\end{equation}
\begin{equation}
K_{j,k}=
\begin{cases}
[l] & j=j_2, \hsp k=1 \\
[l]-\{j_2\} & j=j_2 , \hsp k=2 \\
[j_1+3,l]-\{j_2\} & j=l, \hsp k=1 \\
[j_3+2, l] & j=l, \hsp k=2
\end{cases}
.
\end{equation}
\end{subequations}
Then the corresponding LS-monomials $\{\alpha_{j,k}\}$ and $\{\beta_{j,k}\}$ are
\begin{subequations}\label{eqn6.10}
\begin{equation}
\alpha_{j,k} = 
\begin{cases}
\varpi_{z'}-\epsilon_{y-1}-\epsilon_y + \epsilon_{z'+1}+\epsilon_{z'+2}  & j = j_2 , \hsp k=1 \\
\varpi_{z'}-\epsilon_y + \epsilon_{z'+1}  & j = j_2 , \hsp k=2 \\
\varpi_{z}-\epsilon_{z-1}-\epsilon_{z} + \epsilon_{y}+\epsilon_{y+1}  & j =l , \hsp k=1 \\
\varpi_z - \epsilon_{z} + \epsilon_{y-1} & j=l , \hsp k=2
\end{cases}
,
\end{equation}
\begin{equation}
\beta_{j,k}=
\begin{cases}
\varpi_{z'}-\epsilon_{y-1}-\epsilon_y + \epsilon_{z'+1}+\epsilon_{z'+2}  & j = j_2 , \hsp k=1 \\
\varpi_{z'}-\epsilon_y + \epsilon_{z'+1}  & j = j_2 , \hsp k=2 \\
\varpi_{z}-\epsilon_{z-1}-\epsilon_{z} + \epsilon_{y-1}+\epsilon_{y+1}  & j =l , \hsp k=1 \\
\varpi_z - \epsilon_{z} + \epsilon_{y} & j=l , \hsp k=2
\end{cases}
.
\end{equation}
\end{subequations}
From (Case 1) to (Case 2), we constructed two distinct $I$-dominant LS-monomials $\{\alpha_{j,k}\}$ and $\{\beta_{j,k}\}$ with the same weight based on $u'$ and $\lambda$. Let $w'$ be the maximal representative of $W_I u'$. Due to Theorem \ref{thm4.1}, $H^0(X(w'W_J), L_\lambda)^*$ is not $L_I$-multiplicity-free. From $u'\leq u$, we have $w' \leq w$. Using Lemma \ref{lem2.5}, $H^0(X(wW_J), L_\lambda)^*$ also is not $L_I$-multiplicity-free. Due to Theorem \ref{lem2.4}, $X(wW_J)$ is not $L_I$-spherical.
\end{proof}

\begin{theorem}\label{thm6.5}
Let $G= GL_n$, $W=S_n$, $I, J \subseteq [n-1]$ and $u \in {}^I W^J$. If $[u,I,J]$ has $221211$, $221311$, $331322$ or $331422$ patterns, then $[u,I,J]$ is non-spherical.
\end{theorem}
\begin{proof}
If $[u,I,J]$ has $221311$, $331322$ or $331422$ patterns, we use $J$-expansion Lemma \ref{lem4.5} to convert the pattern into a $221211$ pattern. Therefore, we may assume that $[u,I,J]$ has a $221211$ pattern.

Write $[u,J] = x_1 , \cdots , x_n$. Recall that $h:=n - |J|$. Since $[u,I,J]$ has a $221211$ pattern, choose $a < b< c <d < e < f$ so that $x_c = x_e = x_f < x_a = x_b = x_d$. Use $J$-expansion Lemma \ref{lem4.5} so that $1 , \cdots , x_f$ and $x_f +1  ,\cdots , h$ are equalized. Then $[u,I,J']$ still has a $221211$ pattern. Therefore, we may assume that $[n-1]-J = \{ z \}$.

Due to Lemma \ref{lem2.6}, there exist $y , y' \in [n-1]-I$ such that $b \leq y < c$ and $d \leq y' < e$. Use $I$-expansion Lemma \ref{lem4.4} on $I' := [n-1]-\{y , y'\}$. Then $[u',I',J]$ still has a $221211$ pattern. Therefore, we may assume that $[n-1]-I = \{y , y'\}$. In other words,
\begin{equation}\label{eqn6.29}
[u,I,J] =  \cdots 22i1 \cdots 2i11 \cdots .
\end{equation}
Choose $\lambda \in \Parj$
\begin{equation}\label{eqn6.30}
\lambda := 3\varpi_z.
\end{equation}
From the relation between $y < y'$ and $z$, we have three cases. Due to symmetry, it is enough to cover the cases $z \leq y < y'$ and $y < z <y'$.

\begin{figure}
\centering
\begin{subfigure}[t]{0.45\textwidth}
\centering
\includegraphics[width=\textwidth]{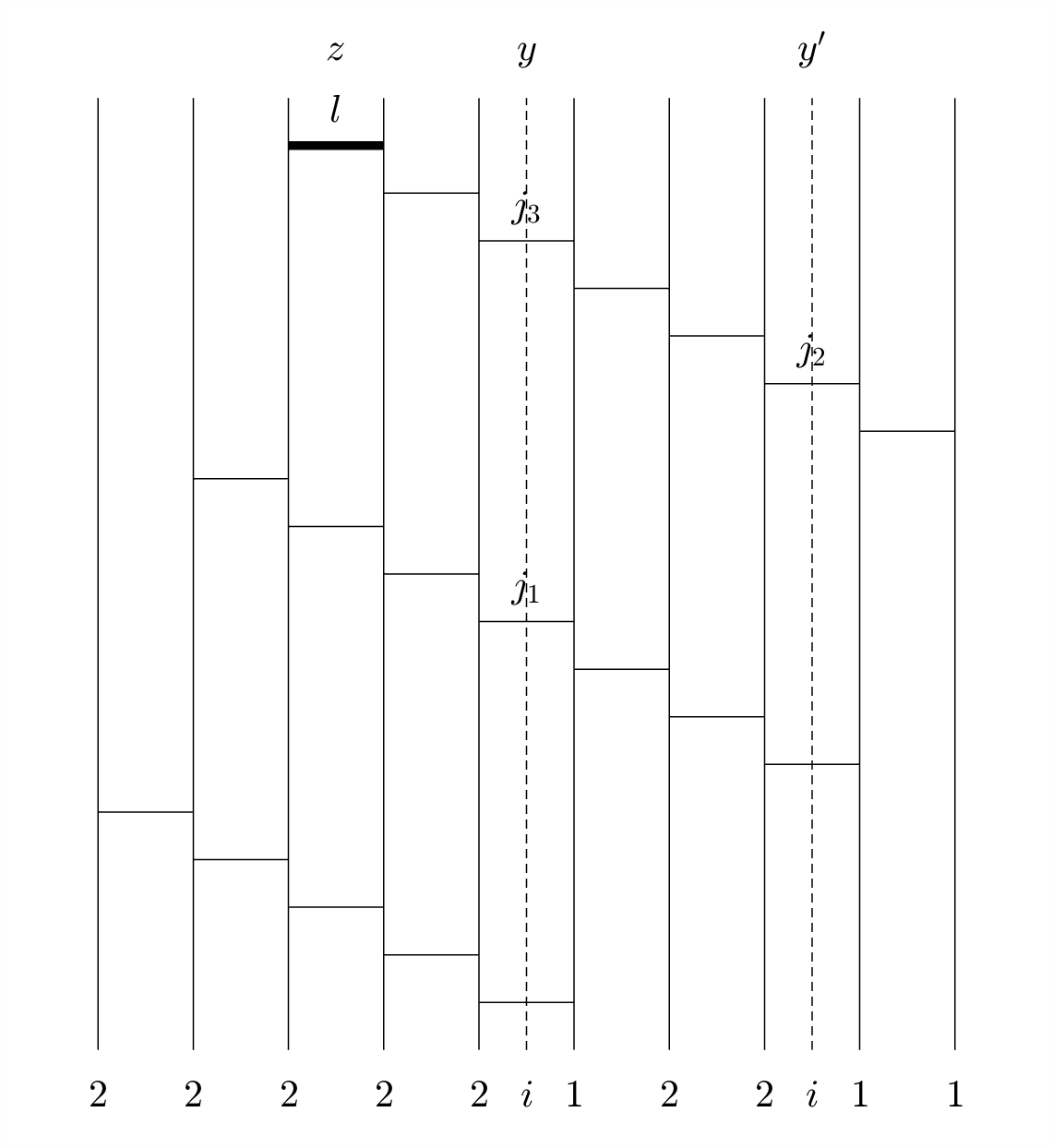}
\caption{(Case 1) $z \leq y < y'$}
\label{fig25}
\end{subfigure}
\hfill
\begin{subfigure}[t]{0.45\textwidth}
\centering
\includegraphics[width=\textwidth]{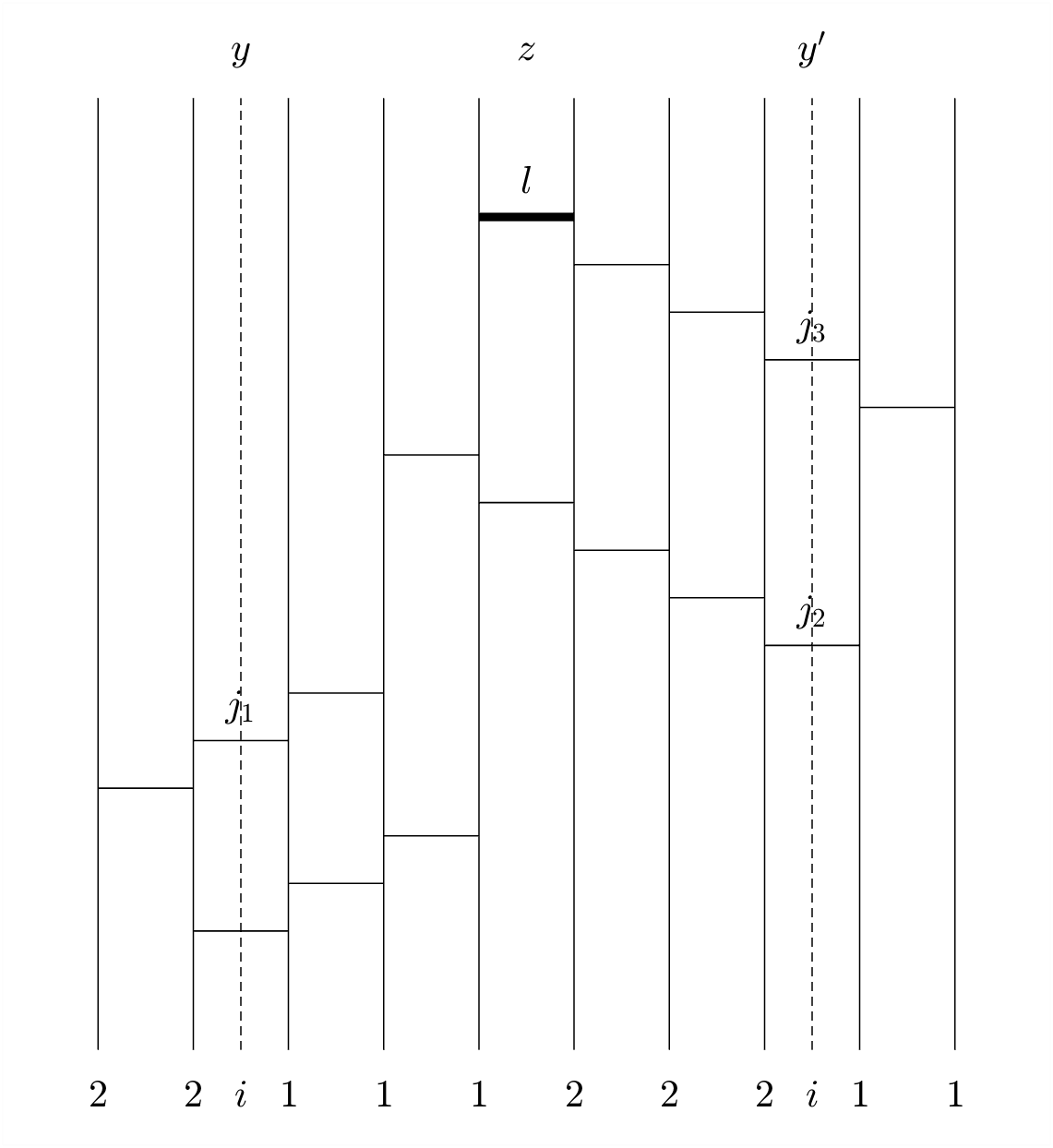}
\caption{(Case 2) $y < z <y'$}
\label{fig26}
\end{subfigure}
\caption{$221211$ patterns in $u'$}
\label{fig02526}
\end{figure}

$\bullet$ (Case 1) $z \leq y< y'$: Due to a $221211$ pattern,
\begin{equation}\label{eqn6.31}
u(z-2) \geq y+1, \quad u(z-1) \geq y'+1, \quad u(z) \geq y'+2.
\end{equation}
As in Figure \ref{fig25}, define
\begin{equation}\label{eqn6.32}
u':=s_{y} s_{y-1} \cdots s_{z-1}s_{z-2} s_{y'} s_{y'-1} \cdots s_z s_{z-1} s_{y'+1} s_{y'} \cdots s_{z+1 }s_z.
\end{equation}
Then $u' \leq u$ by applying Lemma \ref{lem5.5} with $K:=\{z-2 , z-1 , z \}$ and \eqref{eqn6.31}. Also, $u'$ is the minimal representative of $W_Iu'$. Denote $l:= \ell (u')$ and define $(\mathbf{i}, \mathbf{m})$ from the reduced word expression in \eqref{eqn6.32} and $\lambda$ in \eqref{eqn6.30}. We can compute
\begin{equation}\label{eqn6.33}
m_j = 
\begin{cases}
3 & j = l \\
0 & \text{otherwise}
\end{cases}
.
\end{equation}

Write $j_1:= y' - z +4$ , $j_2 := y' + y - 2z +7$ and $j_3 := 2y' - 2z+7$. Define $J_{j,k}$ and $K_{j,k}$ for $1\leq j \leq l$ and $1 \leq k \leq m_j$
\begin{equation}\label{eqn6.34}
J_{j,k} = 
\begin{cases}
[l] & j = l , \hsp k=1 \\
[j_1-1 , j_2-2] \cup [j_2,l] & j=l, \hsp k=2 \\
[j_3 , l] & j=l, \hsp k=3
\end{cases}
, \quad K_{j,k}=
\begin{cases}
[l] & j = l , \hsp k=1 \\
[j_1 , j_2-2] \cup [j_2,l] & j=l, \hsp k=2 \\
[j_3-1 , l] & j=l, \hsp k=3
\end{cases}
.
\end{equation}
Then the corresponding LS-monomials $\{\alpha_{j,k}\}$ and $\{\beta_{j,k}\}$ are
\begin{subequations}\label{eqn6.35}
\begin{equation}
\alpha_{j,k} = 
\begin{cases}
\varpi_z - \epsilon_{z-2} - \epsilon_{z-1} - \epsilon_z + \epsilon_{y+1} + \epsilon_{y' + 1} + \epsilon_{y'+2} & j=l , \hsp k=1 \\
\varpi_z - \epsilon_{z-1} - \epsilon_z + \epsilon_{y+2} + \epsilon_{y' + 1} & j=l , \hsp k=2 \\
\varpi_z - \epsilon_z + \epsilon_{y+1}  & j=l , \hsp k=3
\end{cases}
,
\end{equation}
\begin{equation}
\beta_{j,k}=
\begin{cases}
\varpi_z - \epsilon_{z-2} - \epsilon_{z-1} - \epsilon_z + \epsilon_{y+1} + \epsilon_{y' + 1} + \epsilon_{y'+2} & j=l , \hsp k=1 \\
\varpi_z - \epsilon_{z-1} - \epsilon_z + \epsilon_{y+1} + \epsilon_{y' + 1} & j=l , \hsp k=2 \\
\varpi_z - \epsilon_z + \epsilon_{y+2}  & j=l , \hsp k=3
\end{cases}
.
\end{equation}
\end{subequations}

$\bullet$ (Case 2) $y< z < y'$: Due to a $221211$ pattern,
\begin{equation}\label{eqn6.36}
u(z-1) \geq y'+1 , \quad u(z) \geq y'+2 , \quad u(z+1) \leq y-1, \quad u(z+2) \leq y.
\end{equation}
As in Figure \ref{fig26}, define
\begin{equation}\label{eqn6.37}
u':= s_y s_{y+1} \cdots s_{z-2} s_{z-1} s_{y-1} s_y \cdots s_{z-3} s_{z-2} s_{y'} s_{y'-1} \cdots s_{z} s_{z-1} s_{y'+1} s_{y'} \cdots s_{z+1} s_z.
\end{equation}
Then $u' \leq u$ by applying Lemma \ref{lem5.5} with $K:=\{ z-1 , z , z+1 , z+2 \}$ and \eqref{eqn6.36}. Also, $u'$ is the minimal representative of $W_Iu'$. Denote $l:= \ell (u')$ and define $(\mathbf{i}, \mathbf{m})$ from the reduced word expression in \eqref{eqn6.37} and $\lambda$ in \eqref{eqn6.30}. We can compute
\begin{equation}\label{eqn6.38}
m_j = 
\begin{cases}
3 & j = l \\
0 & \text{otherwise}
\end{cases}
.
\end{equation}

Write $j_1:= z-y +2$, $j_2 := 2z-2y+1$ and $j_3:= z-2y+y'+4$. Define $J_{j,k}$ and $K_{j,k}$ for $1\leq j \leq l$ and $1 \leq k \leq m_j$
\begin{subequations}\label{eqn6.39}
\begin{equation}
J_{j,k} = 
\begin{cases}
[l] & j=l , \hsp k=1 \\
[j_1 , j_2-1]\cup \{ j_3-3 , j_3 - 2 \} \cup [j_3 , l] & j=l , \hsp k=2 \\
\emptyset & j = l , \hsp k=3
\end{cases}
,
\end{equation}
\begin{equation}
K_{j,k}=
\begin{cases}
[l] & j=l , \hsp k=1 \\
[j_1 , j_2-1]\cup \{ j_3 - 2 \} \cup [j_3 , l] & j=l , \hsp k=2 \\
\{l\} & j = l , \hsp k=3
\end{cases}
.
\end{equation}
\end{subequations}
Then the corresponding LS-monomials $\{\alpha_{j,k}\}$ and $\{\beta_{j,k}\}$ are
\begin{subequations}\label{eqn6.40}
\begin{equation}
\alpha_{j,k} = 
\begin{cases}
\varpi_z -\epsilon_{y-1} - \epsilon_y + \epsilon_{y'+1} + \epsilon_{y'+2} & j=l , \hsp k=1 \\
\varpi_z - \epsilon_y -\epsilon_z + \epsilon_{z+1}+ \epsilon_{y'+1}  & j=l , \hsp k=2 \\
\varpi_z & j=l , \hsp k=3 
\end{cases}
,
\end{equation}
\begin{equation}
\beta_{j,k}=
\begin{cases}
\varpi_z -\epsilon_{y-1} - \epsilon_y + \epsilon_{y'+1} + \epsilon_{y'+2} & j=l , \hsp k=1 \\
\varpi_z - \epsilon_y + \epsilon_{y'+1}  & j=l , \hsp k=2 \\
\varpi_z -\epsilon_z + \epsilon_{z+1} & j=l , \hsp k=3 \\
\end{cases}
.
\end{equation}
\end{subequations}
From (Case 1) to (Case 2), we constructed two distinct $I$-dominant LS-monomials $\{\alpha_{j,k}\}$ and $\{\beta_{j,k}\}$ with the same weight based on $u'$ and $\lambda$. Let $w'$ be the maximal representative of $W_I u'$. Due to Theorem \ref{thm4.1}, $H^0(X(w'W_J), L_\lambda)^*$ is not $L_I$-multiplicity-free. From $u'\leq u$, we have $w' \leq w$. Using Lemma \ref{lem2.5}, $H^0(X(wW_J), L_\lambda)^*$ is also not $L_I$-multiplicity-free. Due to Theorem \ref{lem2.4}, $X(wW_J)$ is not $L_I$-spherical.
\end{proof}

\begin{lemma}\label{lem7.1}
Let $G = GL_n$, $W= S_n$, $I,J \subseteq [n-1]$ and $u \in {}^I W^J$. Suppose $[u,I,J]$ avoids
\begin{equation}\label{eqn7.3}
\begin{matrix}
321, & 3412, & 3i312, & 231i1, & 2i21i1, & \\
244113, & 233112, & 221211, & 221311, & 331322, & 331422
\end{matrix}
\end{equation}
patterns. Then $u$ avoids permutation patterns as follows.
\begin{equation}\label{eqn7.4}
321, \quad 246135, \quad 351624, \quad 356124, \quad 415263, \quad 451623.
\end{equation}
\end{lemma}
\begin{proof}
Due to Lemma \ref{lem6.0}, it is sufficient to prove that $u$ does not have $356124$ and $451623$ patterns. Since $356124$ and $451623$ are inverses of each other, we need to prove that $u^{-1}$ avoids these patterns.

Suppose $u^{-1}$ has a $356124$ pattern. Recall that $[u,J]$ is obtained by equalizing consecutive numbers from the left to the right in $u^{-1}$. Since $\underline{35}612\underline{4}$ cannot be equalized and $[u,I,J]$ does not have $3412$, $1$ and $2$ are equalized. So are $5$ and $6$. Then $[u,I,J]$ has a $244113$ or $233112$ pattern, leading to a contradiction.

Suppose $u^{-1}$ has a $451623$ pattern. Since there are descents between $5$, $1$ and $6$, $2$, $i$ should be inserted between them in $[u,I,J]$. From $3412$ in $\underline{45}16\underline{23}$, one of them should be equalized. If both of them are equalized, we have $221211$, $221311$, $331322$ or $331422$ patterns in $[u,I,J]$. If $2$ and $3$ in $4516\underline{23}$ are not equalized, then $4,5,6$ in $\underline{45}1\underline{6}23$ should be equalized due to $3412$. This gives a $3i312$ pattern, leading to a contradiction. The case when $\underline{45}1623$ is not equalized can be treated similarly.
\end{proof}

\subsection{Reduction to bigrassmannians}\label{sub6.2} 

We prove a technical lemma using LS-monomials. Let $G = GL_n$ and $W=S_n$. Let $I,J \subseteq [n-1]$ and $u \in W_I w W_J$ be the minimal representative with a reduced word expression $\mathbf{i} = (i_1 , \cdots , i_l)$ of $u$. Assume that there exists a $p \in [n-2]$
\begin{equation}\label{eqn6.11}
\{ t \in [l] \mid i_t = p\} = \{r\} , \quad \{t \in [l] \mid i_t = p+1 \} =\{s\}.
\end{equation}

From $\mathbf{i} = (i_1 , \cdots , i_l)$, take entries which are less than $p+1$ and define a subsequence $\mathbf{i}' := (i_1 ' , \cdots , i_{l'}')$. Also, take entries which are more than $p$ and define a subsequence $\mathbf{i}'' := (i_1 '' , \cdots , i_{l''}'')$. In other words, $i_j' < p+1$, $i_j '' > p$ and $l' + l'' = l$. Define
\begin{equation}\label{eqn6.23}
u' := s_{i_1'} s_{i_2'} \cdots s_{i_{l'}'},  \quad u'' := s_{i_1''} s_{i_2''} \cdots s_{i_{l''}''}.
\end{equation}

Notice that the constructions of $u'$ and $u''$ do not depend on the choice of $\mathbf{i}$. To show this, suppose $\mathbf{\hat{i}} = (\hat{i}_1 , \cdots \hat{i}_l)$ be another reduced word of $u$. Suppose $r <s$. Then
\begin{equation}\label{eqn5.71}
\hat{r} \in \{ t \in [l] \mid \hat{i}_t = p \} , \hsp \hat{s} \in \{ t \in [l] \mid \hat{i}_t = p+1 \} \quad \Rightarrow \quad \hat{r} < \hat{s}.
\end{equation}
Indeed, $\mathbf{i}$ and $\mathbf{\hat{i}}$ are identified using
\begin{equation}
(\cdots , i , i+1 , i , \cdots ) \sim ( \cdots , i+1 , i , i+1 , \cdots ) , \hsp (\cdots , i , j , \cdots ) \sim ( \cdots , j , i , \cdots ) , \hsp |i-j|>1.
\end{equation}
Therefore, the property \eqref{eqn5.71} is preserved under these relations since an identification $( \cdots , p , p+1 , p , \cdots) \sim ( \cdots , p+1 , p, p+1 , \cdots)$ never occurs. On the other hand, since identifications $(\cdots , i , i+1 , i , \cdots ) \sim (\cdots , i+1 , i , i+1 , \cdots )$ only happen for $i \neq p$, constructions of $u'$ and $u''$ are all the same. $r>s$ case can be treated similarly.

\begin{figure}
\centering
\begin{subfigure}[t]{0.45\textwidth}
\centering
\includegraphics[width=\textwidth]{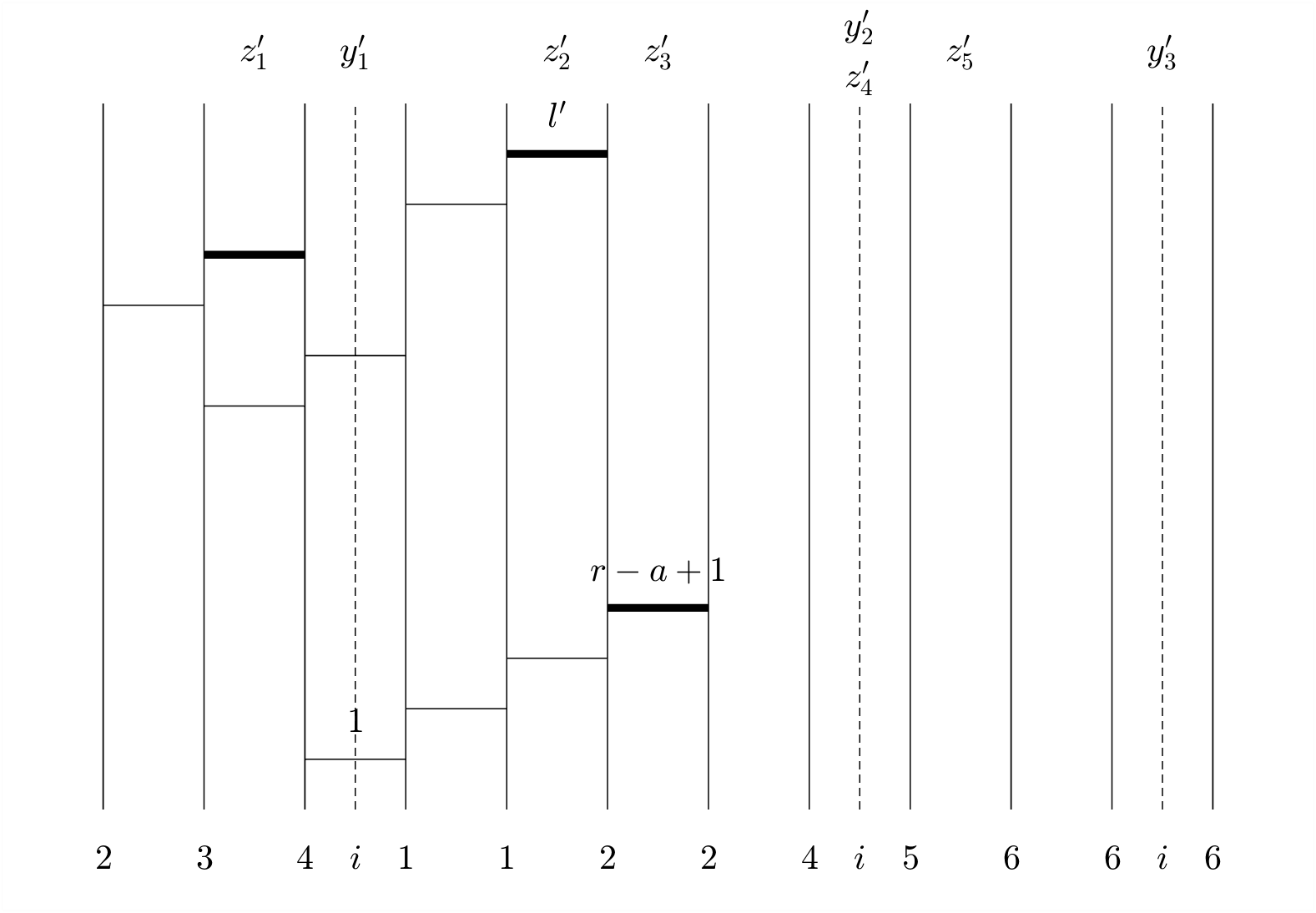}
\caption{$u' , I' , J'$}
\label{fig18}
\end{subfigure}
\hfill
\begin{subfigure}[t]{0.45\textwidth}
\centering
\includegraphics[width=\textwidth]{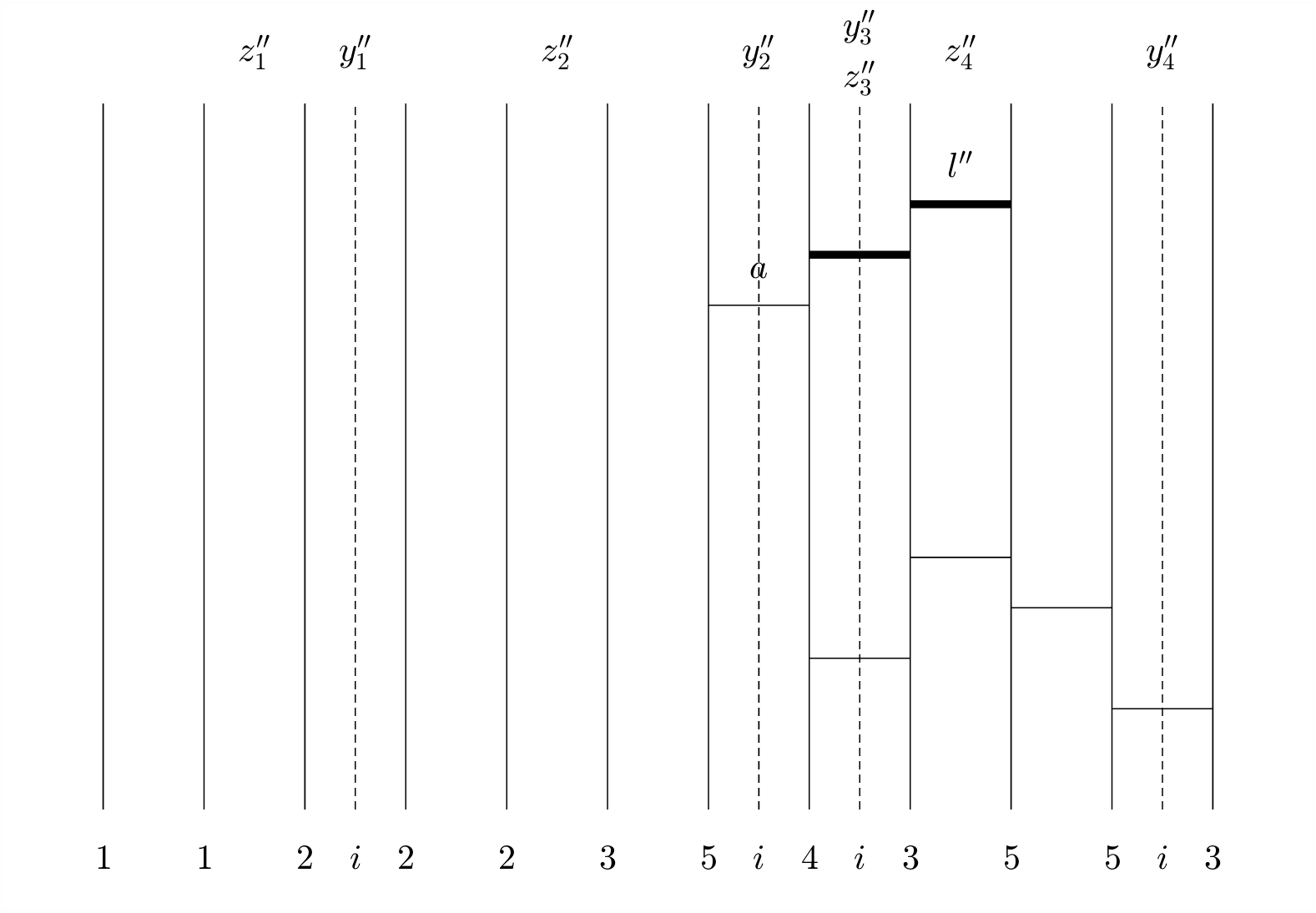}
\caption{$u'', I'', J''$}
\label{fig19}
\end{subfigure}
\begin{subfigure}[t]{0.45\textwidth}
\centering
\includegraphics[width=\textwidth]{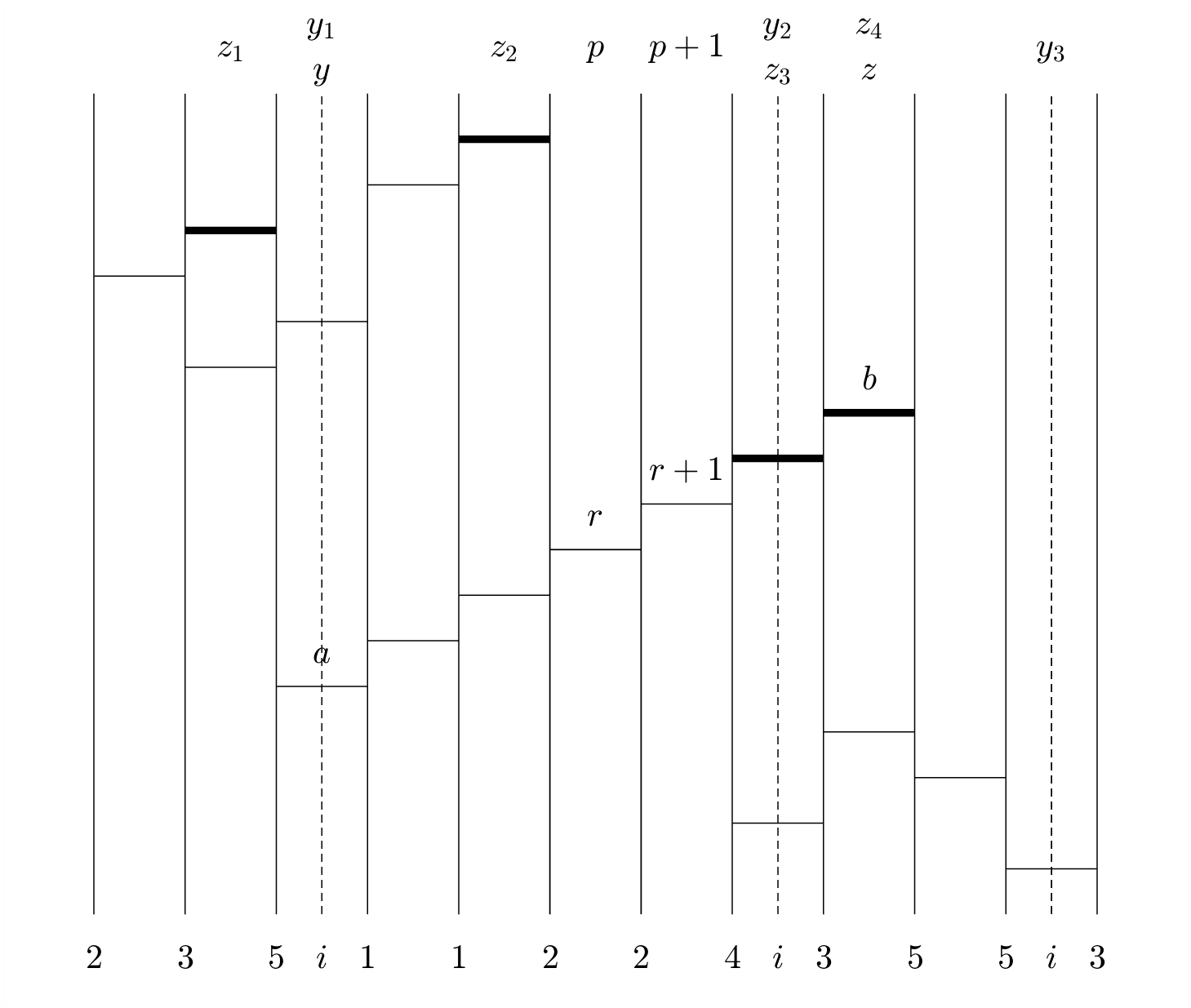}
\caption{$u,I,J$}
\label{fig17}
\end{subfigure}
\caption{$p=6$, $r=8$, $a=5$ and $b=11$}
\label{fig0171819}
\end{figure}

If $r<s$, define
\begin{equation}\label{eqn6.12}
I' := I, \quad  J' := J-\{p\}, \quad  I'' := I-\{p+1\}, \quad  J'' := J.
\end{equation}
If $r>s$, define
\begin{equation}\label{eqn6.13}
I' := I-\{p\}, \quad  J' := J, \quad I'' := I, \quad  J'' := J-\{p+1\}.
\end{equation}
Then $u'$ and $u''$ are minimal representatives of $W_{I'}u' W_{J'}$ and $W_{I''}u''W_{J''}$, respectively. For an example, see Figure \ref{fig0171819}.

\begin{lemma}\label{lem6.2}
Let $u,I,J$ be given and $u',I',J'$ and $u'', I'', J''$ be defined as above.
\begin{itemize}
\item[(a)] $[u,I,J]$ avoids the patterns in \eqref{eqn1.7} if and only if $[u',I',J']$ and $[u'',I'',J'']$ avoid the patterns in \eqref{eqn1.7}.
\item[(b)] $[u,I,J]$ is spherical if and only if $[u',I',J']$ and $[u'',I'',J'']$ are spherical.
\end{itemize}
\end{lemma}
\begin{proof}
Due to symmetry, it is sufficient to prove the case when $r<s$. Note that the constructions of $u'$ and $u''$ do not depend on reduced word expressions of $u$. Without loss of generality, we may choose a reduced expression $\mathbf{i} = (i_1 , \cdots , i_l)$ of $u$ so that $r+1 = s$ and there exist $a \leq r < b$ satisfying the following.
\begin{itemize}
\item $1 \leq j \leq a-1$: $i_j >p+1$.
\item $ a\leq j \leq b$: $i_j = j-r+p$. In particular, $i_r = p$ and $i_{r+1} = p+1$.
\item $b+1 \leq j \leq l$: $i_j <p$.  
\end{itemize}

Indeed, this is achieved by applying a relation $s_{k+1} s_1 s_2 \cdots s_n = s_1 s_2 \cdots s_n s_{k}$ repeatedly. To elaborate, suppose there exists $j$ such that $i_j < p$ but $j <a$. Assume that $j$ is the maximal index with that property. Using $(\cdots , i_j , i_{j+1} , \cdots) \sim (\cdots , i_{j+1} , i_j, \cdots)$ with $|i_j - i_{j+1}|>1$, we have an identified word with $j = a-1$.
\begin{enumerate}
\item If $i_j \in [a-r+p+1 ,  p-1]$, then use
\begin{equation*}
(\cdots , i_j , a-r+p , a-r+p+1 , \cdots , b-r+p, \cdots ) \sim (\cdots , a-r+p , a-r+p+1 , \cdots ,  b-r+p, i_j - 1\cdots ).
\end{equation*}
\item $i_j = a-r+p$ cannot occur since $\mathbf{i}$ is a reduced word.
\item If $i_j = a-r+p-1$, add $i_j$ into the string
\begin{equation*}
(\cdots , a-r+p-1, a-r+p , a-r+p+1 , \cdots ,b-r+p, \cdots).
\end{equation*}
\item If $i_j \in [1 , a-r+p-2]$, then use
\begin{equation*}
(\cdots , i_j , a-r+p , a-r+p+1 , \cdots , b-r+p, \cdots ) \sim (\cdots , a-r+p , a-r+p+1 , \cdots , b-r+p, i_j \cdots ).
\end{equation*}
\end{enumerate}
The opposite case when there exists $j$ such that $i_j > p+1$ but $j >b$ can be treated similarly. Therefore, we may assume that $\mathbf{i}$ has the properties as desired. For an example, see Figure \ref{fig17}. Our first claim is that
\begin{equation}\label{eqn6.14}
y:=i_a \in [n-1]-I , \quad z:=i_b \in [n-1]-J.
\end{equation}
To prove this, write $[u,J] = x_1 \cdots x_n$. As a permutation, $u^{-1}(y)=z+1$ and $u^{-1}(y+1) \leq p$. Therefore, $y \in [n-1]-I$. Similarly, if $z \in J$, then this contradicts the condition that $u$ is the minimal representative of $uW_J$ since $us_z < u$. This proves our claim.

Let $l':=\ell (u')$ and $l'':=\ell (u'')$. Then
\begin{equation}
l'=l-b+r-a+1, \quad l''=b-r+a-1.
\end{equation}
Using the construction in \eqref{eqn6.23}, we have reduced word expressions $\mathbf{i'} = (i_1' , \cdots, i_{l'}')$ and $\mathbf{i''} = (i_1'' , \cdots, i_{l''}'')$ for $u'$ and $u''$, respectively. In particular,
\begin{equation}\label{eqn6.24}
i_j ' =
\begin{cases}
i_{j+a-1} & 1 \leq j \leq r-a+1 \\
i_{j+a+b-r-1} & r-a+2 \leq j \leq l'
\end{cases}
, \quad i_j''=
\begin{cases}
i_j & 1 \leq j \leq a-1 \\
i_{j+r-a+1} & a \leq j \leq l''
\end{cases}
.
\end{equation}
Write $[u' , J'] = x_1 '  \cdots x_n' $ and $[u'', J''] = x_1 '' \cdots x_n ''$.

``$\Rightarrow$ of (a)'':  We prove the contraposition. Suppose $[u',I',J']$ has a pattern at $x_{p_1}' x_{p_2}' \cdots x_{p_k}'$ with some $i$ insertions. If $p_k \geq p+2$, then $x_{p_k}' \geq x_{p_j}'$ for $1 \leq j \leq k-1$. However, all the patterns in \eqref{eqn1.7} do not have this property, leading to $p_k \leq p+1$. Then $x_{p_1} x_{p_2} \cdots x_{p_k}$ with some $i$ insertions in $[u,I,J]$ also has the same pattern.

On the other hand, suppose $[u'',I'',J'']$ has a pattern at $x_{p_1}'' x_{p_2}'' \cdots x_{p_k}''$ with some $i$ insertions. If $p_1 \leq p$, then $x_{p_1}'' \leq x_{p_j}''$ for $2 \leq j \leq k$. However, all the patterns in \eqref{eqn1.7} do not have this property, leading to $p_1 \geq p+1$. If $p_1 \neq p+1$, then $x_{p_1} x_{p_2} \cdots x_{p_k}$ with some $i$ insertions in $[u,I,J]$ also has the same pattern. If $p_1 = p+1$, then $[u,I,J]$ has the same pattern at $x_y x_{p_2}x_{p_3} \cdots x_{p_k}$ with $i$ insertions.

``$\Leftarrow$ of (a)'': We prove the contraposition. Suppose $[u,I,J]$ has a pattern at $x_{p_1} x_{p_2} \cdots x_{p_k}$ with $i$ insertions. If $p_k \leq p+1$, then $[u',I',J']$ has the same pattern at $x_{p_1}' x_{p_2}' \cdots x_{p_k}'$ with $i$ insertions. If $p_1 >p+1$, then $[u'',I'',J'']$ has the same pattern at $x_{p_1}''x_{p_2}'' \cdots x_{p_k}''$ with $i$ insertions. If $p_1 = y$ and $p_2 > p+1$, then $[u'', I'', J'']$ has the same pattern at $x_{1+p}'' x_{p_2}'' x_{p_3}'' \cdots x_{p_k}''$ with $i$ insertions.

Otherwise, $p_1 \leq p+1$ and $p_k >p+1$. Suppose the pattern is $244113$ or $233112$. Either way, $x_{p_2} = x_{p_3} > x_{p_6}$ and $x_{p_1} > x_{p_5}$. From $x_{p_2} > x_{p_6}$, $p_2 = y$ or $p_2 > p+1$. If $p_2 = y$, then $p_3 \neq y$ with $x_{p_3} > x_{p_6}$. Then $p_3 >p+1$, leading to $p_5 > p+1$. This is a contradiction with $x_{p_1} > x_{p_5}$. Therefore, $p_2 > p+1$. From $x_{p_1} > x_{p_5}$, $p_1 = y$, which was covered previously.

Lastly, suppose $p_1 \leq p+1$ and $p_k >p+1$ and the pattern is one of \eqref{eqn1.7} except for $244113$ and $233112$. Then $x_{p_1}>x_{p_k}$ and $x_{p_2} > x_{p_k}$. Therefore, $p_1 = y$ and $p_2 >p+1$. Again, this was covered previously.

``$\Rightarrow$ of (b)'' : We prove the contraposition. There are two cases as follows.

$\bullet$ (Case 1) $[u',I',J']$ is non-spherical: Then there exist a $\lambda' \in {\rm Par}_{J'}$, a corresponding $\mathbf{m'} = (m_1 ' , \cdots , m_{l'}')$ and LS-monomials $\{\alpha_{j,k}'\}$ and $\{\beta_{j,k}' \}$ which are $I'$-dominant with the same weight. Let $\{J_{j,k}'\}$ and $\{K_{j,k}'\}$ be the LS-chains of $\{\alpha_{j,k}'\}$ and $\{\beta_{j,k}'\}$, respectively. For $1 \leq j \leq l$, define $\mathbf{m} = (m_1 , \cdots , m_l)$
\begin{equation}\label{eqn6.15}
m_j := 
\begin{cases}
0 & 1 \leq j \leq a-1 \\
m_{j-a+1}' & a \leq j \leq r-1 \\
0 & r \leq j \leq b-1 \\
m_{j+r-a-b +1}' & j\geq b
\end{cases}
.
\end{equation}

Choose $\lambda \in \Parj$ matching $\mathbf{m}$ and $\mathbf{i}$. Define a function $\varphi_1: [l] \rightarrow [l'] \cup \{0\}$ between index sets
\begin{equation}
\varphi_1 (j):=
\begin{cases}
0 & 1 \leq j \leq a-1 \\
j-a+1 & a \leq j \leq r \\
r-a+1 & r+1 \leq j \leq b \\
j-b+r-a+1 & j \geq b+1
\end{cases}
.
\end{equation}

For each $1 \leq j \leq l$ and $1 \leq k \leq m_j$, define
\begin{equation}
J_{j,k}:=
\begin{cases}
\varphi_1^{-1}(J_{j-a+1,k}') & a \leq j \leq r-1 \\
\varphi_1^{-1}(J_{j+r-a-b+1,k}') & j \geq b
\end{cases}
.
\end{equation}

Then we can compute the corresponding LS-monomials $\{ \alpha_{j,k} \}$ from $\{J_{j,k}\}$ and
\begin{equation}\label{eqn6.18}
\alpha_{j,k} =
\begin{cases}
\alpha_{j-a+1,k}' & a \leq j \leq r-1 \\
s_{p+1,z+1}(\alpha_{r-a+1,k}') + \sum_{t=p+1}^z \epsilon_t & j = b \\
\alpha_{j-a-b+1+r,k}'  & b+1 \leq j \leq l
\end{cases}
.
\end{equation}
Here, $s_{p+1,z+1}$ is a transposition switching $(p+1)$th and $(z+1)$th components. Therefore, we have an $I$-dominant LS-monomial $\{\alpha_{j,k}\}$. Similarly, we construct $\{\beta_{j,k}\}$ from $\{\beta_{j,k}'\}$.

$\bullet$ (Case 2) $[u'',I'',J'']$ is non-spherical: Then there exist a $\lambda'' \in {\rm Par}_{J''}$, a corresponding $\mathbf{m}''$ and LS-monomials $\{\alpha_{j,k}''\}$ and $\{\beta_{j,k}''\}$ which are $I''$-dominant with the same weight. Let $\{J_{j,k}''\}$ and $\{K_{j,k}''\}$ be the LS-chains of $\{\alpha_{j,k}''\}$ and $\{\beta_{j,k}''\}$, respectively. For $1 \leq j \leq l$, define
\begin{equation}\label{eqn6.19}
m_j:=
\begin{cases}
m_j'' & 1 \leq j \leq a-1 \\
0 & a \leq j \leq r \\
m_{j-r+a-1}'' & r+1 \leq j \leq b \\
0 & b+1 \leq j \leq l
\end{cases}
.
\end{equation}

Define a function $\varphi_2 : [l] \rightarrow [l''] \cup \{0\}$ between index sets
\begin{equation}
\varphi_2 (j):=
\begin{cases}
j & 1 \leq j \leq a-1 \\
a & a \leq j \leq r \\
j-r-1+a & r+1 \leq j \leq b \\
0 & b+1 \leq j \leq l
\end{cases}
.
\end{equation}

For each $1 \leq j \leq l$ and $1 \leq k \leq m_j$, define
\begin{equation}
J_{j,k}:=
\begin{cases}
\varphi_2^{-1} (J_{j,k}'') & 1 \leq j \leq a-1 \\
\varphi_2^{-1} (J_{j-r+a-1,k}'') & r+1 \leq j \leq b
\end{cases}
.
\end{equation}
Then we can compute
\begin{equation}\label{eqn6.22}
\alpha_{j,k}=
\begin{cases}
\alpha_{j,k}'' & 1 \leq j \leq a-1 \\
s_{y,p+1}(\alpha_{j-r+a-1,k}'') & r+1 \leq j \leq b
\end{cases}
.
\end{equation}
Therefore, we have an $I$-dominant LS-monomial $\{\alpha_{j,k}\}$. Similarly, we construct $\{\beta_{j,k}\}$ from $\{\beta_{j,k}''\}$.

From (Case 1) and (Case 2), we constructed two distinct $I$-dominant LS-monomials $\{ \alpha_{j,k}\}$ and $\{ \beta_{j,k} \}$ with the same weight based on $u$ and $\lambda$. Let $w$ be the maximal representative of $W_I u$. Due to Theorem \ref{thm4.1}, $H^0(X(wW_J), L_\lambda)^*$ is not $L_I$-multiplicity-free. Due to Theorem \ref{lem2.4}, $X(wW_J)$ is not $L_I$-spherical. Therefore, $[u,I,J]$ is also non-spherical.

``$\Leftarrow$ of (b)'' : We prove the contraposition. Assume that $[u,I,J]$ is non-spherical. Due to Theorem \ref{thm4.1}, there exist a $\lambda \in \Parj$, a corresponding $\mathbf{m} = (m_1 , \cdots , m_l)$ and LS-monomials $\{\alpha_{j,k} \}$ and $\{\beta_{j,k}\}$ which are $I$-dominant with the same weight. Let $\{J_{j,k}\}$ and $\{K_{j,k}\}$ be the LS-chains of $\{\alpha_{j,k}\}$ and $\{\beta_{j,k}\}$, respectively. Since $\{\alpha_{j,k} \}$ and $\{\beta_{j,k}\}$ are distinct, for some $j_0$ and $k_0$, $\alpha_{j_0,k_0} \neq \beta_{j_0,k_0}$ \textit{i.e.} the $q$-th components of $\alpha_{j_0,k_0}$ and $\beta_{j_0,k_0}$ are different. Then we have two cases, $q \leq p+1$ and $q > p+1$. For each $\alpha = (\alpha_1 , \cdots , \alpha_n) \in \Z^n$, denote the size $|\alpha| := \alpha_1 + \cdots + \alpha_n$.

$\bullet$ (Case 1) $q \leq p+1$: Define $\pi_1 : \Z^n \rightarrow \Z^n$, $(\alpha_1 , \cdots , \alpha_n ) \mapsto (\alpha_1 , \cdots , \alpha_{p+1} , 0 , \cdots , 0)$. Then $|\pi_1(\alpha_{j,k})| = p$ if $j=r$. For each $r+1 \leq j \leq b$ and $1 \leq k \leq m_j$, the size $|\pi_1(\alpha_{j,k})|$ is either $p$ or $p+1$. In fact, our claim is that
\begin{equation}\label{eqn8.0}
|\pi_1(\alpha_{r,1})|, \cdots , |\pi_1(\alpha_{r, m_{r}})| , \cdots , |\pi_1(\alpha_{b,1})|, \cdots , |\pi_1(\alpha_{b,m_b})|
\end{equation}
are weakly increasing between $p$ to $p+1$. To show this, let $r +1 \leq j \leq b$ and $1 \leq k \leq m_j$. Then $|\pi_1(\alpha_{j,k})| = p$ if and only if $r+1 , \cdots , j \in J_{j,k}$. Since the sequence of $J_{j,k}$ is nested, this leads to our claim. Denote the number of $p$ which appears in \eqref{eqn8.0} by $\tilde{m}$. 

Choose $\lambda' \in {\rm Par}_{J'}$ from $\mathbf{m'}=(m_1' , \cdots, m_{l'}')$, which is defined
\begin{equation}
m_{j}' :=
\begin{cases}
m_{j+a-1} & 1 \leq j \leq r-a \\
\tilde{m} & j=r-a+1\\
m_{j+b-r+a-1} & r-a+2 \leq j \leq l'
\end{cases}
.
\end{equation}

Define a function $\psi_1 : [l'] \rightarrow [l]$ between index sets
\begin{equation}
\psi_1 (j):=
\begin{cases}
j+a-1 & 1 \leq j \leq r-a+1 \\
j-r+a-1+b & r-a+2 \leq j \leq l'
\end{cases}
.
\end{equation}

For each $1 \leq j \leq l'$ and $1 \leq k \leq m_j'$, define
\begin{equation}
J_{j,k}' :=
\begin{cases}
\psi_1^{-1} (J_{j+a-1,k}) & 1 \leq j \leq r-a \\
\psi_1^{-1} (J_{j+b-r+a-1,k}) & r-a+2 \leq j \leq l'
\end{cases}
.
\end{equation}
For $j = r-a+1$, consider $J_{r, 1} , \cdots , J_{r,m_r} , \cdots , J_{b,1} , \cdots , J_{b,m_b}$. Whenever $|\pi_1(\alpha_{j,k})| = p+1$, erase the corresponding $J_{j,k}$ from the list. Replace the remaining $J_{j,k}$ with $\psi_1^{-1} (J_{j,k})$. By construction, there should be $\tilde{m}$ index sets. Define $J_{r-a+1 , 1}' , J_{r-a+1 , 2}' , \cdots , J_{r-a+1,\tilde{m}}'$ as these index sets. Consequently, $\alpha_{r-a+1,1}' , \alpha_{r-a+1,2}' , \cdots , \alpha_{r-a+1,\tilde{m}}'$ are 
\begin{equation}
\pi_1 (\alpha_{r,1}) , \cdots , \pi_1 (\alpha_{r,m_r}) , \cdots , \pi_1 (\alpha_{b,1}) , \cdots , \pi_1 (\alpha_{b,m_b})
\end{equation}
with $\pi_1 (\alpha_{j,k})$ removed if $|\pi_1 (\alpha_{j,k})| = p+1$.

Note that for each $1 \leq j \leq a-1$, $|\pi_1 (\alpha_{j,k})| = p+1$. For each $a \leq j \leq r-1$ or $b+1 \leq j \leq l$, $|\pi_1 (\alpha_{j,k})| = i_j$. Since $\{ \alpha_{j,k} \}$ and $\{ \beta_{j,k} \}$ have the same weight, we have the same $\tilde{m}$ for $|\pi_1(\beta_{j,k})| = p$. By our claim, between $r \leq j \leq b$, the sequence of $|\pi_1(\beta_{j,k})|$ is also weakly increasing between $p$ to $p+1$. Therefore, $|\pi_1(\alpha_{j,k})| = |\pi_1(\beta_{j,k})|$ if $r \leq j \leq b$. When we construct an LS-monomial $\{ \beta_{j,k}' \}$ similarly from $\{ \beta_{j,k} \}$, the removals of $\pi_1 (\alpha_{j,k})$ and $\pi_1 (\beta_{j,k})$ are the same. This proves that $\{ \alpha_{j,k}' \}$ and $\{ \beta_{j,k}' \}$ have the same weight and are distinct.

$\bullet$ (Case 2) $q > p+1$: Define $\pi_2 : \Z^n \rightarrow \Z^n$, $(\alpha_1 , \cdots , \alpha_n ) \mapsto (0 , \cdots , 0 , \alpha_{p+2}  , \cdots , \alpha_n)$. Choose $\lambda '' \in {\rm Par}_{J''}$ from $\mathbf{m''} = (m_1 '' , \cdots , m_{l''}'')$, which is defined
\begin{equation}
m_{j}'':=
\begin{cases}
m_j & 1 \leq j \leq a-1 \\
m_{j+r-a+1} & a \leq j \leq l''
\end{cases}
.
\end{equation}

Define a function $\psi_2 : [l''] \rightarrow [l]$ between index sets
\begin{equation}
\psi_2 (j):= 
\begin{cases}
j & 1\leq j \leq a-1 \\
j+r-a+1 & a \leq j \leq l''
\end{cases}
.
\end{equation}

For $1 \leq j \leq l''$ and $1 \leq k \leq m_j ''$, define
\begin{equation}
J_{j,k}'' :=
\begin{cases}
\psi_2^{-1}(J_{j,k}) & 1 \leq j \leq a-1 \\
\psi_2^{-1} (J_{j+r-a+1,k}) & a \leq j \leq l''
\end{cases}
.
\end{equation}
Then if $1 \leq j \leq a-1$, $\alpha_{j,k}'' = \alpha_{j,k}$. If $a \leq j \leq l''$,
\begin{equation}
\alpha_{j,k}'' = 
\begin{cases}
\alpha_{j+r-a+1 , k} & |\alpha_{j+r-a+1,k}-\pi_2 (\alpha_{j+r-a+1,k})| = p+1 \\
\pi_2(\alpha_{j+r-a+1,k})+ \varpi_p  & |\alpha_{j+r-a+1,k}-\pi_2 (\alpha_{j+r-a+1,k})| = p
\end{cases}
.
\end{equation}

For each $r+1 \leq j \leq b$ and $1 \leq k \leq m_j$, count the number of $\alpha_{j,k}$ such that $|\alpha_{j,k}-\pi_2 (\alpha_{j,k})|$ is $p$. Denote this number by $\tilde{m}$. Note that for $1 \leq j \leq a-1$, $|\alpha_{j,k}-\pi_2 (\alpha_{j,k})| = p+1$. For $a \leq j \leq r$ or $b+1 \leq j \leq l$, $|\alpha_{j,k}-\pi_2 (\alpha_{j,k})| = i_j$. Since $\{ \alpha_{j,k} \}$ and $\{ \beta_{j,k} \}$ have the same weight, we have the same $\tilde{m}$ for the number of $|\beta_{j,k}-\pi_2 (\beta_{j,k})| = p$. Therefore, from $\{ \beta_{j,k} \}$, we may construct an LS-monomial $\{ \beta_{j,k} '' \}$, which has the same weight as $\{ \alpha_{j,k} ''\}$.

From (Case 1), we constructed a pair of $I'$-dominant LS-monomials $\{ \alpha_{j,k} ' \}$ and $\{ \beta_{j,k} ' \}$ with the same weight based on $u'$ and $\lambda'$. They are distinct due to the previous assumption $\alpha_{j_0 , k_0} \neq \beta_{j_0 , k_0}$. Let $w'$ be the maximal representative of $W_{I'} u'$. Due to Theorem \ref{thm4.1}, $H^0(X(w'W_{J'}), L_{\lambda'})^*$ is not $L_{I'}$-multiplicity-free. Due to Theorem \ref{lem2.4}, $X(w'W_{J'})$ is not $L_{I'}$-spherical. Therefore, $[u',I',J']$ is also non-spherical. Similarly, for (Case 2), $[u'',I'',J'']$ is non-spherical.
\end{proof}

\begin{remark}\label{rmk6.1}
At this point, we can reprove Theorem \ref{thm2.3}, not as a corollary to Theorem \ref{thm1.2}. Suppose $u$ has $321$ or $3412$ permutation patterns. From Theorem \ref{thm5.1} and \ref{thm5.2}, $X(w)$ is not $L_I$-spherical. Conversely, suppose $u$ does not have $321$ or $3412$ \textit{i.e.} $u$ is a standard Coxeter element. Due to Lemma \ref{lem6.2}, we may assume that $u = s_j$ for some $1 \leq j \leq n-1$. Here, we can easily check that $X(w)$ is $L_I$-spherical. 
\end{remark}

\begin{figure}
\centering
\begin{subfigure}[t]{0.45\textwidth}
\centering
\includegraphics[width=\textwidth]{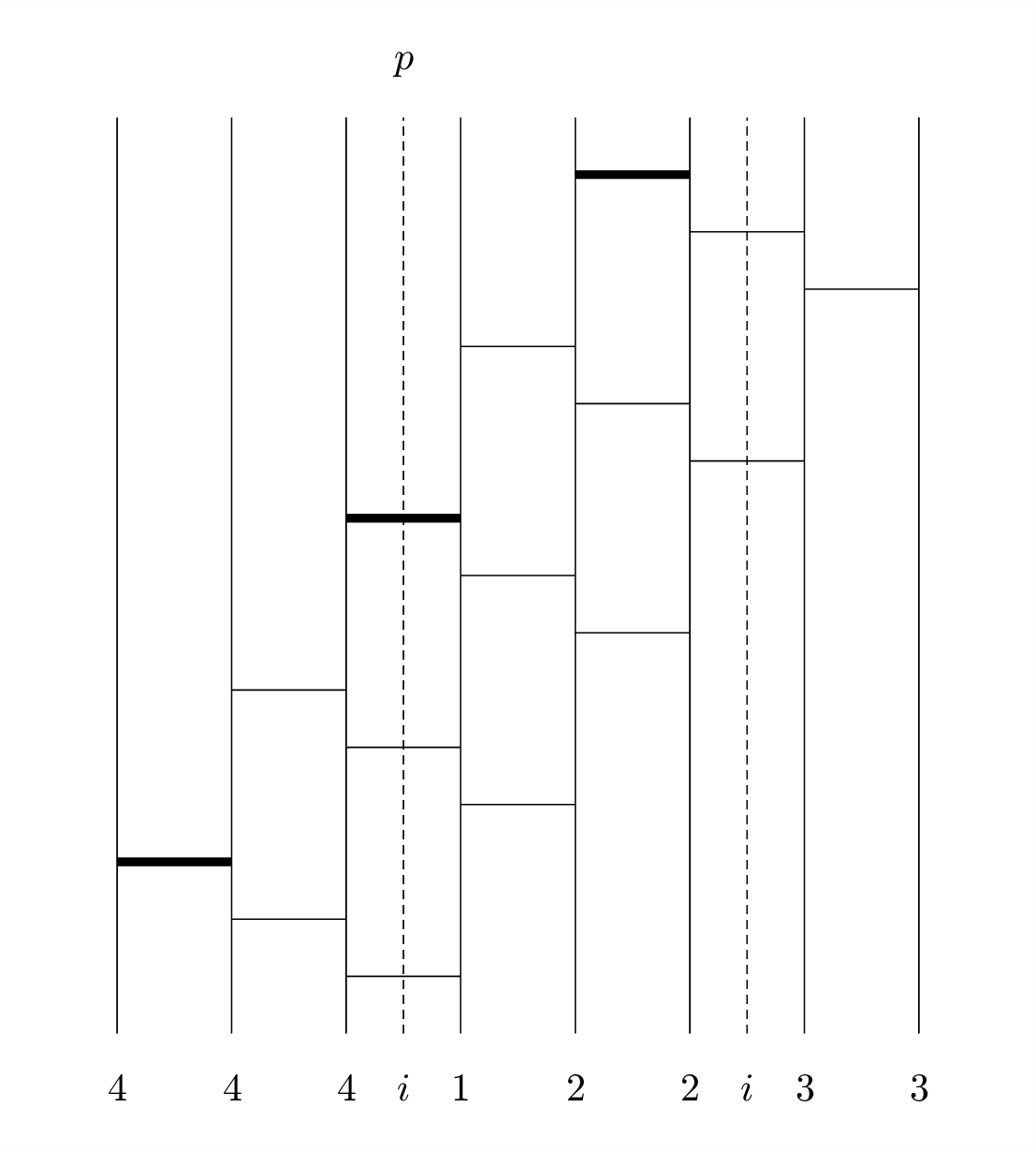}
\caption{Before Grassmannianization}
\label{fig20}
\end{subfigure}
\hfill
\begin{subfigure}[t]{0.45\textwidth}
\centering
\includegraphics[width=\textwidth]{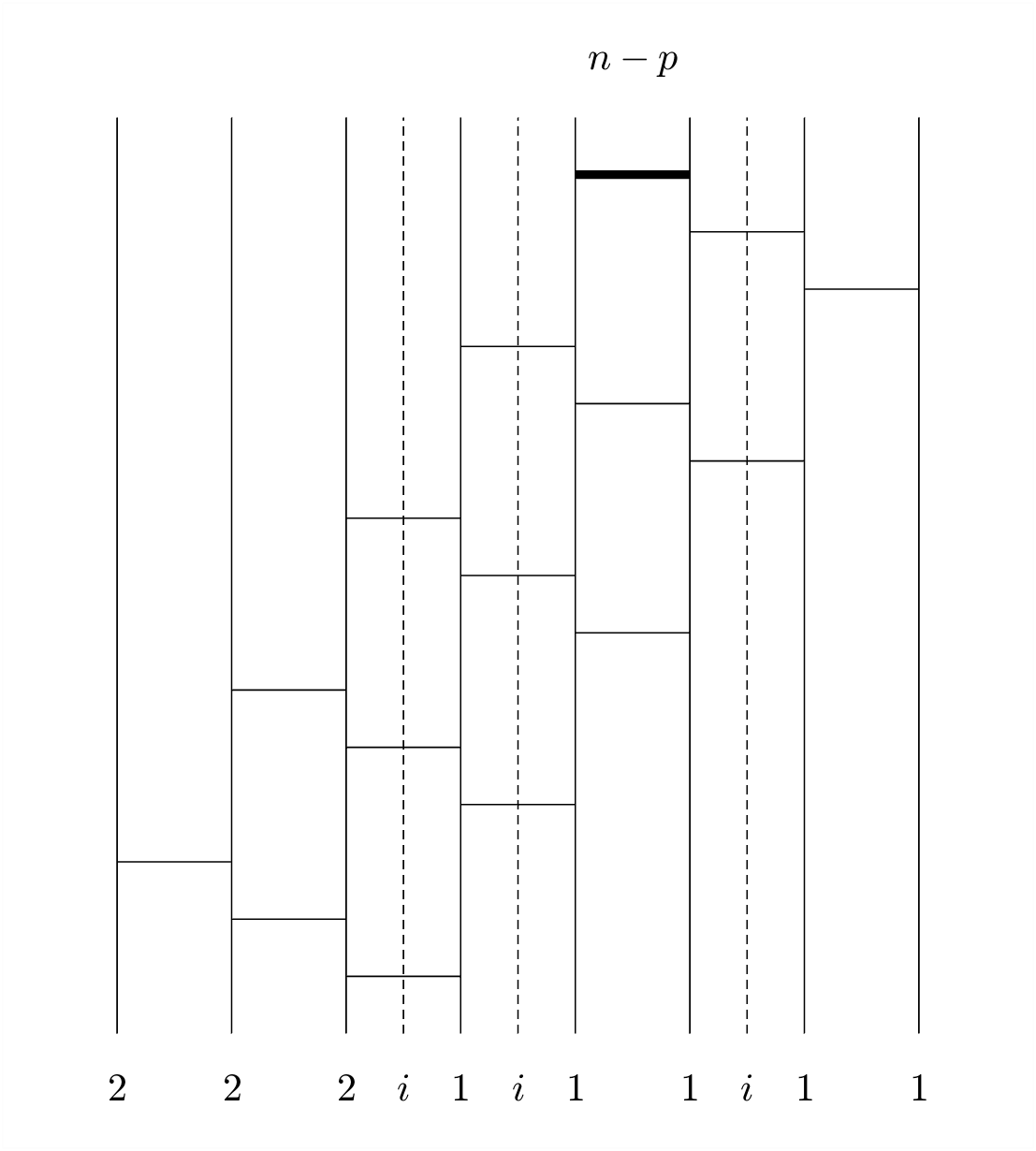}
\caption{After Grassmannianization}
\label{fig21}
\end{subfigure}
\caption{A left rectangle type}
\label{fig02021}
\end{figure}

Let $I, J \subseteq [n-1]$. Write $h:=n-|J|$. Let $u \in W$ and $[u,J] = x_1 x_2 \cdots x_n$. Suppose there exists a $2 \leq p \leq n-2$ such that 
\begin{enumerate}
\item $x_1 = x_2 = \cdots = x_{p} = h$.
\item $1 = x_{p+1} \leq x_{p+2} \leq \cdots \leq x_{n} = h-1$.
\item $1 , 2 , \cdots , p-1 \in I$.
\end{enumerate}
See Figure \ref{fig20} for an example. Automatically, $x_p > x_{p+1}$ leads to $p \notin I$ due to Lemma \ref{lem2.6}. In this case, $[u,I,J]$ is called a \textbf{left rectangle type}. Indeed, if $u$ is the minimal representative of $uW_J$, then
\begin{equation}\label{eqn6.25}
u = s_{p}s_{p-1} \cdots s_2 s_1 s_{p+1} s_{p} \cdots s_3 s_2 \cdots s_{n-1} s_{n-2} \cdots s_{n-p+1}s_{n-p},
\end{equation}
which looks like a rectangle in the wiring diagram.

From $p+1 \leq j \leq n-1$, whenever $x_{j} < x_{j+1}$, remove $j$ from the index set $I$. Let $I'$ be the set after subtraction. Let $J' = [n-1]-\{n-p\}$. In other words, from $[u,I,J] = x_1 \cdots x_n$, we insert $i$ between $x_{j}$ and $x_{j+1}$ whenever $x_j < x_{j+1}$ and then replace $x_1 , \cdots , x_p$ with $2$ and $x_{p+1} , \cdots , x_{n}$ with $1$. For instance, if $[u,I,J] = 444i122i33$ as in Figure \ref{fig20}, we have $[u,I',J'] = 222i1i11i11$ as in Figure \ref{fig21}. $[u,I',J']$ is called \textbf{Grassmannianization} of $[u,I,J]$.

Similarly, if there exists a $2 \leq p \leq n-2$ with the following properties, $[u,I,J]$ is called a \textbf{right rectangle type}.
\begin{enumerate}
\item $x_{p+1} = x_{p+2} = \cdots = x_n =1$.
\item $2 = x_1 \leq x_2 \leq \cdots x_{p} = h$.
\item $p+1 , p+2 , \cdots , n-1 \in I$.
\end{enumerate}
Again, from $1 \leq j \leq p-1$, whenever $x_j < x_{j+1}$, subtract $j$ from $I$ to construct $I'$. Let $J':= [n-1] - \{n-p\}$. $[u,I',J']$ is also called Grassmannianization of $[u,I,J]$.

\begin{lemma}\label{lem8.2}
Let $w \in S_n$, ${\sf supp}(w) \subseteq [a, b]$ and $c \in [a, b]$. If $b<n$, let $u = s_{b} s_{b-1} \cdots s_c$. Then $\ell (wu) = \ell (w) + \ell (u)$. Similarly, if $a >1$, let $v = s_{a-1} s_{a} \cdots s_{c-1}$. Then $\ell (wv) = \ell (w) + \ell (v)$.
\end{lemma}

\begin{proof}
Let $(i_1 , \cdots , i_l)$ and $(i_1 ' , \cdots , i_{l'}')$ be reduced words of $w$ and $u$, respectively. Suppose $\ell (wu) \neq \ell (w) + \ell (u)$. Then $(i_1 , \cdots , i_l , i_1' , \cdots , i_{l'}')$ is not reduced. There exists $p \in [l]$ and $p' \in [l']$ such that
\begin{equation}
wu = s_{i_1} \cdots s_{i_{p-1}} s_{i_{p+1}} \cdots s_{i_{l}} s_{i_1'} \cdots s_{i_{p'-1}'} s_{i_{p'+1}'} \cdots s_{i_{l'}'}.
\end{equation}
This contradicts $(wu)^{-1} (b+1) = c$. Hence, $\ell (wu) = \ell (w) + \ell (u)$. $\ell (wv) = \ell (w) + \ell (v)$ can be treated similarly.
\end{proof}

\begin{lemma}\label{lem6.3}
Let $G= GL_n$, $W=S_n$, $I,J \subseteq [n-1]$ and $u \in {}^IW^J$. Suppose $[u,I,J]$ is either of left rectangle type or of right rectangle type and let $[u,I',J']$ be the Grassmannianization of $[u,I,J]$. Then $[u,I,J]$ is spherical if and only if $[u,I',J']$ is spherical.
\end{lemma}
\begin{proof}
Due to symmetry, we may assume that $[u,I,J]$ is of left rectangle type with $[u,J] = x_1 \cdots x_n$ satisfying $x_1 = \cdots = x_p = h$ and $1 = x_{p+1 } \leq \cdots \leq x_{n} = h-1$. Let $[u,I',J']$ be the Grassmannianization of $[u,I,J]$. Let $w$ and $w'$ be the maximal representatives of $W_I u$ and $W_{I'} u$, respectively.

``$\Rightarrow$'': We prove the contraposition. Suppose $X(w'W_{J'})$ is not $L_{I'}$-spherical. Due to Theorem \ref{lem2.4}, there exists a $\lambda' \in {\rm Par}_{J'}$ such that $H^0 (X(w'W_{J'}), L_{\lambda'})$ is not $L_{I'}$-multiplicity-free. Since determinant twists $(1, \cdots , 1)$ preserve multiplicities, we may assume that $\lambda_n' = 0$. Let $N \in \Z_{> 0}$ be chosen so that $\lambda ' = N \varpi_{n-p}$. Using the reduced word expression of $u$ as in \eqref{eqn6.25} and $\lambda'$, we have $\mathbf{i'}=(i_1' , \cdots , i_{l'}')$ and $\mathbf{m'} = (m_1' , \cdots , m_{l'}')$ where $m_{l'}' = N$ and $m_1' = \cdots = m_{l'-1}' = 0$.

Due to Theorem \ref{thm4.1}, there exist distinct $I'$-dominant LS-monomials $\{\alpha_{j,k}'\}$ and $\{\beta_{j,k}'\}$ with the same weight corresponding to $(\mathbf{i'},\mathbf{m'})$. Let $\{J_{j,k}'\}$ and $\{K_{j,k}'\}$ be the corresponding LS-chains. Choose $\lambda \in \Parj$
\begin{equation}
\lambda:= \sum_{j \in [n-1]-J} N\varpi_j.
\end{equation}
Again, from \eqref{eqn6.25}, we have $\mathbf{i} = \mathbf{i'}$ and $l=l'$. For $\mathbf{m} = (m_1 , \cdots , m_l)$, we have $m_l = N$ and $m_j$ is either $0$ or $N$. Define
\begin{equation}
J_{j,k}:=
\begin{cases}
[l] & 1 \leq j \leq l-1 \\
J_{l',k}' & j = l
\end{cases}
, \quad K_{j,k} :=
\begin{cases}
[l] & 1 \leq j \leq l-1 \\
K_{l',k}' & j = l
\end{cases}
.
\end{equation}
Define $\{ \alpha_{j,k} \}$ and $\{\beta_{j,k}\}$ from $J_{j,k}$ and $K_{j,k}$. Due to the construction, they are distinct but with the same weight. Also, they are $I'$-dominant. Consider $q \in I-I'$. Then $x_{q} < x_{q+1}$ and $q-p \in [n-1]-J$. Due to the construction, there are $N$ copies of $\alpha_{j,k} = \varpi_{q} - \varpi_p$. This guarantees that $\{\alpha_{j,k}\}$ is $I$-dominant. So is $\{\beta_{j,k}\}$. Using Theorem \ref{thm4.1}, $H^0(X(wW_J), L_\lambda)$ is not $L_I$-multiplicity-free. Due to Theorem \ref{lem2.4}, $X(wW_J)$ is not $L_I$-spherical.

``$\Leftarrow$'': We prove the contraposition. Suppose $X(wW_J)$ is not $L_I$-spherical. Due to Theorem \ref{lem2.4}, there exists a $\lambda \in \Parj$ such that $H^0(X(wW_J),L_\lambda)$ is not $L_{I}$-multiplicity-free. Choose $\mathbf{i} = (i_1 , \cdots , i_l)$ from \eqref{eqn6.25} and $\mathbf{m} = (m_1 , \cdots, m_l)$ from $\lambda$. Due to Theorem \ref{thm4.1}, there exist distinct $I$-dominant LS-monomials $\{\alpha_{j,k}\}$ and $\{\beta_{j,k}\}$ with corresponding LS-chains $\{J_{j,k}\}$ and $\{K_{j,k}\}$. Define $M:= \sum_{1 \leq j \leq l} {m_j}$. Let
\begin{equation}
\lambda' := M \varpi_{n-p}.
\end{equation}
Choose $\mathbf{i'} = \mathbf{i}$ where $l = l'$ and $\mathbf{m'} = (m_1 ' , \cdots , m_{l'}')$ from $\lambda'$. Then $m_{j}' = 0$ except for $m_{l'}' = M$.

From $J_{j,k}$ construct $J_{j,k} \cup [j+1 , l']$. We then have a nested sequence
\begin{equation}
J_{1,1} \cup [2,l'] , \cdots , J_{1,m_1} \cup [2, l'] , J_{2,1} \cup [3,l'] , \cdots, J_{2, m_2} \cup [3,l'] , \cdots 
\end{equation}
Define $J_{l',1}' , J_{l',2}' , \cdots , J_{l',M}'$ by reading off the index sets above. Our claim is that this gives us $I'$-dominant LS-monomial $\{\alpha_{j,k}'\}$. To show this, let $1 \leq j \leq l$ such that $m_j \neq 0$. $\mathbf{i} (J_{j,k} \cap [j])$ is a reduced word. Due to Lemma \ref{lem8.2}, $\mathbf{i} (J_{j,k} \cup [j+1, l '])$ is also a reduced word. 

From \eqref{eqn6.25}, there exists $q \in \Z_{>0}$ such that $j = pq$. If $p+q <n$, then
\begin{align}\label{eqn8.9}
u(\mathbf{i} (J_{j,k} \cup [j+1, l '])) \cdot \varpi_{n-p} &= u(\mathbf{i} (J_{j,k} \cap [j])) \cdot (\varpi_k + \epsilon_{p+q+1} + \epsilon_{p+q+2} + \cdots + \epsilon_n) \\
&= \alpha_{j,k} + \epsilon_{p+q+1} + \epsilon_{p+q+2} + \cdots + \epsilon_n.
\end{align}
$x_{p+q} < x_{p+q+1}$ follows from $m_j \neq 0$. This implies that $p+q \notin I'$. Therefore, $\{ \alpha_{j , k}' \}$ is still $I'$-dominant. If $p+q = n$, then $j = l'$. Since $J_{j,k} \cup [j+1, l '] = J_{j,k}$, \eqref{eqn8.9} is $\alpha_{j,k}$.

Similarly, construct $K_{j,k} \cup [j+1 , l']$ and $\beta_{j,k}'$. Then again, $\beta_{j,k}$ is replaced with $\beta_{j,k} + \epsilon_{p+q+1} + \cdots + \epsilon_n$. This shows that $\{ \alpha_{j,k}' \}$ and $\{ \beta_{j,k}' \}$ have the same weight and are distinct. Using Theorem \ref{thm4.1}, $H^0(X(w'W_{J'}), L_{\lambda'})$ is not $L_{I'}$-multiplicity-free. Due to Theorem \ref{lem2.4}, $X(w'W_{J'})$ is not $L_{I'}$-spherical.
\end{proof}

\begin{theorem}\label{thm6.4}
Let $G = GL_n$, $W=S_n$, $I,J \subseteq [n-1]$ and $u \in {}^I W^J$. Suppose $u$ is bigrassmannian. Then $[u,I,J]$ is spherical if and only if $[u,I,J]$ avoids the patterns in \eqref{eqn1.7}.
\end{theorem}
\begin{proof}
If $u$ is a standard Coxeter element, apply Lemma \ref{lem3.3} and \ref{lem6.2}. Then we may assume that $u = s_i$ for some $1 \leq i \leq n-1$ or $u$ is the identity. Use Theorem \ref{thm4.1} to prove that $[u,I,J]$ is spherical. Also, $[u,I,J]$ does not have patterns in \eqref{eqn1.7} since we have only two cases.
\begin{enumerate}
\item $[u,I,J]$ is weakly increasing.
\item $[u,I,J]$ is weakly increasing except for two consecutive entries.
\end{enumerate}
This proves that $[u,I,J]$ is always spherical and avoids patterns in \eqref{eqn1.7}.

Assume that $u$ is not a standard Coxeter element. Since $u$ is a bigrassmannian, there exist $1 \leq q < p < r \leq n-1$ such that
\begin{equation}
u = s_p s_{p-1} \cdots s_{q+1} s_q s_{p+1} s_p \cdots s_{q+2} s_{q+1} \cdots s_r s_{r-1} \cdots s_{r-p+q +1} s_{r-p+q}.
\end{equation}

Our claim is that it is enough to prove the case $q=1$ and $r= n-1$. To show this, construct a permutation $\bar{u} \in S_{r-q+2}$ by replacing $s_i$ with $s_{i-q+1}$
\begin{equation}
\bar{u} := s_{p-q+1} s_{p-q} \cdots s_2 s_1 s_{p-q+2} s_{p-q+1} \cdots s_3 s_2 \cdots s_{r-q+1} s_{r-q} \cdots s_{r-p+2} s_{r-p+1}.
\end{equation}
Also, define
\begin{equation}
\bar{I} := ((1-q) + I) \cap [r-q+1], \quad \bar{J} := ((1-q) + J) \cap [r-q+1].
\end{equation}
Then we may see that $[u,I,J]$ avoids patterns \eqref{eqn1.7} if and only if $[\bar{u}, \bar{I}, \bar{J}]$ avoids patterns \eqref{eqn1.7}. Also, since $s_1 , \cdots s_{q-1}$ and $s_{r+1} , \cdots s_{n-1}$ do not appear in the reduced word expression of $u$, they play no role in the construction of LS-monomials. Applying Theorem \ref{thm4.1} and Theorem \ref{lem2.4}, $[u,I,J]$ is spherical if and only if $[\bar{u} , \bar{I}, \bar{J}]$ is spherical. Therefore, we may assume that $q=1$ and $r = n-1$ without loss of generality.

``$\Rightarrow$'': We prove the contraposition. Suppose $[u,I,J]$ has one of the patterns in \eqref{eqn1.7}. Due to Theorem \ref{thm5.1}, \ref{thm5.3} and \ref{thm5.4}, we may exclude the patterns $3412$, $3i312$, $231i1$ and $2i21i1$. Then the remaining cases are $[u,I,J]$ being left rectangle or right rectangle types. Without loss of generality, assume that $[u,I,J]$ is of left rectangle type. Then the possible patterns that may appear are $3331122$, $22211i11$, $22ii1i1$, $33ii12$, $44123$ and $3312i2$. After Grassmannization, they are converted into $22211i11$, $22211i11$, $22ii1i1$, $22ii1i1$, $221i1i1$ and $221i1i1$ in $[u,I',J']$. Due to Theorem \ref{thm2.2}, $[u,I',J']$ is non-spherical. Due to Lemma \ref{lem6.3}, this implies that $[u,I,J]$ is non-spherical.

``$\Leftarrow$'': Suppose $[u,I,J]$ is non-spherical. If $[u,I,J]$ has one of the patterns $3412$, $3i312$, $231i1$ and $2i21i1$, then we are done. Suppose $[u,I,J]$ avoids them all. Then $[u,I,J]$ is of left rectangle or right rectangle type. Without loss of generality, assume that $[u,I,J]$ is of left rectangle type. Let $[u,I',J']$ be the Grassmannianization of $[u,I,J]$. Due to Lemma \ref{lem6.3}, $[u,I',J']$ is non-spherical. From Theorem \ref{thm2.2}, $[u,I',J']$ has one of the patterns in \eqref{eqn2.1}. Since we are assuming that $[u,I,J]$ is of left rectangle type, either $22ii1i1$ or $22211i11$ appears in $[u,I',J']$.

$\bullet$ (Case 1) $[u,I',J']$ has a $22ii1i1$ pattern: The second $i$ or the third $i$ might have been inserted during Grassmannianization. If the second $i$ was inserted, then $[u,I,J]$ has $44123$ or $3312i2$. If the third $i$ was inserted, then $[u,I,J]$ has $44123$ or $33ii12$. If both of them were inserted, $[u,I,J]$ has $44123$. If none of them were inserted, then $[u,I,J]$ has $22ii1i1$ or $33ii12$.

$\bullet$ (Case 2) $[u,I',J']$ has a $22211i11$ pattern: If $i$ was inserted during Grassmannianization, then $[u,I,J]$ has a $3331122$ or $44123$ pattern. If $i$ was not inserted, then $[u,I,J]$ has a $22211i11$, $3312i2$, $33ii12$, $3331122$ or $44123$ pattern.

From (Case 1) to (Case 2), $[u,I,J]$ has one of the patterns in \eqref{eqn1.7}.
\end{proof}

\begin{proof}[Proof of Theorem \ref{thm1.2}]
If $[u,I,J]$ has one of the patterns in \eqref{eqn7.3}, then $X(wW_J)$ is not $L_I$-spherical due to Theorem \ref{thm5.1}, \ref{thm5.2}, \ref{thm5.3}, \ref{thm5.4}, \ref{thm6.1} and \ref{thm6.5}. Therefore, we only need to determine sphericality for the case when $[u,I,J]$ avoids the patterns in \eqref{eqn7.3}. From Lemma \ref{lem7.1}, $u$ avoids permutation patterns in \eqref{eqn7.4}. Using Lemma \ref{lem8.1}, we may decompose $u$ into a product of bigrassmannian permutations. Due to Lemma \ref{lem3.3} and \ref{lem6.2}, we may assume that $u$ is bigrassmannian without loss of generality. From Theorem \ref{thm6.4}, $X(wW_J)$ is $L_I$-spherical if and only if $[u,I,J]$ avoids the patterns in \eqref{eqn1.7}.
\end{proof}

To prove Corollary \ref{coro1.3}, we need to show that for each $u \in {}^I W^J$, if $[u,I,J]$ is spherical, then so is $[u^{-1}, J, I]$. We explain here how to convert $[u,I,J]$ into $[u^{-1},J,I]$. Read the numerical entries in $[u,I,J]$ in weakly increasing order. If there are entries with the same value, read them from the left to the right. As we read the number, record it into another number $k$, so that the number was between the $(k-1)$th and the $k$th $i$. Whenever the number we are reading increases, insert $i$.

For instance, consider $[u,I,J] = 13i4i223i4$. From $\underline{1}3i4i223i4$, read off $1$ which is between the $1$st and the $0$th (which actually does not exist) $i$. Therefore, record it into $1$. Before reading off $2$, record $i$ since the number increases. From $13i4i\underline{22}3i4$, read off two $2$s between the second and the third $i$. Record them into two $3$s. In this way, we have $1i33i13i24$, which is $[u^{-1}, J, I]$.

\begin{proof}[Proof of Corollary \ref{coro1.3}]
Assume that $[u,I,J]$ is not spherical. Due to Theorem \ref{thm1.2}, $[u,I,J]$ has one of the patterns in \eqref{eqn1.7}. Then we proceed as follows.

\begin{itemize}
\item[\textit{(Step 1)}] Assume that $[u,I,J]$ has a pattern $P_1$ in \eqref{eqn1.7}.
\item[\textit{(Step 2)}] Using Lemma \ref{lem2.6}, insert $i$ to $P_1$ whenever there is a descent and construct $P_2$.
\item[\textit{(Step 3)}] Use $I$-expansion to eliminate $i$ in $[u,I,J]$ except for ones included in $P_2$ and construct $[u',I',J]$. Then $[u',I',J]$ still has $P_2$.
\item[\textit{(Step 4)}] From $P_2$, $[u'^{-1} , J , I']$ has a pattern $P_3$ in \eqref{eqn1.7}, which means it is non-spherical.
\item[\textit{(Step 5)}] Due to Lemma \ref{lem4.5}, $[u^{-1}, J , I]$ is non-spherical since its $J$-expansion  $[u'^{-1} , J , I']$ is non-spherical.
\end{itemize}

\begin{table}[h!]
\centering
\begin{tabular}{||c c c || c c c || c c c ||} 
 \hline
 $P_1$ & $P_2$ & $P_3$ & $P_1$ & $P_2$ & $P_3$ & $P_1$ & $P_2$ & $P_3$ \\ [0.5ex] 
 \hline\hline
$321$ & $3i2i1$ & $321$ & {}  $221211$ & $22i12i11$ & $233112$ & {} $331322$ & $33i13i22$ & $233112$ \\
$3412$ & $34i12$ & $2i21i1$ & {} $244113$ & $244i113$ & $221211$ {} & $22211i11$ & $222i11i11$ & $2233111$\\
$231i1$ & $23i1i1$ & $231i1$ & {} $233112$ & $233i112$ & $221211$ & {} $2233111$ & $2233i111$ & $22211i11$ {} \\
$3i312$ & $3i3i12$ & $3i312$ & {} $3312i2$ & $33i12i2$ & $2i2311$ {} &  $3331122$ & $333i1122$ & $22i22111$\\
$2i21i1$ & $2i2i1i1$ & $3412$ & {}  $44123$ & $44i123$ & $2i2ii11$ & {} $22i22111$ & $22i22i111$ & $3331122$  {} \\
$23411$ & $234i11$ & $22ii1i1$ {} & $221311$ & $22i13i11$ & $233112$ & {} $331422$ & $33i14i22$ & $233112$ {} \\
$2i2311$ & $2i23i11$ & $3312i2$ & {} $\underline{22ii1i1}$ & $22i1i1i1$ & $23411$ & {} $\underline{2i2ii11}$ & $2i2i2i11$ & $44123$\\
$\underline{33ii12}$ & $33i1i12$ & $23ii11$ & {} $\underline{23ii11}$ & $23i3i11$ & $33ii12$ & {} & & {}\\ [1ex]
 \hline
\end{tabular}
\caption{The patterns $P_1$, $P_2$ and $P_3$ from Step 1 to Step 5}
\label{tab1}
\end{table}
There are exceptions when $P_1 = 22ii1i1, 33ii12, 2i2ii11 , 23ii11$; see underlined patterns in Table \ref{tab1}. We may take $P_2 = 22i1i1i1 , 33i1i12 , 2i2i2i11 , 23i3i11$. Otherwise, we take other patterns as $P_1$ and start again. For instance, if $[u,I,J]$ has $22ii1i1$, then it has either $33i12i2$, $22i1i1i1$, $321$, $2i21i1$ or $231i1$. Therefore, we may assume that $[u,I,J]$ has $22i1i1i1$.
\end{proof}

\begin{proof}[Proof of Corollary \ref{coro1.4}]
Suppose for some $I,J \subseteq [n-1]$, $u \in {}^I W^J$ and $[u,I,J]$ is spherical. Due to Theorem \ref{thm1.2}, $[u,I,J]$ avoids patterns in \eqref{eqn1.7}. From Lemma \ref{lem7.1}, $u$ avoids permutation patterns in \eqref{eqn7.4}. Suppose $u$ has a permutation pattern $2561374$. Then $u^{-1}$ has a permutation pattern $4157236$. When computing $[u,I,J]$, $i$ should be inserted when it is decreasing. In other words, $4i157i236$. Note that $3$ and $4$ cannot be equalized since the ordering is reversed. Since $[u,I,J]$ avoids $3412$ or $3i312$ patterns, $2$ and $3$ in $4i157i236$ should be equalized. Also, due to $6$, $5$ and $7$ cannot be equalized. Therefore, from underlined entries $\underline{4i}1\underline{57}i\underline{23}6$, $[u,I,J]$ has either $2i2311$ or $23411$ patterns. From Theorem \ref{thm1.2}, this implies that $[u,I,J]$ is non-spherical, leading to a contradiction. Hence, $u$ avoids a permutation pattern $2561374$. $4157236$ can be treated in a similar fashion.

Conversely, suppose $u$ avoids permutation patterns in \eqref{eqn1.12}. From Lemma \ref{lem4.11}, we have bigrassmannian permutations $w^{(1)} , \cdots , w^{(k)}$ with $u = w^{(1)} \cdots w^{(k)}$. We want to construct $I , J \subseteq [n-1]$. We first start with $I = J= [n-1]$. For each $1 \leq i \leq k$, we do the following.
\begin{itemize}
\item Remove ${\sf maxs}(w^{(i)}) - {\sf se}(w^{(i)})$ from $I$. Remove ${\sf maxs}(w^{(i)}) - {\sf sw}(w^{(i)})$ from $J$.
\item Suppose there exists $j \neq i$ such that ${\sf maxs}(w^{(j)}) = {\sf mins}(w^{(i)})$.
\begin{itemize}
\item If $j <i$, then remove ${\sf mins}(w^{(i)})$ from $I$.
\item If $j>i$, then remove ${\sf mins}(w^{(i)})$ from $J$.
\end{itemize}
\item Suppose there exists $j \neq i$ such that ${\sf mins}(w^{(j)}) = {\sf maxs}(w^{(i)})$.
\begin{itemize}
\item If $j <i$, then remove ${\sf maxs}(w^{(i)})-1$ from $I$.
\item If $j>i$, then remove ${\sf maxs}(w^{(i)})-1$ from $J$.
\end{itemize}
\end{itemize}
As a result, $u \in {}^I W^J$. Moreover, when applying Lemma \ref{lem3.3} and \ref{lem6.2} to decompose $u$ into $w^{(1)} , \cdots , w^{(k)}$, $I$ and $J$ do not change. 

Our claim is that $[u,I,J]$ is spherical. Applying Lemma \ref{lem3.3} and \ref{lem6.2}, we need to show that for each $1 \leq i \leq k$, $[w^{(i)}, I, J]$ is spherical. Write ${\sf sw}(w^{(i)}) = p$ and ${\sf se}(w^{(i)}) = q$. If $p = 1$ or $q = 1$, then it is spherical since $w^{(i)}$ is standard Coxeter. Suppose $p, q \geq 2$. Let $[w^{(i)} , J] = x_1 \cdots x_n$. Recall that $[w^{(i)}, I , J]$ is $x_1 \cdots x_n$ with $i$ insertions. If $[w^{(i)}, I , J]$ is non-spherical, then it has a pattern between $x_a$ and $x_b$ where $a = {\sf mins}(w^{(i)})$ and $b = {\sf maxs}(w^{(i)})$. However, by taking the pattern between $x_a$ and $x_b$ in $[w^{(i)}, I , J]$, it should be one of
\begin{subequations}
\begin{equation}
\overbrace{2\cdots 2}^{\# p}i\overbrace{1 \cdots 1}^{\# q}, \quad \overbrace{2\cdots 2}^{\# p-1}3i\overbrace{1 \cdots 1}^{\# q}, \quad 2i\overbrace{2\cdots 2}^{\# p-1}i\overbrace{1 \cdots 1}^{\# q} ,\quad 2i3i \overbrace{1 \cdots 1}^{\# q},
\end{equation}
\begin{equation}
\overbrace{3\cdots 3}^{\# p}i1\overbrace{2 \cdots 2}^{\# q-1}, \quad \overbrace{2\cdots 2}^{\# p}i\overbrace{1 \cdots 1}^{\# q-1}i1, \quad \overbrace{3\cdots 3}^{\# p}i1i2.
\end{equation}
\end{subequations}
From Theorem \ref{thm1.2} and Lemma \ref{lem4.11}, these are all spherical. Therefore, $[u,I,J]$ is spherical.

Let $I_{max}$ and $J_{max}$ be the left ascent set and the right ascent set of $u$. Then $I_{max} \supseteq I$ and $J_{max} \supseteq J$. From $I$-expansion and $J$-expansion Lemma \ref{lem4.4} and \ref{lem4.5}, $[u,I_{max}, J_{max}]$ is also spherical.
\end{proof}

\section*{Acknowledgements and AI disclosure}

We thank Yibo Gao, Reuven Hodges and Alexander Yong for helpful remarks on drafts of this preprint. Earlier, we reported Theorem \ref{thm1.2} and Theorem \ref{thm1.1} in the proceedings of the 38th International Conference on ``Formal Power Series and Algebraic Combinatorics'' without using generative AI tools. For this preprint, we used ChatGPT Plus in the final stage to help check our arguments and proofread. ChatGPT Plus identified a minor error in the proof of Corollary \ref{coro1.4}, which we corrected ourselves. The author takes full responsibility for all content. This work was partially supported by UIUC Campus Research Board grant RB24025, NSF RTG in Combinatorics DMS 1937241 and a Simons Foundation Dissertation Fellowship.

\end{document}